\documentclass[10pt]{article}
\usepackage{amsmath}
\usepackage{hyperref}
\usepackage{amsmath}
\usepackage{amsfonts}
\usepackage{graphicx}
\usepackage{amssymb}
\usepackage{leftidx}
\usepackage{color}
\usepackage{booktabs}
\DeclareRobustCommand{\func}[1]{\operatorname{#1}}
\allowdisplaybreaks
\newtheorem{theorem}{Theorem}[section]
\newtheorem{lemma}[theorem]{Lemma}
\newtheorem{proposition}[theorem]{Proposition}

\newtheorem{example}[theorem]{Example}
\newtheorem{definition}[theorem]{Definition}
\newtheorem{corollary}[theorem]{Corollary}
\newtheorem{remark}[theorem]{Remark}

\newcommand{\dif}{\mathrm{d}}
\newcommand{\Proj}{\operatorname{Proj}}
\newcommand{\esssup}{\operatorname{ess\,sup}}
\newcommand{\E}{\mathbb E}

\newcommand{\R}{\mathbb R}

\newcommand{\vecop}{\operatorname{vec}}
\newcommand{\eps}{\varepsilon}
\newcommand{\rhobar}{\bar m}

\newcommand{\Var}{\operatorname{Var}}
\newcommand{\Htwo}{\mathbb H^2}
\newcommand{\Stwo}{\mathbb S^2}

\newcommand{\norm}[1]{\left\lVert #1\right\rVert}

\newcommand{\trans}{\mathsf T}

\newcommand{\Pz}{P_{0}}
\numberwithin{equation}{section}

\newenvironment{proof}[1][Proof]{\noindent\textbf{#1.} }{\ \rule{0.5em}{0.5em}}

\begin{document}

\title{ A Regularized BSDE Method for Constrained Mean-Variance Portfolio Selection}
\author{Liangquan Zhang$^{1}\thanks{
The corresponding author: L. Zhang acknowledges the financial support partly by the National Nature
Science Foundation of China (Grant No. 12571486 and 12171053) and the
Fundamental Research Funds for the Central Universities, and the Research
Funds of Renmin University of China (No. 23XNKJ05). E-mail: xiaoquan51011@163.com.}, \ \ $Xun Li$^{2}$\thanks{X. Li acknowledges the financial support by the Hong Kong General Research Fund under grant 15226922. E-mail:
li.xun@polyu.edu.hk.} \\
{\small  1. School of Mathematics}\\
{\small \ Renmin University of China, Beijing 100872, China}\\
{\small 2. Department of Applied Mathematics }\\
{\small \  The Hong Kong Polytechnic University, Hong Kong, China }}
\maketitle

\begin{abstract}

This paper develops a regularized backward stochastic differential
equation (BSDE) framework for continuous-time mean-variance portfolio selection with random coefficients under convex portfolio constraints that need not be bounded. Unlike the classical formulation, the volatility is not required to be non-degenerate. Terminal wealth has the affine form \(X_T^u=X_T^0+\mathcal L u\), with \(\mathcal L\) a bounded linear operator. With redundant assets, \(\mathcal L\) has a nontrivial null space: controls differing by a kernel element yield the same wealth and objective, so the solution set is a closed convex set rather than a single point and the optimal portfolio is not unique. Using the Lions derivative, we express the cost via an adjoint mean-field BSDE and formulate optimality as a forward-backward stochastic differential variational inequality (FBSDVI), characterized by projection onto the feasible set (componentwise clipping for box constraints). A Tikhonov \(L^2\) penalty renders each regularized problem strictly convex and uniquely solvable. As \(\varepsilon\downarrow0\), the controls converge strongly in \(\mathbb H^2\) to the minimum-norm optimizer while the original objective value is recovered.

Two redundant assets illustrate the mechanism: wealth depends only on the aggregate \(v=u_1+u_2\), while the long-short \(w=u_1-u_2\) is the null space. The adjoint BSDE fixes $v$, and regularization forces \(w=0\), selecting the equal-split minimum-energy portfolio. The infinite-dimensional HJB reduces to finite-dimensional Riccati equations, recovering the BSDE feedback as \(\varepsilon\downarrow0\); a numerical example confirms the convergence.

\end{abstract}
\newpage
\tableofcontents

\noindent \textbf{AMS subject classifications:} 93E20, 60H15, 60H30.

\noindent\textbf{Keywords:} Mean-variance portfolio selection; forward-backward stochastic differential
variational inequality; redundant assets; Tikhonov regularization;
minimum-norm solution.

\section{Introduction}

Continuous-time mean-variance portfolio selection occupies a central
position at the intersection of mathematical finance and stochastic control (see \cite{YongZhou1999}).
Since the single-period formulation of Markowitz (see \cite{m59}), the mean--variance principle
has remained the standard paradigm for balancing expected return against risk,
trading a prescribed variance against a target mean. The continuous-time
extension inherits the distinctive feature of the criterion. Rather than being
the expectation of a state-dependent utility, the objective depends on both the
first and the second moments of terminal wealth, and is therefore not a standard
separable expected-utility functional. This global dependence on the law of
terminal wealth has two immediate consequences: the pointwise Hamiltonian is
not of the classical expected-utility form, and the preferences are dynamically
inconsistent. A strategy that is optimal at the initial time need not remain
optimal from the viewpoint of an intermediate date, because the variance term
is not additive over time.

The pre-commitment problem, where the investor commits at time zero to a single
control rule, can nevertheless be embedded in a stochastic linear--quadratic
(LQ) framework. The key device, introduced in the foundational treatment of
Zhou and Li \cite{zl2000}, is a scalar multiplier attached to the target
mean. The constrained mean--variance problem is then converted into a family of
unconstrained problems indexed by this multiplier, and the efficient frontier
is recovered by varying it. For unconstrained complete markets this reduction
yields explicit Riccati-based feedbacks. Portfolio constraints, however,
complicate the reduction at every stage. The pointwise stationarity condition
can no longer be written as an algebraic equation solved for the control; it
becomes a variational inequality (VI) in which the control is projected onto
the feasible set. When the feasible set is a box, active bounds lead to
complementarity between the switching component and the bound multipliers, and
the coupled optimality system, the forward wealth equation, backward adjoint
equation, and algebraic constraint, is genuinely forward-backward rather than
a single backward equation.

This paper develops the constrained optimality system as a forward-backward
stochastic differential variational inequality (FBSDVI). The terminology is
intended to make explicit three coupled components: the forward wealth SDE that
propagates the state, the backward stochastic adjoint equation that carries the
shadow price and its martingale integrand \(Z\), and the algebraic VI that
selects the portfolio at each time. In doing so, we connect the theory of
constrained stochastic LQ control with the mature theory of deterministic
differential variational inequalities \cite{PangStewart2008,HanEtAl2012}, and
with the more recent numerical literature on stochastic differential
variational inequalities \cite{ZhangEtAl2023}. The stochastic formulation also
brings into view a term that is easy to lose in a deterministic analogy.
Because the portfolio enters the diffusion coefficient of the wealth equation,
the martingale integrand \(Z\) of the adjoint process contributes a
diffusion-induced term \(\sigma Z\) to the switching vector, alongside the
drift contribution \(B Y\). Omitting this term destroys the variational
structure and, in particular, the orthogonality identities on which the
well-posedness analysis rests.

A central theme of this paper is that the framework remains valid under
assumptions substantially weaker than those customarily imposed. We work with
random, uniformly bounded predictable coefficients, and we do not require the
volatility matrix to be non-degenerate: no uniform ellipticity condition
\(\xi^\top\sigma\sigma^\top\xi\ge c|\xi|^2\) is assumed. Discounting and
It\^o's isometry yield the affine representation
\(X_T^u=X_T^0+\mathcal L u\), where
\(\mathcal L:\mathbb H^2\to L^2(\mathcal F_T)\) is a bounded linear map from
the control space to the space of terminal contingent claims. When the market
contains redundant assets, \(\mathcal L\) has a nontrivial null space. Controls
that differ by an element of \(\ker\mathcal L\) generate the same terminal
wealth and hence the same objective value; the solution set is a closed convex
set, an affine subspace in the unconstrained case, rather than a single point,
and the optimal portfolio is not unique. From the Hamiltonian viewpoint this is
the singular-arc regime: in the absence of a running quadratic portfolio
penalty, the pointwise Hamiltonian is affine in the control. A nonzero switching
component forces the corresponding bound, whereas a zero switching component
leaves the first-order conditions unable to select a unique portfolio. Thus the
familiar ``bang--bang'' description is a local active-set characterization, not
a universal characterization of the optimizer.

To resolve the non-uniqueness in a principled way, we introduce a
Tikhonov-type squared control penalty. The regularized objective is strictly
convex and admits a unique optimizer for every positive penalty parameter, and
the corresponding projection VI is well defined on feasible sets that are not
required to be bounded. The main analytical result concerns the
vanishing-regularization limit. As the penalty parameter tends to zero, the
regularized controls converge strongly, not merely weakly, in \(\mathbb H^2\)
to the minimum-\(\mathbb H^2\)-norm element of the original solution set, while
the unregularized optimal objective value is recovered. This supplies a
canonical, equal-split portfolio in the redundant-asset setting. A
complementary measure-valued HJB analysis (cf. \cite{PW2017, PW18, gpw23}) using the Lions derivative reduces the
infinite-dimensional equation to finite-dimensional Riccati equations that
reproduce the explicit feedback and minimum-energy selection in the limit; a
numerical example illustrates convergence of the control paths, expected
terminal wealth, and energy.

Several qualifications should be stated at the outset. We study the
pre-commitment problem. Because mean-variance preferences are dynamically
inconsistent, an equilibrium strategy for a succession of future selves would
require a separate equilibrium-control system and is not treated here. Within
the pre-commitment setting, our contributions are the following. First, we
formulate the constrained problem as an FBSDVI and identify the full switching
vector, including the diffusion-induced term \(\sigma Z\). Second, we establish
the framework with random coefficients and without non-degeneracy, and we
characterize the non-uniqueness induced by the null space of the
control-to-state map. Third, we prove that vanishing Tikhonov regularization
selects the minimum-norm optimizer by strong convergence, on feasible sets that
may be unbounded, and we recover the same selection through the measure-valued
HJB--Riccati reduction.


The stochastic LQ approach to continuous-time mean-variance portfolio
selection has evolved along a clear line. Zhou and Li \cite{zl2000} first
embedded the unconstrained mean--variance problem in a stochastic LQ framework
through a scalar multiplier attached to the target mean, recovering the
efficient frontier through Riccati-based feedback. Lim and Zhou \cite{lz2002} studied the continuous-time mean-variance portfolio selection with random interest rates, appreciation, and volatility in a complete market. They formulate it as a linearly constrained stochastic LQ problem and solve it through a backward stochastic Riccati equation, deriving explicit efficient portfolios, the efficient frontier, and a mutual fund theorem. The constrained case was
addressed by Li, Zhou, and Lim \cite{LZL2002}, who treated a
short-selling prohibition as the non-negative orthant; their analysis, however,
relies on deterministic coefficients and a uniform non-degeneracy condition
\(\sigma\sigma^\top\ge\delta I\), and proceeds through the viscosity solution of
an HJB equation because the constrained Riccati equation has no smooth
solution. The random-coefficient case was later resolved by Hu and Zhou
\cite{hz05}, who allowed the control-weighting matrix to be singular and,
using two nonlinear extended stochastic Riccati equations together with
Tanaka's formula, obtained a closed-form bang-bang feedback on a cone. That
result is nevertheless restricted to a one-dimensional state and to conic
feasible sets, whose dual cone permits an explicit feasibility characterization;
for non-conic constraints only an implicit feasibility condition is available,
and the multi-dimensional extension is left open.

Our work goes beyond these restrictions in four directions.

\begin{itemize}
  \item First, we formulate optimality as a
forward-backward stochastic differential variational inequality (FBSDVI) over
a general nonempty closed convex feasible set that is not required to be
bounded, thereby covering box constraints and the unconstrained case rather
than cones alone.
  \item Second, we do not impose non-degeneracy of the volatility (see \cite{zl2000, lz2002, LZL2002, hz05}):
the control-to-state operator \(\mathcal L\) mapping the control space into the
space of terminal contingent claims may have a nontrivial null space, which we
identify as the precise source of non-uniqueness, an issue that the ESRE
feedback does not address.

 \item Third, the whole argument is formulated at the
Hilbert-space operator level, through the control-to-state map \(\mathcal L\),
the orthogonal projection \(\Pi_C\) with convex and closed set $C$, and the decomposition along
\(\ker\mathcal L\), and never invokes coordinates or the dimension of the
state. It therefore applies unchanged to finite-, multi-, and
infinite-dimensional state and control spaces, and does not rely on the scalar
state required by explicit Riccati-based feedback.

  \item Finally, and most importantly, we provide a selection
mechanism absent in the conic setting. A Tikhonov-type \(L^{2}\) penalty
renders every regularized problem strictly convex and uniquely solvable, and as
the penalty vanishes the regularized controls converge strongly in \(\mathbb
H^{2}\) to the minimum-\(\mathbb H^{2}\)-norm element of the original solution
set while the unregularized optimal value is recovered. A complementary
measure-valued HJB analysis (\cite{PW2017, PW18, gpw23}) based on the Lions derivative reduces the
infinite-dimensional equation to finite-dimensional Riccati equations and
reproduces the same minimum-energy selection, thus yielding a canonical
portfolio precisely in the degenerate, redundant-asset regime where the conic
feedback remains silent.
\end{itemize}


\paragraph{Two-redundant-risky-asset example.}
\label{sec:financial-interpretation}

We provide the two-redundant-risky-asset example to validate the theoretical result stated above. Although the model is deliberately simple, it captures a common situation in portfolio management: several traded securities generate the same exposure to the underlying sources of risk, so the portfolio construction problem does not uniquely determine how the optimal exposure should be allocated among them.

We consider two risky assets with identical excess return \(b>0\) and identical volatility \(\sigma\). The wealth process is
\begin{equation}
\label{eq:financial-wealth}
dX_t
=
b\bigl(u_{1,t}+u_{2,t}\bigr)\,dt
+
\sigma\bigl(u_{1,t}+u_{2,t}\bigr)\,dW_t,
\qquad X_0=x_0,
\end{equation}
where \(u_{i,t}\) denotes the dollar amount invested in asset \(i\). Define the aggregate risky position and the long--short position by $
v_t:=u_{1,t}+u_{2,t},
w_t:=u_{1,t}-u_{2,t}.
$
Then the wealth equation becomes
\begin{equation}
\label{eq:wealth-v-only}
dX_t=bv_t\,dt+\sigma v_t\,dW_t.
\end{equation}
Thus, only the aggregate position \(v_t\) has an effect on future wealth. The long--short position \(w_t\) is a self-financing spread position: it changes the individual holdings but leaves the total exposure to the Brownian shock unchanged. The two assets are redundant in the sense that they provide the same risk exposure. Indeed, the control-to-state map depends on \(u\) only through the linear functional $
u=(u_1,u_2)\longmapsto u_1+u_2.
$
Its kernel is $
\label{eq:financial-kernel}
\ker \mathcal L
=
\operatorname{span}\{(1,-1)\}.
$
A perturbation of the form $
q_t=(a_t,-a_t)
$
changes neither the aggregate position nor the wealth process: $
q_{1,t}+q_{2,t}=0.
$
Consequently, if a control \(u=(u_1,u_2)\) is optimal, then, for any square-integrable adapted process \(a_t\), the perturbed control $
u_1+a_t, u_2-a_t
$
generates the same terminal wealth and the same mean--variance objective value. This nonuniqueness is not a mathematical artifact alone. It reflects a genuine portfolio-management issue: when several instruments are equivalent with respect to the relevant risk factors, the mean--variance criterion alone determines the factor exposure but not the instrument-level allocation.

Let $
\label{eq:market-price-risk}
\theta:=\frac{b}{\sigma}
$
be the market price of risk, and define the state-price density
$
\label{eq:rho-financial}
\rho_t
=
\exp\left(
-\theta W_t-\frac12\theta^2t
\right).
$
It satisfies $
\label{eq:rho-sde}
d\rho_t=-\theta\rho_t\,dW_t.
$
The adjoint BSDE associated with the mean-variance terminal objective is
\begin{equation}
\label{eq:adjoint-bsde-financial}
dY_t=Z_t\,dW_t,
\qquad
Y_T=\lambda\bigl(X_T-\mathbb E[X_T]\bigr)-1, \qquad \lambda>0.
\end{equation}
The process \(Y_t\) is the marginal valuation of terminal wealth along the optimal trajectory. In the present complete-market example, it is proportional to the state-price density
$Y_t=-\rho_t,
Z_t=\theta\rho_t.
$
The Hamiltonian gradient with respect to each individual control is
$
G_t=bY_t+\sigma Z_t.
$
Because the two assets have the same excess return and the same volatility, their marginal contributions are identical. Hence the adjoint BSDE determines only the aggregate position \(v_t\). It contains no information that distinguishes asset \(1\) from asset \(2\).

The state-price density also gives the market value of the optimal terminal payoff. Since \(\rho_tX_t\) is a martingale under the model coefficients, the replicating wealth at time \(t\) is
$
X_t^*
=
\mathbb E\left[
\frac{\rho_T}{\rho_t}X_T^*
\,\middle|\,\mathcal F_t
\right].
$
For the optimal terminal payoff
\begin{equation}
\label{eq:optimal-terminal-payoff}
X_T^*
=
m^*+\frac1\lambda-\frac1\lambda\rho_T,
\qquad
m^*=\mathbb E[X_T^*],
\end{equation}
the conditional Gaussian identities
$
\mathbb E\!\left[\frac{\rho_T}{\rho_t}\,\middle|\,\mathcal F_t\right]=1,
\qquad
\mathbb E\!\left[
\left(\frac{\rho_T}{\rho_t}\right)^2
\,\middle|\,\mathcal F_t
\right]
=
e^{\theta^2(T-t)}
$
give
\begin{equation}
\label{eq:optimal-wealth-financial}
X_t^*
=
m^*+\frac1\lambda
-
\frac1\lambda e^{\theta^2(T-t)}\rho_t.
\end{equation}
Thus the optimal terminal payoff is an affine function of the state-price density. This is the standard mean--variance efficiency result: the optimal strategy holds a riskless or deterministic component and a short position in the state-price density, which is equivalent to holding a mean-variance tangent portfolio with a stochastic, state-dependent scale.
Although \(w_t\) does not affect wealth, it is not costless from a portfolio-construction perspective. Indeed,
$
u_{1,t}^2+u_{2,t}^2
=
\frac12\bigl(v_t^2+w_t^2\bigr).
$
Therefore the squared control norm decomposes as
\begin{equation}
\label{eq:H2-decomposition-financial}
\|u\|_{\mathbb H^2}^2
=
\frac12\mathbb E\int_0^T v_t^2\,dt
+
\frac12\mathbb E\int_0^T w_t^2\,dt.
\end{equation}
The first term is the energy required to obtain the optimal factor exposure. The second term is purely unremunerated spread exposure: it neither increases expected terminal wealth nor improves the variance of terminal wealth.

For a fixed optimal aggregate position \(v_t^*\), minimizing \(\|u\|_{\mathbb H^2}^2\) is therefore equivalent to minimizing $
\mathbb E\int_0^T w_t^2\,dt.
$
The unique minimizer is $
w_t^\dagger=0,
$
and hence $
\label{eq:equal-split-financial}
u_{1,t}^\dagger=u_{2,t}^\dagger=\frac{v_t^*}{2}.
$
The unregularized optimal aggregate position is $
v_t^*
=
\frac{b}{\lambda\sigma^2}
e^{\theta^2(T-t)}\rho_t,
$
so the minimum-energy strategy is
\begin{equation}
\label{eq:minimum-energy-controls}
u_{1,t}^\dagger=u_{2,t}^\dagger
=
\frac{b}{2\lambda\sigma^2}
e^{\theta^2(T-t)}\rho_t.
\end{equation}
Financially, \(\|u\|_{\mathbb H^2}^2\) can be interpreted as a measure of the sensitivity to model-estimation errors, market impact in a quadratic transaction-cost approximation, funding or margin requirements, etc.

Under any of these interpretations, a nonzero long--short position in two redundant assets consumes capital or risk budget without being compensated by additional expected utility. The minimum-norm rule therefore allocates the required exposure as evenly as possible across symmetric substitutes.

It is important to note that equal weighting is not an additional diversification assumption in this example. It is the consequence of minimizing control energy after the optimal aggregate risk exposure has been determined. If the two assets have different transaction costs, liquidity measures, margin requirements, or position limits, then the symmetric rule is replaced by an appropriately weighted minimum-cost allocation.

Financially, the vanishing regularization acts as a tie-breaking rule. At first, for large \(\varepsilon\), it substantially reduces risk taking. As the penalty becomes smaller, the portfolio's expected return and risk exposure approach their unregularized mean--variance optimum. Nevertheless, throughout the entire path, the regularization suppresses redundant long--short exposure. In the limit, the optimal factor exposure is restored, while the selected instrument-level portfolio remains the minimum-capital solution.

The remainder of the paper is organized as follows. Section~2 sets out the
problem formulation and notation. Section~3 describes the market, information
structure, and admissible portfolios, and derives the portfolio FBSDVI and its
equivalent algebraic forms. Section~4 develops the regularization and response
map, establishing the reduced FBSDEs for both bounded and unbounded control
sets. Section~5 introduces the model setup, uses the Lions derivative and
adjoint BSDE to obtain the explicit aggregate solution and the minimum-energy
split, and analyzes the finite-$\varepsilon$ feedback, the
vanishing-regularization limit, and a numerical illustration. Section~6
concludes, and the appendix collects auxiliary proofs and projection
properties.

\section{Problem Formulation and Preliminaries}

\label{sect2}

\subsection{Notation}

We make use of the following notation:
\begin{equation*}
\begin{tabular}{crl}
$M^{\top }$ & : & the transpose of any matrix or vector $M$; \\
$\Vert M\Vert $ & : & $\sqrt{\sum_{i,j}m_{ij}^{2}}$ for any matrix or vector
$M=(m_{ij})$; \\
$\mathbb{R}^{n}$ & : & $n$ dimensional real Euclidean space; \\
$\mathbb{R}_{+}^{n}$ & : & the subset of $\mathbb{R}^{n}$ consisting of
elements with nonnegative components.%
\end{tabular}%
\end{equation*}%
Consider a finite time horizon $[0,T]$ for fixed $T>0$. We assume $(\Omega ,%
\mathcal{F},\{\mathcal{F}_{t}\}_{0\leq t\leq T},P)$ is a complete, filtered
probability space on which a standard $m$-dimensional Brownian motion $%
W(t)\equiv (W^{1}(t),\cdots ,W^{m}(t))^{\top }$ is defined. Moreover, we will adopt the following notations, unless
otherwise specified.

\begin{itemize}
\item Define
\[
\mathbb H^2(0,T;\R^m)
=
\left\{
v:\ v\ \text{predictable},\quad
\E\int_0^T |v_t|^2\,dt<\infty
\right\},
\]
and
\[
\mathbb S^2(0,T;\R)
=
\left\{
X:\ X\ \text{adapted and continuous},\quad
\E\sup_{0\le t\le T}|X_t|^2<\infty
\right\}.
\]

\item For any normed space $\left( E,\left\vert \cdot \right\vert \right) $,
$\mathcal{P}(E)$ is the set of all probability measures on $\left(
E,\left\vert \cdot \right\vert \right) $. For any $p\in N$, $\mathcal{P}%
_{p}(E)$ is the set of all probability measures of $p$th order on $\left(
E,\left\vert \cdot \right\vert \right) $, define by
\begin{equation*}
\mathcal{P}_{p}(E)=\left\{ \mu \in \mathcal{P}(E)\left\vert \left\Vert \mu
\right\Vert _{p}=\left( \int_{E}\left\vert x\right\vert ^{p}\mu \left(
\mathrm{d}x\right) \right) ^{\frac{1}{p}}<\infty \right. \right\} .
\end{equation*}%
For any probability measures $\mu ,\mu ^{\prime }$ in $\mathcal{P}_{p}(E)$,
the $p$th order Wasserstein distance on $\mathcal{P}_{p}(E)$ is defined as
\begin{equation*}
\mathcal{W}_{p}\left( \mu ,\mu ^{\prime }\right) =\inf_{\pi }\left(
\int_{E\times E}\left\vert z-z^{\prime }\right\vert ^{p}\pi \left(
dz,dz^{\prime }\right) \right) ^{\frac{1}{p}},
\end{equation*}%
where $\pi $ is a coupling of $\mu $ and $\mu ^{\prime }$ in the sense that $%
\pi \in \mathcal{P}(E\times E)$ with marginals $\mu $ and $\mu ^{\prime }$. $%
L^{2}(E)$ is the space of all square-integrable functions $f:E\rightarrow
\mathbb{R}$.

\item The Lions derivative of a functional $f$, introduced in \cite{lions06}%
, is defined through the lift of $f.$ Take any function $f:\mathcal{P}_{2}(%
\mathbb{R}^{n})\rightarrow \mathbb{R}$, and let $\tilde{f}:L^{2}(\mathcal{F};%
\mathbb{R}^{n})\rightarrow \mathbb{R}$ be a life of $f$ such that $\tilde{f}%
\left( \xi \right) =f\left( \mathbb{P}_{\xi }\right) $ for any $\xi \in
L^{2}(\mathcal{F};\mathbb{R}^{n}).$ Then, $\tilde{f}$ is differentiable in
the Fr\'{e}chet sense at $\xi _{0}$ if there exists a linear continuous
mapping $D\tilde{f}\left( \xi _{0}\right) :L^{2}(\mathcal{F};\mathbb{R}%
^{n})\rightarrow \mathbb{R}$ such that%
\begin{equation*}
\tilde{f}\left( \xi \right) -\tilde{f}\left( \xi _{0}\right) =\mathbb{E}%
\left[ D\tilde{f}\left( \xi _{0}\right) \centerdot \left( \xi -\xi
_{0}\right) \right] +o\left( \left\Vert \xi -\xi _{0}\right\Vert
_{L^{2}}\right) ,
\end{equation*}%
as $\left\Vert \xi -\xi _{0}\right\Vert _{L^{2}}\rightarrow 0.$Whenever $%
\tilde{f}$ is the lift of a function $f$ in $\mathcal{P}_{2}(\mathbb{%
R}^{n})$, the law of $D\tilde{f}\left( \xi _{0}\right) $ depends on $\xi
_{0} $ only via its law $\mathbb{P}_{\xi _{0}}$, and
\begin{equation}
D\tilde{f}\left( \xi _{0}\right) =g_{0}\left( \xi _{0}\right) ,  \label{lg}
\end{equation}%
for some Borel function $g_{0}:\mathbb{R}^{n}\rightarrow \mathbb{R}^{n}$ .
(See e.g., \cite{cd18}, Chapter 5 and \cite{gt19}). The Lions derivative is
thus well defined.
\end{itemize}

\begin{definition}
\label{d1}We say $f$ is differentiable at $\mu _{0}=\mathbb{P}_{\xi _{0}}\in
\mathcal{P}_{2}(\mathbb{R}^{n})$ if its lift function $\tilde{f}$ is Fr\'{e}%
chet differentiable at $\xi _{0}$; and in this case, the function $g_{0}$ in
(\ref{lg}) is called the Lions derivative of $f$ at $\mu _{0}$, and denoted
as $\partial _{\mu }f\left( \mu _{0},\cdot \right) $.
\end{definition}

\section{Market, information, and admissible portfolios}

\label{sect3} Let $(\Omega ,\mathcal{F},(\mathcal{F}_{t})_{0\leq t\leq T},%
\mathbb{P})$ satisfy the usual conditions and support a $d$-dimensional
Brownian motion $W$. The filtration is the augmented Brownian filtration
whenever martingale representation is used. There is one money-market
account and $m$ risky assets:
\begin{align}
dS_{t}^{0}& =r_{t}S_{t}^{0}dt,  \label{eq:bank} \\
dS_{t}^{i}& =S_{t}^{i}\left( b_{t}^{i}dt+\sum_{j=1}^{d}\sigma
_{t}^{ij}dW_{t}^{j}\right) ,\qquad i=1,\ldots ,m.  \label{eq:risky}
\end{align}%
Let $B_{t}=(b_{t}^{1}-r_{t},\ldots ,b_{t}^{m}-r_{t})^{\top}$ and $%
\sigma _{t}=(\sigma _{t}^{ij})\in \mathbb{R}^{m\times d}$. If $u_{t}^{i}$ is
the dollar amount in risky asset $i$, the residual $X_{t}-\mathbf{1}^{%
\top}u_{t}$ is invested in the money-market account. Self-financing
wealth therefore satisfies
\begin{equation}
dX_{t}=(r_{t}X_{t}+B_{t}^{\top}u_{t})dt+u_{t}^{\top}\sigma
_{t},dW_{t}, \qquad X_{0}=x_{0}.  \label{eq:wealth}
\end{equation}


\begin{enumerate}
\item[\textbf{(A1)}] \label{ass:market}Assume that $r$, $B$, and $\sigma $
are bounded $\left\{ \mathcal{F}_{t}\right\} _{t\geq 0}$-adapted processes. The $\left\{ \mathcal{F}_{t}\right\} _{t\geq 0}$-adapted bounds $\ell ,h\in
\mathbb{H}^{2}(\mathbb{R}^{m})$ satisfy $\ell _{t}\leq h_{t}$ componentwise
and, for some deterministic $C_{K}$, $|\ell _{t}|+|h_{t}|\leq C_{K}$. The
initial wealth $x_{0}$ is deterministic.
\end{enumerate}

Define the random box $K_{t}(\omega )=[\ell _{t}(\omega ),h_{t}(\omega )]$
and
\begin{equation}
\mathcal{U}_{ad}:=\{u\in \mathbb{H}^{2}(\mathbb{R}^{m}):u_{t}\in K_{t}\quad
dt\otimes d\mathbb{P}\text{-a.e.}\}.  \label{eq:U}
\end{equation}%
Under (A1), $\mathcal{U}_{ad}$ is nonempty, convex, closed, bounded, and hence
weakly compact in the Hilbert space $\mathbb{H}^{2}(\mathbb{R}^{m})$. For
every $u\in \mathcal{U}_{ad}$, \eqref{eq:wealth} has a unique $X^{u}\in \mathbb{S}%
^{2}(\mathbb{R})$. We first analyze the scenario in which the control set $K_t$ is bounded, and subsequently extend our analysis to the unbounded case where $K_t = \mathbb{R}_+^m$. In both cases, the terminal state depends affinely on the control variable. More precisely, there is a
bounded affine map
\begin{equation}  \label{eq:endpoint-map}
\mathcal{T}:\mathbb{H}^2(\mathbb{R}^m)\to L^2_{\mathcal{F}_T}(\Omega),\qquad
\mathcal{T}(u)=X_T^u.
\end{equation}
Its linear part will be denoted by $\mathcal{L}$, so $X_T^u=X_T^0+\mathcal{L}%
u$. The terminal-mean functional
$
a(u):=\mathbb{E}[X_T^u]
$
is continuous and affine. Indeed, the wealth process satisfies
\begin{equation}  \label{eq:wealth-SDE}
dX_t = \bigl(r_tX_t+B_t^\top u_t\bigr)\,dt + u_t^\top\sigma_t\,dW_t, \qquad
X_0=x_0.
\end{equation}
For simplicity, we consider one wealth in SDE (\ref{eq:wealth-SDE}). Our method can actually deal with general case (see Remark \ref{rem:dimension-independent}).
Define the bank-account growth factor and the discount factor by
$
\beta_t = \exp\left(\int_0^t r_s\,ds\right), D_t=\beta_t^{-1},
$
respectively. Since $
dD_t=-r_tD_t\,dt,
$
applying It\^o's product rule to $D_tX_t$ gives
\begin{equation*}
\begin{aligned} d(D_tX_t) &=D_t\,dX_t+X_t\,dD_t\\ &=
D_t\bigl(r_tX_t+B_t^\top u_t\bigr)\,dt +D_tu_t^\top\sigma_t\,dW_t
-r_tD_tX_t\,dt\\ &= D_tB_t^\top u_t\,dt + D_tu_t^\top\sigma_t\,dW_t.
\end{aligned}
\end{equation*}
Thus the terms involving $r_tX_t$ cancel. Integrating from $0$ to $T$, and using $D_0=1$, yields
\begin{equation}  \label{eq:discounted-terminal-wealth}
D_TX_T = x_0 +\int_0^T D_tB_t^\top u_t\,dt +\int_0^T
D_tu_t^\top\sigma_t\,dW_t.
\end{equation}
Multiplying by $\beta_T=D_T^{-1}$, we obtain
\begin{equation}  \label{eq:terminal-wealth-expanded}
X_T^u = \beta_Tx_0 +\beta_T\int_0^T D_tB_t^\top u_t\,dt +\beta_T\int_0^T
D_tu_t^\top\sigma_t\,dW_t.
\end{equation}
Set
$
X_T^0=\beta_Tx_0
$
and define
\begin{equation}  \label{eq:linear-control-operator}
\mathcal{L}u = \beta_T\int_0^T D_tB_t^\top u_t\,dt + \beta_T\int_0^T
D_tu_t^\top\sigma_t\,dW_t.
\end{equation}
Here $X_T^0$ is the terminal wealth obtained by holding no risky assets and
investing only in the money-market account. Indeed, when $u\equiv0$, the
wealth equation reduces to $
dX_t=r_tX_t\,dt,
$
whose solution is $X_t=\beta_tx_0$.

The operator $\mathcal{L}$ is linear in $u$. Consequently, equation~%
\eqref{eq:terminal-wealth-expanded} can be written as
\begin{equation}  \label{eq:affine-terminal-map}
X_T^u=X_T^0+\mathcal{L}u.
\end{equation}
Therefore, the terminal-wealth map
\begin{equation}  \label{eq:terminal-map}
\mathcal{T}(u):=X_T^u=X_T^0+\mathcal{L}u
\end{equation}
is affine: it consists of a linear part plus a constant term. It is a linear
map in the strict sense only when $x_0=0$, in which case the constant term $%
X_T^0$ vanishes.

We now present a general result.

\begin{lemma}
\label{lem:affine-n}
Let the $n$-dimensional wealth vector $X$ satisfy
\begin{equation}
dX_t=\bigl(r_t I_n X_t+U_t\beta_t\bigr)\,dt+U_t\Sigma_t\,dW_t,
\qquad X_0=x_0,
\end{equation}
where $U_t\in\mathbb R^{n\times m}$, $\beta_t\in\mathbb R^m$,
$\Sigma_t\in\mathbb R^{m\times d}$, and $r_t$ is scalar. Assume that $r$, $\beta$ and $\Sigma$ are adapted and bounded, and that
$U\in H^2_{n\times m}:=L^2_\mathcal F([0,T];\mathbb R^{n\times m})$ (if $n=1$, set
$\mathbb H^2=\mathbb H^2_{1\times m}$). Define the adapted discount factor and the terminal accumulation
factor by
\begin{equation}
\rho_t:=\exp\!\left(-\int_0^t r_\tau\,d\tau\right),\qquad
R(T,0):=\rho_T^{-1},\qquad X_T^0:=R(T,0)\,x_0 .
\end{equation}
Then the terminal wealth admits the affine representation
\begin{equation}
X_T=X_T^0+\mathcal LU,
\end{equation}
where the bounded linear map $\mathcal L:H^2_{n\times m}\to L^2(F_T;\mathbb R^n)$ is
\begin{equation}
\mathcal{L}U=R(T,0)\int_{0}^{T}\rho _{s}\,U_{s}\beta
_{s}\,ds+R(T,0)\int_{0}^{T}\rho _{s}\,U_{s}\Sigma _{s}\,dW_{s}.  \label{Lu}
\end{equation}
The stochastic integral has the adapted integrand $\rho_s U_s\Sigma_s$, and the terminal
factor $R(T,0)$ is kept outside. When $r$ is deterministic, $R(T,0)\rho_s=R(T,s)$ and
this reduces to
$LU=\int_0^T R(T,s)U_s\beta_s\,ds+\int_0^T R(T,s)U_s\Sigma_s\,dW_s$.
\end{lemma}

\begin{proof}
The adapted discount factor obeys $d\rho_t=-r_t\rho_t\,dt$, $\rho_0=1$, hence has finite
variation and zero quadratic variation. Applying It\^o's product rule to $\rho_t X_t$,
\begin{equation}
d(\rho_t X_t)=\rho_t\,dX_t+X_t\,d\rho_t
=\rho_t U_t\beta_t\,dt+\rho_t U_t\Sigma_t\,dW_t,
\end{equation}
because the two terms containing $r_t I_n X_t$ cancel. Integrating and using $\rho_0=1$,
\begin{equation}
\rho_T X_T=x_0+\int_0^T\rho_s U_s\beta_s\,ds
+\int_0^T\rho_s U_s\Sigma_s\,dW_s .
\end{equation}
The stochastic integral is a well-defined It\^o integral since $\rho_s U_s\Sigma_s$ is
predictable. Multiplying by $\rho_T^{-1}=R(T,0)$ as an overall factor gives the asserted
representation. It remains to verify boundedness. Since $|r_t|\le\bar r$,
\begin{equation}
e^{-\bar rT}\le\rho_t\le e^{\bar rT},\qquad
\|R(T,0)\|_\infty=\|\rho_T^{-1}\|_\infty\le e^{\bar rT}.
\end{equation}
Writing $C_\beta,C_\Sigma$ for the uniform bounds on $\beta,\Sigma$, Cauchy--Schwarz and
It\^o's isometry give
\begin{equation}
\mathbb E\!\left|R(T,0)\int_0^T\rho_s U_s\beta_s\,ds\right|^2
\le \|R(T,0)\|_\infty^2\,T\,
\mathbb E\!\int_0^T\rho_s^2|U_s\beta_s|^2\,ds
\le C_1\|U\|_{H^2}^2,
\end{equation}
and
\begin{equation}
\mathbb E\!\left|R(T,0)\int_0^T\rho_s U_s\Sigma_s\,dW_s\right|^2
\le \|R(T,0)\|_\infty^2\,
\mathbb E\!\int_0^T\rho_s^2|U_s\Sigma_s|^2\,ds
\le C_2\|U\|_{H^2}^2 .
\end{equation}
Combining them yields $\mathbb E|\mathcal LU|^2\le C\|U\|_{H^2}^2$, with $C$ independent of
$U,n,m,d$. Hence $\mathcal L$ is bounded and the construction is dimension independent.
\end{proof}
\begin{remark}
\label{rem:dimension-independent}

From Lemma \ref{lem:affine-n}, we see \(\mathcal L:\mathbb H^{2}(0,T;\mathbb R^{m})\to L^{2}(\mathcal F_T;\mathbb R^{n})\)
is bounded linear. The entire analysis uses only the Hilbert-space geometry of
\(\mathcal L\), namely, its kernel, the lower bound on \((\ker\mathcal L)^\perp\), the
adjoint \(\mathcal L^{*}\), and the projection onto the feasible set, and never
invokes coordinates or the value of \(n\). By contrast, the explicit
Riccati/ESRE approach (see \cite{hz05}) requires an \(n\times n\) matrix-valued stochastic Riccati
equation, whose multi-dimensional existence, comparison and matrix positive
part constitute the open problem already noted in the scalar-state literature,
and it typically imposes the non-degeneracy condition (see \cite{LZL2002, hz05})
\(\sigma\sigma^\top\ge\delta I\). The operator argument applies unchanged for
every finite \(n\) (and for infinite-dimensional state spaces), accommodates a
nontrivial \(\ker\mathcal L\) without non-degeneracy, and selects the
minimum-norm minimizer through the vanishing Tikhonov penalty (see Section~4); this is precisely
why the formulation is dimension-independent.
\end{remark}

\begin{remark}
We have another expression for operator $\mathcal L$ via the Kronecker product as follows.
Set
$
u_t := \vecop(U_t) \in \R^{nm}.
$
Applying our identity to the drift control term gives
$
\vecop(U_t\beta_t) = (\beta_t^\top \otimes I_n)\,u_t.
$
For the diffusion term, proceed column by column. The $j$-th column of
$U\Sigma$ is $U\Sigma_{\cdot j}$; replacing $\beta$ in the identity by
$\Sigma_{\cdot j}$ yields
\begin{equation}
\vecop\!\big(U_t(\Sigma_t)_{\cdot j}\big)
= \big((\Sigma_t)_{\cdot j}^\top \otimes I_n\big)\,u_t.
\end{equation}
Using $U\Sigma\,dW = \sum_{j=1}^{d}(U\Sigma)_{\cdot j}\,dW^j$, the
vectorized wealth equation involves only the vector $u$:
\begin{equation}
dX_t = r_t I_n X_t\,dt
     + (\beta_t^\top \otimes I_n)\,u_t\,dt
     + \sum_{j=1}^{d}\big((\Sigma_t)_{\cdot j}^\top \otimes I_n\big)\,u_t\,dW_t^j.
\end{equation}
Let
$\rho_s:=\exp\!\big(-\int_0^s r_\tau\,d\tau\big)$ be the adapted discount factor and
$R(T,0):=\rho_T^{-1}$ the terminal accumulation factor. Applying the discounting
method (as in the preceding lemma) to the vector equation above gives
\begin{equation}
X_T = X_T^0 + \mathcal L u,
\qquad X_T^0 = R(T,0)\,x_0,
\end{equation}
\begin{equation}
\mathcal L u
= R(T,0)\int_0^T \rho_s\,(\beta_s^\top \otimes I_n)\,u_s\,ds
+ R(T,0)\sum_{j=1}^{d}\int_0^T \rho_s\,
  \big((\Sigma_s)_{\cdot j}^\top \otimes I_n\big)\,u_s\,dW_s^j.
\end{equation}
Each stochastic integrand
$\rho_s\big((\Sigma_s)_{\cdot j}^\top \otimes I_n\big)u_s$ is adapted, while the
terminal factor $R(T,0)$ is kept outside. When $r$ is deterministic, writing
$R(T,s)=\exp\!\big(\int_s^T r_\tau\,d\tau\big)$ gives $R(T,0)\rho_s=R(T,s)$, and the
operator reduces to
\begin{equation}
\mathcal L u
= \int_0^T R(T,s)\,(\beta_s^\top \otimes I_n)\,u_s\,ds
+ \sum_{j=1}^{d}\int_0^T R(T,s)\,
  \big((\Sigma_s)_{\cdot j}^\top \otimes I_n\big)\,u_s\,dW_s^j.
\end{equation}
Here $\mathcal L$ acts strictly on the vector $u$; since $u = \vecop(U)$, it is fully
equivalent to the matrix form (\ref{Lu}), but the notation is no longer mixed.
\end{remark}

For notational convenience, we continue to consider a single wealth variable. For the general case, we put in Subsection \ref{ma} in the end of Section~4. We define the cost functional
\begin{equation}
J_\lambda\left( u\left( \cdot \right) \right) =\mathbb{E}\bigg [\frac{\lambda }{2}%
\left( X_{T}^{2}-\left( \mathbb{E}X_{T}\right) ^{2}\right) -X_{T}\bigg ].
\label{cost1}
\end{equation}%
where $\lambda >0$. When $\lambda $ is a positive constant, the above
criterion reduces to the classical mean-variance portfolio selection, as
studied e.g. in (Li and Zhou 2000 \cite{zl2000}).

We set
\begin{equation*}
\bar{\mu}_{2}=\int_{\mathbb{R}^{d}}x^{2}\mu \left( \mathrm{d}x\right) ,\text{
Var}\left( \mu \right) =\bar{\mu}_{2}-\bar{\mu}^{2},\text{ }\mu \in \mathcal{%
P}_{2}\left( \mathbb{R}\right) ,
\end{equation*}%
where $\bar{\mu}=\int_{\mathbb{R}^{n}}x\mu \left( \mathrm{d}x\right) .$
Therefore, the cost (\ref{cost1}) can be rewritten as
\begin{equation}
J\left( u\left( \cdot \right) \right) =\frac{\lambda }{2}\mathbb{V}\text{ar}%
\left( \mathbb{P}_{X_{T}}\right) -\bar{\mu}\left( \mathbb{P}_{X_{T}}\right) .
\label{cost2}
\end{equation}
\begin{lemma}
\label{lem:mv-quadratic-form}
Let the terminal wealth be the affine function
$
X_T(u)=X_T^{0}+\mathcal L u,
$
where \(\mathcal L\) is a bounded linear operator from the control space into
\(L^{2}\) defined in (\ref{Lu}), and let \(P_0Z:=Z-\E Z\) be the orthogonal projection onto centered
random variables. Then the mean--variance objective
\[
J_\lambda(u)=\frac{\lambda}{2}\operatorname{Var}\!\bigl(X_T(u)\bigr)-\E[X_T(u)]
\]
admits the decomposition
\begin{equation}
J_{\lambda }(u)=\frac{\lambda }{2}\Vert (P_{0}\mathcal{L})u\Vert
_{L^{2}}^{2}+\ell (u)+C,  \label{jop}
\end{equation}
where
\begin{equation}
\ell (u)=\lambda \,\mathbb{E}\!\left[ P_{0}X_{T}^{0}\,(P_{0}\mathcal{L})u%
\right] -\mathbb{E}[\mathcal{L}u]  \label{lop}
\end{equation}
is a bounded linear functional and
\[
C=\frac{\lambda}{2}\|P_0X_T^{0}\|_{L^{2}}^{2}-\E[X_T^{0}].
\]
\end{lemma}
\begin{proof}
For any square-integrable random variable \(Z\),
\[
\operatorname{Var}(Z)=\E[(Z-\E Z)^2]=\|P_0 Z\|_{L^2}^{2}.
\]
Using the linearity of \(P_0\) and \(X_T(u)=X_T^0+\mathcal L u\),
\[
\operatorname{Var}(X_T(u))
=\|P_0(X_T^0+\mathcal L u)\|^2
=\|P_0X_T^0+P_0\mathcal L u\|^2.
\]
By the inner-product identity \(\|a+b\|^2=\|a\|^2+\|b\|^2+2\langle a,b\rangle\),
\[
\operatorname{Var}(X_T(u))
=\|P_0X_T^0\|^2+\|P_0\mathcal L u\|^2
+2\,\E\!\left[P_0X_T^0\cdot P_0\mathcal L u\right].
\]
\[
\E[X_T(u)]=\E[X_T^0]+\E[\mathcal L u].
\]
Substituting above into the objective gives
\[
\begin{aligned}
J_\lambda(u)
=\;&\frac{\lambda}{2}\|P_0\mathcal L u\|^2\\
&+\Bigl(\lambda\,\E[P_0X_T^0\cdot P_0\mathcal L u]-\E[\mathcal L u]\Bigr)\\
&+\Bigl(\frac{\lambda}{2}\|P_0X_T^0\|^2-\E[X_T^0]\Bigr).
\end{aligned}
\]
The first line is the quadratic term \(\frac{\lambda}{2}\|(P_0\mathcal L)u\|^2\),
the second line is \(\ell(u)\), and the third line is the constant \(C\), which
proves the claimed decomposition.

Linearity is immediate. Since \(P_0\) is an orthogonal projection,
\(\|P_0\|=1\); the Cauchy--Schwarz inequality then yields
\[
\big|\E[P_0X_T^0\cdot P_0\mathcal L u]\big|
\le \|P_0X_T^0\|\,\|P_0\mathcal L u\|
\le \|P_0X_T^0\|\,\|\mathcal L\|\,\|u\|,
\]
and
\[
|\E[\mathcal L u]|\le\|\mathcal L u\|\le\|\mathcal L\|\,\|u\|.
\]
Hence \(|\ell(u)|\le\bigl(\lambda\|P_0X_T^0\|+1\bigr)\|\mathcal L\|\,\|u\|\), so
\(\ell\) is a bounded linear functional on the control space. The proof is completed.
\end{proof}
\begin{remark}
\label{rem:kernel-hierarchy}
Three nested notions of redundancy should be distinguished. A perturbation $q$ leaves
the whole state path unchanged only if it annihilates both the drift and the diffusion
coefficients,
\begin{equation}
\mathcal K_{\mathrm{path}}
:=\big\{q\in \mathbb H^2(0,T;\mathbb R^m):
B_t^\top q_t=0,\ \sigma_t^\top q_t=0\ \text{a.e.}\big\}.
\end{equation}
A control in $\ker \mathcal L$ changes no terminal wealth, whereas one in $\ker(P_0\mathcal L)$ produces a
merely deterministic terminal contribution. Hence
\begin{equation}
\mathcal K_{\mathrm{path}}\subseteq\ker \mathcal L\subseteq\ker(P_0\mathcal L),
\end{equation}
and the inclusions may be strict. In the two-redundant-asset example the three spaces
coincide with $\operatorname{span}\{(1,-1)\}$.
\end{remark}

In Hu and Zhou \cite{hz05} and Li et al. \cite{LZL2002}, the expected terminal wealth is fixed in advance through the explicit
constraint $\mathbb E[X_T]=z$. A scalar multiplier then embeds the problem
into a stochastic LQ problem, making $z$ (equivalently, the multiplier) an
exogenous parameter on which the resulting ESREs depend; a parameterized
family is solved and the efficient frontier is recovered afterwards.

Our FBSDVI formulation below presupposes no such target. The forward wealth
equation, the adjoint BSDE (for classical case see \cite{PardouxPeng1990}), and the portfolio variational inequality form a
self-contained optimality system that characterizes the optimal policy
directly, without imposing $\mathbb E[X_T]=z$ at the outset. The
mean-variance trade-off is instead recovered from the adjoint variables and
enters as a parameter of the efficient frontier rather than as an assumption
on the problem. The analysis is therefore independent of any preassigned
level of expected terminal wealth and, when the control-to-state map has a
nontrivial null space, describes the entire solution set rather than a single
target-constrained solution.

\subsection{First variation and It\^o duality}

For any $v\in\mathcal{U}$ and $\varepsilon\in[0,1]$, the convex perturbation
$u^\varepsilon=\bar u+\varepsilon(v-\bar u)$ remains feasible. Put $\delta
u=v-\bar u$. Because the state equation is affine,
$
X^{u^\varepsilon}=\bar X+\varepsilon\eta,
$
where
\begin{equation}  \label{eq:variation}
d\eta_t=(r_t\eta_t+B_t^\mathsf{T}\delta u_t)dt+\delta u_t^\mathsf{T}%
\sigma_t,dW_t,\qquad \eta_0=0.
\end{equation}
The one-sided derivative is
\begin{equation}  \label{eq:first-variation}
DJ(\bar u)[\delta u]=\mathbb{E}[(\lambda (X_{T}-\mathbb{E}X_{T})-1)\eta _{T}].
\end{equation}
Define the Hamiltonian
\begin{equation}
H_{t}(x,u,y,z)=y(r_{t}x+B_{t}^{\mathsf{T}}u)+z^{\mathsf{T}}\sigma _{t}^{%
\mathsf{T}}u,  \label{eq:Hamiltonian}
\end{equation}%
where $y\in \mathbb{R}$, $z\in \mathbb{R}^{d}$. Its derivatives are
\begin{equation}
H_{x}=r_{t}y,\qquad H_{u}=B_{t}y+\sigma _{t}z.  \label{eq:H-derivatives}
\end{equation}%
Let $(Y,Z)\in \mathbb{S}^{2}(\mathbb{R})\times \mathbb{H}^{2}(\mathbb{R}%
^{d}) $ solve the linear mean-field BSDE (see \cite{blj16})
\begin{equation}
dY_{s}=-AY_{s}\,ds+\sum_{j=1}^{d}Z_{s}^{j}\,dW_{s}^{j},Y_{T}=\lambda (X_{T}-%
\mathbb{E}X_{T})-1.  \label{eq:adjoint}
\end{equation}%
Existence and uniqueness follow from standard linear mean field BSDE theory
(cf. \cite{blj16}).

Applying It\^{o}'s formula to $Y_{t}\eta _{t}$ gives
\begin{align}
d(Y_{t}\eta _{t})& =Y_{t}d\eta _{t}+\eta _{t}dY_{t}+d[Y,\eta ]_{t}  \notag \\
& =\left[ Y_{t}(r_{t}\eta _{t}+B_{t}^{\mathsf{T}}\delta
u_{t})-r_{t}Y_{t}\eta _{t}+Z_{t}^{\mathsf{T}}\sigma _{t}^{\mathsf{T}}\delta
u_{t}\right] dt  \notag \\
& \quad +\left[ Y_{t}\delta u_{t}^{\mathsf{T}}\sigma _{t}+\eta _{t}Z_{t}^{%
\mathsf{T}}\right] dW_{t}  \notag \\
& =\delta u_{t}^{\mathsf{T}}(B_{t}Y_{t}+\sigma _{t}Z_{t})dt+dM_{t}.
\label{eq:Ito-product}
\end{align}%
The term $Z_{t}^{\mathsf{T}}\sigma _{t}^{\mathsf{T}}\delta u_{t}dt$ is the
quadratic covariation of the martingale parts. After localization, square
integrability makes $M$ a true martingale. Integrating and taking
expectations yields
\begin{equation}
\mathbb{E}[(\lambda (X_{T}-\mathbb{E}X_{T})-1)\eta _{T}]=\mathbb{E}%
\int_{0}^{T}\delta u_{t}^{\mathsf{T}}(B_{t}Y_{t}+\sigma _{t}Z_{t})dt.
\label{eq:duality}
\end{equation}%
Combining \eqref{eq:first-variation} and \eqref{eq:duality},
\begin{equation}
DJ_{\lambda }(\bar{u})[v-\bar{u}]=\mathbb{E}\int_{0}^{T}(v_{t}-\bar{u}_{t})^{%
\mathsf{T}}G_{t}dt,\qquad G_{t}:=B_{t}Y_{t}+\sigma _{t}Z_{t}.
\label{eq:gradient}
\end{equation}

\subsection{Maximum principle and the portfolio FBSDVI}

\begin{theorem}
\label{the1}
Under assumption \textup{(A1)}, consider the linear controlled wealth
process (\ref{eq:wealth}) with the cost functional (\ref{cost2}).
A control \(\bar u\in\mathcal U\) minimizes \eqref{cost2} over
\[
\mathcal U_{ad}
=
\left\{
u\in\mathbb H^2(t,T;\mathbb R^m):
u_s\in K_s\ \text{for }dt\otimes d\mathbb P\text{-a.e.}
\right\}
\]
if and only if its wealth process \(\bar X\), together with an adjoint pair
\((Y,Z)\), satisfies \eqref{eq:wealth}, the linear mean-field BSDE
\begin{equation}
\label{eq:adjoint}
dY_s
=
-r_sY_s\,ds
+
Z_s^\top dW_s
=
-r_sY_s\,ds
+
\sum_{j=1}^{d}Z_s^j\,dW_s^j,
\qquad
Y_T
=
\lambda\bigl(\bar X_T-\mathbb E[\bar X_T]\bigr)-1,
\end{equation}
and
\begin{equation}
\label{eq:VI}
\bar u_s\in K_s,
\qquad
\left\langle G_s,v-\bar u_s\right\rangle\ge0
\quad\text{for every }v\in K_s,
\qquad
dt\otimes d\mathbb P\text{-a.e.},
\end{equation}
where
\begin{equation}
\label{eq:gradient-process}
G_s
:=
B_sY_s+\sigma_sZ_s
=
B_sY_s+\sum_{j=1}^{d}\sigma_{\cdot j,s}Z_s^j.
\end{equation}
\end{theorem}

\begin{proof}
We prove the result through the Lions derivative, a linearized state
equation, an adjoint It\^o integration by parts, and a predictable
localization argument.

\paragraph{Step 1: Lions derivative of the mean-field terminal cost.}
For \(\mu\in\mathcal P_2(\mathbb R)\), define
\[
\Phi(\mu)
:=
\frac{\lambda}{2}\operatorname{Var}(\mu)-m(\mu),
\qquad
m(\mu):=\int_{\mathbb R}x\,\mu(dx).
\]
Then
$
J(u)=\Phi(\mathbb P_{X_T^u}).
$
The lift of \(\Phi\) to \(L^2(\mathcal F_T;\mathbb R)\) is
\[
\widetilde\Phi(\xi)
=
\frac{\lambda}{2}
\left(
\mathbb E[\xi^2]-(\mathbb E[\xi])^2
\right)
-\mathbb E[\xi].
\]
For \(h\in L^2(\mathcal F_T;\mathbb R)\),
\[
\begin{aligned}
D\widetilde\Phi(\xi)[h]
&=
\lambda\left(
\mathbb E[\xi h]-\mathbb E[\xi]\mathbb E[h]
\right)-\mathbb E[h]\\
&=
\mathbb E\left[
\bigl(\lambda(\xi-\mathbb E[\xi])-1\bigr)h
\right].
\end{aligned}
\]
Hence
\begin{equation}
\label{eq:Lions-derivative}
\partial_\mu\Phi(\mu)(x)
=
\lambda\bigl(x-m(\mu)\bigr)-1.
\end{equation}
Evaluating this Lions derivative at \(\mu=\mathbb P_{\bar X_T}\) and \(x=\bar X_T\)
gives exactly the terminal condition
\[
Y_T
=
\partial_\mu\Phi(\mathbb P_{\bar X_T})(\bar X_T)
=
\lambda\bigl(\bar X_T-\mathbb E[\bar X_T]\bigr)-1.
\]
Thus the mean-field dependence enters the adjoint equation through the
terminal condition.

\paragraph{Step 2: Linearized wealth equation and directional derivative.}
Let \(q\) be a control increment and set \(u^\theta=\bar u+\theta q\). The
linearity of \eqref{eq:wealth} implies
$
X_s^{u^\theta}=\bar X_s^u+\theta\widehat X_s^q,
$
where \(\widehat X^q\) solves
\begin{equation}
\label{eq:linearized-state}
d\widehat X_s^q
=
\bigl(A_s\widehat X_s^q+B_s^\top q_s\bigr)\,ds
+
q_s^\top\sigma_s\,dW_s,
\qquad
\widehat X_t^q=0.
\end{equation}
Writing \(h=\widehat X_T^q\), one has
\[
\operatorname{Var}(X_T+\theta h)
=
\operatorname{Var}(X_T)
+
2\theta\,\mathbb E\bigl[(X_T-\mathbb E[X_T])h\bigr]
+
\theta^2\operatorname{Var}(h).
\]
Therefore
\begin{equation}
\label{eq:exact-cost-expansion}
J(u+\theta q)
=
J(u)
+
\theta\,\mathbb E[Y_T\widehat X_T^q]
+
\frac{\lambda\theta^2}{2}
\operatorname{Var}(\widehat X_T^q).
\end{equation}
It follows that
\begin{equation}
\label{eq:directional-derivative}
DJ(u)[q]
=
\mathbb E[Y_T\widehat X_T^q].
\end{equation}
Since the final term in \eqref{eq:exact-cost-expansion} is nonnegative,
\(J\) is convex in \(u\).

\paragraph{Step 3: Adjoint It\^o identity.}
Apply It\^o's product formula to \(Y_s\widehat X_s^q\). Using
\eqref{eq:adjoint} and \eqref{eq:linearized-state}, we obtain
\[
\begin{aligned}
d\bigl(Y_s\widehat X_s^q\bigr)
={}&
\widehat X_s^q\,dY_s
+
Y_s\,d\widehat X_s^q
+
d[Y,\widehat X^q]_s\\
={}&
-A_sY_s\widehat X_s^q\,ds
+
\widehat X_s^qZ_s^\top dW_s\\
&+
Y_s\bigl(A_s\widehat X_s^q+B_s^\top q_s\bigr)\,ds
+
Y_sq_s^\top\sigma_s\,dW_s\\
&+
Z_s^\top\sigma_s^\top q_s\,ds.
\end{aligned}
\]
The terms containing \(A_sY_s\widehat X_s^q\) cancel. Moreover, $
Y_sB_s^\top q_s
=
(B_sY_s)^\top q_s,
$ and, by the definition of the matrix-vector product, $
Z_s^\top\sigma_s^\top q_s
=
(\sigma_sZ_s)^\top q_s
=
\sum_{j=1}^{d}
(\sigma_{\cdot j,s}Z_s^j)^\top q_s.
$ Consequently,
\[
d\bigl(Y_s\widehat X_s^q\bigr)
=
\langle G_s,q_s\rangle\,ds+dM_s,
\]
where $
G_s=B_sY_s+\sigma_sZ_s
$ and $
dM_s
=
\left(
\widehat X_s^qZ_s^\top
+
Y_sq_s^\top\sigma_s
\right)dW_s
$ is a square-integrable martingale increment under \textup{(A1)}.

Integrating from \(t\) to \(T\), taking expectations, and using
\(\widehat X_t^q=0\) and \(\mathbb E[M_T-M_t]=0\), we get
\begin{equation}
\label{eq:gradient}
DJ(u)[q]
=
\mathbb E[Y_T\widehat X_T^q]
=
\mathbb E\int_t^T
\langle G_s,q_s\rangle\,ds.
\end{equation}
Thus the control gradient in the \(k\)-th coordinate is explicitly
\[
G_s^k
=
B_s^kY_s
+
\sum_{j=1}^{d}\sigma_s^{kj}Z_s^j,
\qquad k=1,\ldots,m.
\]

\paragraph{Step 4: Necessity and pointwise variational inequality.}
Suppose that \(\bar u\) minimizes \eqref{cost2}. For any \(v\in\mathcal U\)
and \(\theta\in(0,1]\), convexity of the feasible set gives $
\bar u+\theta(v-\bar u)\in\mathcal U.
$ Optimality and \eqref{eq:gradient} therefore imply
\begin{equation}
\label{eq:integral-VI}
\mathbb E\int_t^T
\langle G_s,v_s-\bar u_s\rangle\,ds\ge0
\qquad
\text{for every }v\in\mathcal U.
\end{equation}
To obtain the pointwise inequality, let \(A\subseteq[t,T]\times\Omega\) be
an arbitrary predictable set and let \(w\) be a predictable process with
\(w_s\in K_s\). Define
\[
v_s=
\begin{cases}
w_s, &(s,\omega)\in A,\\
\bar u_s, &(s,\omega)\notin A.
\end{cases}
\]
Since the box constraints are decomposable, \(v\in\mathcal U\). Substitution
into \eqref{eq:integral-VI} gives
$
\mathbb E\int_A
\langle G_s,w_s-\bar u_s\rangle\,ds\ge0.
$
Because the predictable set \(A\) is arbitrary, the integrand must be
nonnegative \(dt\otimes d\mathbb P\)-almost everywhere. Otherwise, a
measurable selector and localization to a sufficiently small predictable set
would produce a strictly negative integral, contradicting the preceding
inequality. Hence
\[
\langle G_s,v-\bar u_s\rangle\ge0
\qquad
\text{for every }v\in K_s,
\quad
dt\otimes d\mathbb P\text{-a.e.}
\]

\paragraph{Step 5: Sufficiency.}
Conversely, suppose that \(\bar X\), \((Y,Z)\), and \(\bar u\) satisfy
\eqref{eq:wealth}, \eqref{eq:adjoint}, and \eqref{eq:VI}. For any
\(v\in\mathcal U\), set \(q=v-\bar u\). Integrating \eqref{eq:VI} and using
\eqref{eq:gradient} yields
\[
DJ(\bar u)[v-\bar u]
=
\mathbb E\int_t^T
\langle G_s,v_s-\bar u_s\rangle\,ds
\ge0.
\]
Taking \(\theta=1\) in \eqref{eq:exact-cost-expansion} gives
\[
J(v)-J(\bar u)
=
DJ(\bar u)[v-\bar u]
+
\frac{\lambda}{2}
\operatorname{Var}\bigl(\widehat X_T^{v-\bar u}\bigr).
\]
Both terms on the right-hand side are nonnegative, and therefore $
J(v)\ge J(\bar u)
\text{ for every }v\in\mathcal U.
$
Thus \(\bar u\) is optimal, completing the proof.
\end{proof}

Now we are able to present the complete pre-commitment FBSDVI (omitting the superscript $\bar X$, etc)
\begin{equation*}
\begin{cases}
dX_{t}=\bigl(r_{t}X_{t}+B_{t}^{\top }u_{t}\bigr)\,dt+u_{t}^{\top }\sigma
_{t}\,dW_{t}, \\
dY_{t}=-r_{t}Y_{t}\,dt+Z_{t}^{\top }dW_{t}, \\
Y_{T}=\lambda (X_{T}-\mathbb{E}X_{T})-1, \\
\langle B_{t}Y_{t}+\sigma _{t}Z_{t},v_{t}-u_{t}\rangle \geq 0,\quad \forall
v_{t}\in K_{t}.%
\end{cases}%
\end{equation*}

\noindent The inequality $\langle G,v-u\rangle\ge0$ is the first-order
optimality condition for a differentiable convex function in the special case
where the gradient is identically equal to $G$ (see Lemma \ref{firor}). The objective whose gradient is
identically equal to $G$ is precisely the linear function
$f(u)=\langle G,u\rangle$. Hence satisfying the variational inequality is
equivalent to minimizing $\langle G,u\rangle$ over $K$, and this is in turn
equivalent to the projection condition
$
u=\Pi_K(u-G).
$

The variational inequality in the pre-commitment
FBSDVI only characterizes $u$ as an optimal projection. The control-to-state map has a nontrivial null space (redundant assets),
moving $u$ along this null space leaves both the state and the adjoint
processes unchanged while preserving the variational inequality, so the
optimal control is an affine subspace rather than a single point and is
therefore not unique. The optimal control is unique only when the
control-to-state map is injective (no redundant assets) or when a Tikhonov
term with $\varepsilon>0$ is added, which restores strict convexity. Hence the
variational inequality describes the whole optimal solution set, not a single
optimal control.
\subsection{Equivalent algebraic forms}

For a closed convex set $K$, use
\begin{equation*}
N_K(u):=\{\eta:\left\langle \eta,v-u\right\rangle\le0, \ \forall v\in K\}.
\end{equation*}

%
%
%
%
%
\begin{proposition}
\label{prop:equivalence}
Fix \((t,\omega)\), write
\[
G:=B_tY_t+\sigma_tZ_t,
\qquad
K_t:=[\ell_t,h_t]
=\prod_{k=1}^{m}[\ell_t^k,h_t^k],
\]
and assume \(\ell_t^k\le h_t^k\). Then, for \(u\in K_t\), the following
statements are equivalent:

\begin{enumerate}
\item[(i)] The variational inequality
\begin{equation}
\label{eq:VI-box}
\langle G,v-u\rangle\ge 0
\qquad\text{for every }v\in K_t;
\end{equation}

\item[(ii)] The normal-cone inclusion
\begin{equation}
\label{eq:normal-inclusion}
0\in G+N_{K_t}(u);
\end{equation}

\item[(iii)] There exist multipliers
\(\alpha,\beta\in\mathbb R_+^m\) such that
\begin{equation}
\label{eq:KKT}
G-\alpha+\beta=0,
\end{equation}
with the complementarity conditions
\begin{equation}
\label{eq:complementarity}
0\le \alpha\perp u-\ell_t\ge0,
\qquad
0\le \beta\perp h_t-u\ge0;
\end{equation}

\item[(iv)] Componentwise,
\begin{equation}
\label{eq:switching}
u^k=
\begin{cases}
\ell_t^k, & G^k>0,\\[2mm]
h_t^k, & G^k<0,\\[2mm]
\text{any element of }[\ell_t^k,h_t^k], & G^k=0.
\end{cases}
\end{equation}
\end{enumerate}
\end{proposition}

\begin{proof}
We prove the equivalences step by step.

\paragraph{Step 1: Equivalence of \textup{(i)} and \textup{(ii)}.}
Recall that the Euclidean normal cone of a nonempty closed convex set
\(K\subset\mathbb R^m\) at \(u\in K\) is
\begin{equation}
\label{eq:normal-cone-def}
N_K(u)
=
\left\{
n\in\mathbb R^m:
\langle n,v-u\rangle\le0
\ \text{for every }v\in K
\right\}.
\end{equation}
Thus the variational inequality \eqref{eq:VI-box} is equivalent to $
\langle -G,v-u\rangle\le0
\qquad\text{for every }v\in K_t.
$
By the definition of the normal cone, this means $
-G\in N_{K_t}(u).
$
Adding \(G\) to both sides gives $
0\in G+N_{K_t}(u),
$
which is precisely \eqref{eq:normal-inclusion}. Hence $
\text{(i)}\iff\text{(ii)}.
$

\paragraph{Step 2: The normal cone of a box.}
Because \(K_t\) is a Cartesian product of intervals, $
K_t=\prod_{k=1}^{m}[\ell_t^k,h_t^k],
$
its normal cone is the Cartesian product of the one-dimensional normal
cones:
\[
N_{K_t}(u)
=
\prod_{k=1}^{m}N_{[\ell_t^k,h_t^k]}(u^k).
\]
For a scalar interval \([\ell,h]\) with \(\ell<h\),
\[
N_{[\ell,h]}(x)
=
\begin{cases}
(-\infty,0], & x=\ell,\\
\{0\}, & \ell<x<h,\\
[0,\infty), & x=h.
\end{cases}
\]
Indeed:
\begin{itemize}
\item at the lower endpoint \(x=\ell\), the outward normal points in the
negative direction;
\item in the interior, every nonzero direction is feasible to first order,
so the normal cone is \(\{0\}\);
\item at the upper endpoint \(x=h\), the outward normal points in the
positive direction.
\end{itemize}

Consequently, a vector \(n\) belongs to \(N_{K_t}(u)\) if and only if there
exist \(\alpha,\beta\in\mathbb R_+^m\) such that $
\label{eq:normal-cone-multipliers}
n=-\alpha+\beta,
$
and
\begin{equation}
\label{eq:box-complementarity}
\alpha^k(u^k-\ell_t^k)=0,
\qquad
\beta^k(h_t^k-u^k)=0,
\qquad k=1,\dots,m.
\end{equation}
Here \(\alpha^k\) is the multiplier associated with the lower constraint
\(u^k\ge\ell_t^k\), while \(\beta^k\) is the multiplier associated with the
upper constraint \(u^k\le h_t^k\).

Using \eqref{eq:normal-cone-multipliers} with \(n=-G\), we obtain $
-G=-\alpha+\beta,
$
or equivalently $
G=\alpha-\beta,
$
which is the same as $
G-\alpha+\beta=0.
$
Together with \eqref{eq:box-complementarity}, this proves $
\text{(ii)}\implies\text{(iii)}.
$

\paragraph{Step 3: The KKT conditions imply the variational inequality.}
Suppose that \(\alpha,\beta\ge0\) satisfy
\eqref{eq:KKT}--\eqref{eq:complementarity}. Then $
G=\alpha-\beta.
$
For any \(v\in K_t\), using
\[
\alpha^\top(u-\ell_t)=0,
\qquad
\beta^\top(h_t-u)=0,
\]
we obtain
\[
\begin{aligned}
\langle G,v-u\rangle
&=
\alpha^\top(v-u)-\beta^\top(v-u)\\
&=
\alpha^\top(v-\ell_t)
+\beta^\top(h_t-v).
\end{aligned}
\]
Since \(v\in K_t\), $
v-\ell_t\ge0,
h_t-v\ge0.
$
Together with \(\alpha,\beta\ge0\), this give $
\langle G,v-u\rangle\ge0.
$
Thus $
\text{(iii)}\implies\text{(i)}.
$

\paragraph{Step 4: Componentwise characterization.}
Because the coordinates of \(v\) can be varied independently, the
variational inequality \eqref{eq:VI-box} is equivalent to the \(m\) scalar
inequalities
\begin{equation}
\label{eq:scalar-VI}
G^k(v^k-u^k)\ge0
\qquad
\text{for every }v^k\in[\ell_t^k,h_t^k],
\quad k=1,\dots,m.
\end{equation}
Indeed, \eqref{eq:scalar-VI} follows from \eqref{eq:VI-box} by varying only
the \(k\)-th coordinate and keeping all other coordinates fixed.
Conversely, summing \eqref{eq:scalar-VI} over \(k\) gives
\eqref{eq:VI-box}.

Fix \(k\) and consider the scalar interval
\([\ell_t^k,h_t^k]\).
\begin{itemize}
\item If \(G^k>0\), choose \(v^k=\ell_t^k\) in \eqref{eq:scalar-VI}. Then $
G^k(\ell_t^k-u^k)\ge0.
$
Since \(G^k>0\), this implies \(u^k\le\ell_t^k\). Because
\(u^k\ge\ell_t^k\), it follows that $
u^k=\ell_t^k.
$

\item If \(G^k<0\), choose \(v^k=h_t^k\). Then $
G^k(h_t^k-u^k)\ge0.
$
Since \(G^k<0\), we must have \(h_t^k-u^k\le0\), and hence $
u^k=h_t^k.
$

\item If \(G^k=0\), the scalar inequality in \eqref{eq:scalar-VI} becomes
\(0\ge0\), which imposes no restriction on \(u^k\). Therefore $
u^k\in[\ell_t^k,h_t^k].
$
\end{itemize}
Conversely, suppose that \eqref{eq:switching} holds.
\begin{itemize}
\item If \(G^k>0\) and \(u^k=\ell_t^k\), then \(v^k-u^k\ge0\), and hence
\(G^k(v^k-u^k)\ge0\).
\item If \(G^k<0\) and \(u^k=h_t^k\), then \(v^k-u^k\le0\), and hence
again \(G^k(v^k-u^k)\ge0\).
\item If \(G^k=0\), the inequality is identically zero.
\end{itemize}
Thus \eqref{eq:scalar-VI} holds for every \(k\), proving $
\text{(i)}\iff\text{(iv)}.
$

\paragraph{Step 5: Explicit multipliers and conclusion.}
In the nondegenerate case \(\ell_t^k<h_t^k\), the multipliers may be chosen
as
\[
\alpha^k=
\begin{cases}
G^k, & u^k=\ell_t^k\ \text{and }G^k\ge0,\\
0, & \text{otherwise},
\end{cases}
\]
and
\[
\beta^k=
\begin{cases}
-G^k, & u^k=h_t^k\ \text{and }G^k\le0,\\
0, & \text{otherwise}.
\end{cases}
\]
These multipliers satisfy \eqref{eq:KKT}--\eqref{eq:complementarity}. If
\(\ell_t^k=h_t^k\), the feasible set is the singleton
\(\{\ell_t^k\}\), and the corresponding coordinate condition is automatic.

We have therefore established
$
\text{(i)}\iff\text{(ii)}\iff\text{(iii)}\iff\text{(iv)}.
$
This completes the proof.
\end{proof}
\begin{remark}
\label{rem:singular} Equation \eqref{eq:switching} forces a bound only when $%
G^k\ne0$. On $\{G^k=0\}$, the first-order system is singular and does not
determine $u^k$. Interior portfolios in classical mean--variance solutions
are consistent with this condition: the coupled variables $(Y,Z)$ adjust so
that the switching vector vanishes. Calling the entire optimizer
``bang--bang'' without excluding positive-measure switching sets is
therefore too strong.
\end{remark}

\section{Regularization, response map, and reduced FBSDEs}

For $\varepsilon >0$, consider
\begin{equation}
J_{\lambda ,\varepsilon }(u)=\mathbb{E}\bigg [\frac{\lambda }{2}\left(
X_{T}^{2}-\left( \mathbb{E}X_{T}\right) ^{2}\right) -X_{T}\bigg ]+\frac{%
\varepsilon }{2}\mathbb{E}\int_{0}^{T}|u_{t}|^{2}dt.
\label{prob:regularized}
\end{equation}

\begin{lemma}
\label{lem:cost-convexity}
Let the wealth process satisfy
\begin{equation}
dX_t=\bigl(r_tX_t+B_t^\top u_t\bigr)\,dt+u_t^\top\sigma_t\,dW_t,
\qquad X_0=x_0.
\label{eq:wealth-SDE-lemma}
\end{equation}
For \(\lambda>0\) and \(\varepsilon\ge 0\), define
\begin{equation}
J_{\lambda,\varepsilon}(u)
=\mathbb{E}\!\left[\frac{\lambda}{2}
\Big(X_T^2-(\mathbb{E}X_T)^2\Big)-X_T\right]
+\frac{\varepsilon}{2}\mathbb{E}\!\int_0^T|u_t|^2\,dt.
\label{eq:cost-lemma}
\end{equation}
Then \(J_{\lambda,\varepsilon}\) is convex in the control \(u\), and strictly
convex when \(\varepsilon>0\).
\end{lemma}

\begin{proof}
We give the proof into four steps:
\paragraph{Step 1: The control-to-state map is affine.}
Equation \eqref{eq:wealth-SDE-lemma} is a linear SDE in \(u\): its
coefficients \(r_t,B_t,\sigma_t\) do not depend on the control. Let
\(\Phi_{t,s}\) denote the adjoint (discount) factor, i.e.
\(d\Phi_{t,s}=r_t\Phi_{t,s}\,dt\) in the variable \(t\) with
\(\Phi_{s,s}=1\) (when \(r\) is deterministic,
\(\Phi_{t,s}=\exp\!\int_s^t r_\tau\,d\tau\)). The variation-of-constants
formula gives
\begin{equation}
X_t=\Phi_{t,0}x_0
+\int_0^t\Phi_{t,s}B_s^\top u_s\,ds
+\int_0^t\Phi_{t,s}u_s^\top\sigma_s\,dW_s.
\label{eq:state-affine}
\end{equation}
Hence \(u\mapsto X_T\) is affine: \(X_T(u)=a+Lu\), where
\(a=\Phi_{T,0}x_0\) is independent of \(u\) and \(L\) is a linear operator.

\paragraph{Step 2: Rewrite the terminal term as a variance.}
Since \((\mathbb{E}X_T)^2\) is a deterministic constant,
\[
\mathbb{E}\!\left[X_T^2-(\mathbb{E}X_T)^2\right]
=\mathbb{E}[X_T^2]-(\mathbb{E}X_T)^2
=\operatorname{Var}(X_T).
\]
Thus
\begin{equation}
J_{\lambda,\varepsilon}(u)
=\frac{\lambda}{2}\operatorname{Var}\bigl(X_T(u)\bigr)
-\mathbb{E}[X_T(u)]
+\frac{\varepsilon}{2}\mathbb{E}\!\int_0^T|u_t|^2\,dt.
\label{eq:cost-split}
\end{equation}

\paragraph{Step 3: The variance term is convex.}
Take two controls \(u^{(1)},u^{(2)}\), set \(X_i=X_T(u^{(i)})\), and let
\(\theta\in[0,1]\). By affinity,
\[
X_T\bigl(\theta u^{(1)}+(1-\theta)u^{(2)}\bigr)
=\theta X_1+(1-\theta)X_2.
\]
Writing \(v_i=\operatorname{Var}(X_i)\) and
\(c=\operatorname{Cov}(X_1,X_2)\),
\[
\operatorname{Var}(\theta X_1+(1-\theta)X_2)
=\theta^2v_1+(1-\theta)^2v_2+2\theta(1-\theta)c.
\]
The Cauchy--Schwarz inequality for covariances gives
\(c\le\sqrt{v_1v_2}\), so
\[
\operatorname{Var}(\theta X_1+(1-\theta)X_2)
\le\bigl(\theta\sqrt{v_1}+(1-\theta)\sqrt{v_2}\bigr)^2
\le\theta v_1+(1-\theta)v_2,
\]
where the last inequality follows from the concavity of \(z\mapsto\sqrt z\).
Therefore \(u\mapsto\operatorname{Var}(X_T(u))\) is convex, and multiplication
by \(\lambda/2>0\) preserves convexity.

\paragraph{Step 4: The remaining terms, and strict convexity.}
From \eqref{eq:state-affine}, \(\mathbb{E}[X_T(u)]\) is affine in \(u\); hence
\(-\mathbb{E}[X_T]\) is affine and therefore convex. The map
\(u\mapsto |u|^2\) is convex, so
\(\frac{\varepsilon}{2}\mathbb{E}\!\int_0^T|u_t|^2\,dt\) is convex for
\(\varepsilon\ge0\), and strictly convex for \(\varepsilon>0\).

A sum of convex functions is convex; moreover, adding a strictly convex
function to a convex one yields a strictly convex function. Hence
\(J_{\lambda,\varepsilon}\) is convex in \(u\), and strictly convex for
\(\varepsilon>0\). This completes the proof.
\end{proof}

The first-order VI becomes
\begin{equation}
u_{t}\in K_{t},\qquad \left\langle \varepsilon u_{t}+B_{t}Y_{t}+\sigma
_{t}Z_{t},v-u_{t}\right\rangle \geq 0\quad (v\in K_{t}).
\label{eq:regularized-VI}
\end{equation}

\begin{proposition}
\label{prop:response} The VI \eqref{eq:regularized-VI} has the unique
response
\begin{equation}  \label{eq:projection}
u_t=\Phi_{\varepsilon,t}(Y_t,Z_t) :=\func{Proj}_{K_t}\left[%
-\varepsilon^{-1}(B_tY_t+\sigma_tZ_t)\right].
\end{equation}
It is predictable and satisfies
\begin{equation}  \label{eq:response-Lip}
|\Phi_{\varepsilon,t}(y,z)-\Phi_{\varepsilon,t}(y^{\prime },z^{\prime })| \le%
\frac{\left\| B_t\right\||y-y^{\prime }|+\left\|
\sigma_t\right\||z-z^{\prime }|}{\varepsilon}.
\end{equation}
For a box, \eqref{eq:projection} is componentwise clipping.
\end{proposition}

\begin{proof}
We prove the assertion pointwise in $(t,\omega)$. We therefore fix $%
(t,\omega)$ outside a null set and suppress the dependence on $(t,\omega)$.

\paragraph{Step 1: Variational inequality as a quadratic minimization.}

Set $
g:=By+\sigma z
$
and define the quadratic functional
\begin{equation*}
Q(u):=\frac{\varepsilon}{2}|u|^2+g^\top u, \qquad u\in K.
\end{equation*}
For any $u,v\in\mathbb{R}^m$, a direct expansion gives
\begin{equation}  \label{eq:quadratic-expansion}
Q(v)-Q(u) = \langle \varepsilon u+g,v-u\rangle + \frac{\varepsilon}{2}%
|v-u|^2.
\end{equation}
Since $\varepsilon>0$, the functional $Q$ is strongly convex on $\mathbb{R}^m
$, with modulus $\varepsilon$.

If $u\in K$ satisfies the variational inequality
\begin{equation}  \label{eq:pointwise-VI}
\langle \varepsilon u+g,v-u\rangle\ge0, \qquad \forall v\in K,
\end{equation}
then \eqref{eq:quadratic-expansion} implies
\begin{equation*}
Q(v)-Q(u)\ge \frac{\varepsilon}{2}|v-u|^2\ge0, \qquad \forall v\in K.
\end{equation*}
Thus $u$ minimizes $Q$ over $K$. Conversely, if $u$ minimizes $Q$ over $K$,
the first-order optimality condition for the convex differentiable
functional $Q$ is precisely \eqref{eq:pointwise-VI}. Hence solving %
\eqref{eq:regularized-VI} is equivalent to minimizing $Q$ over $K$.

Strong convexity also implies uniqueness. Indeed, if $u_1$ and $u_2$ are two
solutions, applying \eqref{eq:quadratic-expansion} in both directions gives
\begin{equation*}
Q(u_2)-Q(u_1)\ge \frac{\varepsilon}{2}|u_2-u_1|^2,
\end{equation*}
and
\begin{equation*}
Q(u_1)-Q(u_2)\ge \frac{\varepsilon}{2}|u_1-u_2|^2.
\end{equation*}
Adding these inequalities yields $
0\ge \varepsilon |u_1-u_2|^2,
$
so $u_1=u_2$.

\paragraph{Step 2: Completion of the square.}

We complete the square:
\begin{equation}  \label{eq:complete-square}
\begin{aligned} Q(u) &= \frac{\varepsilon}{2}|u|^2+g^\top u\\ &=
\frac{\varepsilon}{2} \left|u+\varepsilon^{-1}g\right|^2 -
\frac{1}{2\varepsilon}|g|^2. \end{aligned}
\end{equation}
The second term on the right-hand side is independent of $u$. Therefore
minimizing $Q(u)$ over $K$ is equivalent to minimizing $
\left|u+\varepsilon^{-1}g\right|^2
$
over $K$, or equivalently, to finding the Euclidean projection of $%
-\varepsilon^{-1}g$ onto $K$. Thus
\begin{equation*}
u = \func{Proj}_{K}\left(-\varepsilon^{-1}g\right) = \func{Proj}_{K} \left[%
-\varepsilon^{-1}(By+\sigma z)\right].
\end{equation*}
For completeness, the Euclidean projection $p=P_K(a)$ is characterized by
the variational inequality
\begin{equation}  \label{eq:projection-VI}
\langle p-a,v-p\rangle\ge0, \qquad \forall v\in K.
\end{equation}
Taking $a=-\varepsilon^{-1}g,$
condition \eqref{eq:projection-VI} becomes
\begin{equation*}
\left\langle p+\varepsilon^{-1}g,v-p\right\rangle\ge0, \qquad \forall v\in
K.
\end{equation*}
Multiplying by $\varepsilon>0$, we obtain
\begin{equation*}
\langle \varepsilon p+g,v-p\rangle\ge0, \qquad \forall v\in K,
\end{equation*}
which is exactly \eqref{eq:pointwise-VI}. This proves \eqref{eq:projection}.

%

\paragraph{Step 3: Lipschitz estimate.}

For two pairs $(y,z)$ and $(y^{\prime },z^{\prime })$, define
\begin{equation*}
a:=-\varepsilon^{-1}(By+\sigma z), \qquad a^{\prime -1}(By^{\prime }+\sigma
z^{\prime }).
\end{equation*}
Then
\begin{equation*}
\Phi_{\varepsilon,t}(y,z)=P_K(a), \qquad \Phi_{\varepsilon,t}(y^{\prime
},z^{\prime })=P_K(a^{\prime }).
\end{equation*}
By Proposition \ref{pro2}, the triangle inequality,
the submultiplicativity of the Euclidean matrix norm, and the fact that $y$
is scalar, we get
\begin{equation*}
\begin{aligned}
\left|\Phi_{\varepsilon,t}(y,z)
-\Phi_{\varepsilon,t}(y',z')\right|
&\le |a-a'|\\
&=
\frac{1}{\varepsilon}
\left|B(y-y')+\sigma(z-z')\right|\\
&\le
\frac{1}{\varepsilon}
\left(
|B|\,|y-y'|
+
\|\sigma\|\,|z-z'|
\right).
\end{aligned}
\end{equation*}
This proves \eqref{eq:response-Lip}.

\paragraph{Step 4: Predictability.}

The coefficients $B$ and $\sigma$, the adjoint process $Y$, and a
predictable version of the martingale coefficient $Z$ are predictable. Hence $
G_t=B_tY_t+\sigma_tZ_t
$
is predictable. Euclidean projection is continuous, and therefore Borel
measurable, as a function of the point being projected. For a predictable
closed-convex-valued multifunction $K_t$, the measurable projection theorem
implies that the unique projection $\func{Proj}_{K_t}(\cdot)$ is jointly
measurable in $(t,\omega)$. Consequently, the composition
\begin{equation*}
(t,\omega)\longmapsto \func{Proj}_{K_t(\omega)} \left[-\varepsilon^{-1}G_t(%
\omega)\right]
\end{equation*}
is predictable. For the box constraint below, this also follows directly
from the explicit min--max formula.

\paragraph{Step 5: The box case.}

Suppose
\begin{equation*}
K_t = \prod_{i=1}^{m}[\ell_t^i,h_t^i], \qquad \ell_t^i\le h_t^i.
\end{equation*}
Projection onto a Cartesian product separates coordinatewise: $
\left[P_{K_t}(a)\right]^i = P_{[\ell_t^i,h_t^i]}(a^i).
$
For a one-dimensional interval $[\ell,h]$, $
P_{[\ell,h]}(x) = \max\bigl\{\ell,\min\{h,x\}\bigr\}.
$
Therefore, writing $
a_t = -\varepsilon^{-1}\bigl(B_tY_t+\sigma_tZ_t\bigr),
$
the regularized control is
\begin{equation}  \label{eq:box-clip}
\left[\Phi_{\varepsilon,t}(Y_t,Z_t)\right]^i = \max\left\{ \ell_t^i,
\min\left\{ h_t^i, -\varepsilon^{-1} \left[B_tY_t+\sigma_tZ_t\right]^i
\right\} \right\},
\end{equation}
for $i=1,\dots,m$. Thus the projection is precisely componentwise clipping.

Finally, when the box bounds are uniformly bounded in $(t,\omega)$, the
projected control takes values in a uniformly bounded set, and hence $%
\Phi_{\varepsilon,\cdot}(Y_\cdot,Z_\cdot)\in\mathbb{H}^2$. This completes
the proof.
\end{proof}


Substitution gives the reduced FBSDEs
\begin{equation}  \label{eq:reduced-FBSDE}
\left\{\begin{aligned}
dX_t={}&\left[r_tX_t+B_t^\top\Phi_{\varepsilon,t}(Y_t,Z_t)\right]dt
+\Phi_{\varepsilon,t}(Y_t,Z_t)^\top\sigma_t,dW_t,\\
dY_t={}&-r_tY_tdt+Z_t^\trans dW_t,\\ X_0={}&x_0,\qquad Y_T=\lambda (X_{T}-\mathbb{E}X_{T})-1.
\end{aligned}\right.
\end{equation}

In order to study the existence and uniqueness of FBSDEs (\ref{eq:reduced-FBSDE}), we first give some sharp estimates.

\begin{lemma}
\label{lem:linear-forward-stability}
Let \(r\) be a bounded predictable scalar process, let
\(B\in\mathbb{R}^{m}\) and \(\sigma\in\mathbb{R}^{m\times d}\) be
bounded predictable processes, and let
\(q\in\mathbb{H}^{2}(0,T;\mathbb{R}^{m})\). Suppose that
\(\widetilde{X}\in\mathbb{S}^{2}(0,T;\mathbb{R})\) solves
\begin{equation}
\label{eq:linear-forward-increment-est}
d\widetilde{X}_{t}
=
r_{t}\widetilde{X}_{t}\,dt
+B_{t}^{\top}q_{t}\,dt
+q_{t}^{\top}\sigma_{t}\,dW_{t},
\qquad
\widetilde{X}_{0}=0.
\end{equation}
Set
\[
\bar r=\operatorname*{ess\,sup}_{t,\omega}|r_{t}|,
\qquad
\bar B=\operatorname*{ess\,sup}_{t,\omega}|B_{t}|,
\qquad
\bar\sigma
=
\operatorname*{ess\,sup}_{t,\omega}\|\sigma_{t}\|.
\]
Then
\begin{equation}
\label{eq:forward-S2-estimate}
\|\widetilde{X}\|_{\mathbb{S}^{2}}
\le C_{X}\|q\|_{\mathbb{H}^{2}},
\end{equation}
where one may take
\begin{equation}
\label{eq:explicit-CX}
C_{X}
=
\sqrt{2}\,e^{2\bar rT}
\sqrt{T\bar B^{2}+4\bar\sigma^{2}}.
\end{equation}
In particular,
\begin{equation}
\label{eq:terminal-forward-estimate}
\|\widetilde{X}_{T}\|_{L^{2}}
\le
\|\widetilde{X}\|_{\mathbb{S}^{2}}
\le
C_{X}\|q\|_{\mathbb{H}^{2}}.
\end{equation}
\end{lemma}

\begin{proof}
Define the discount factor
\begin{equation}
\label{eq:discount-factor}
\Gamma_{t}
=
\exp\left(-\int_{0}^{t}r_{s}\,ds\right).
\end{equation}
Then $
d\Gamma_{t}=-r_{t}\Gamma_{t}\,dt.
$ Applying It\^o's product rule to
\(\Gamma_{t}\widetilde{X}_{t}\), we obtain
\[
\begin{aligned}
d(\Gamma_{t}\widetilde{X}_{t})
&=
\Gamma_{t}\,d\widetilde{X}_{t}
+\widetilde{X}_{t}\,d\Gamma_{t}\\
&=
\Gamma_{t}
\left(
r_{t}\widetilde{X}_{t}\,dt
+B_{t}^{\top}q_{t}\,dt
+q_{t}^{\top}\sigma_{t}\,dW_{t}
\right)
-r_{t}\Gamma_{t}\widetilde{X}_{t}\,dt\\
&=
\Gamma_{t}B_{t}^{\top}q_{t}\,dt
+\Gamma_{t}q_{t}^{\top}\sigma_{t}\,dW_{t}.
\end{aligned}
\]
Since \(\widetilde{X}_{0}=0\), integration gives
\begin{equation}
\label{eq:discounted-integral-form}
\Gamma_{t}\widetilde{X}_{t}
=
A_{t}+M_{t},
\end{equation}
where
\[
A_{t}
=
\int_{0}^{t}\Gamma_{s}B_{s}^{\top}q_{s}\,ds,
\qquad
M_{t}
=
\int_{0}^{t}\Gamma_{s}q_{s}^{\top}\sigma_{s}\,dW_{s}.
\]
Because \(|r_{t}|\le \bar r\),
\begin{equation}
\label{eq:discount-bounds}
e^{-\bar rT}\le \Gamma_{t}\le e^{\bar rT},
\qquad
\Gamma_{t}^{-1}\le e^{\bar rT}.
\end{equation}
Write $
R:=e^{\bar rT}.
$

For the finite-variation term, the Cauchy--Schwarz inequality gives
\[
\begin{aligned}
\mathbb{E}\sup_{0\le t\le T}|A_{t}|^{2}
&\le
\mathbb{E}
\left(
\int_{0}^{T}
\Gamma_{s}|B_{s}^{\top}q_{s}|\,ds
\right)^{2}\\
&\le
T\,\mathbb{E}\int_{0}^{T}
\Gamma_{s}^{2}|B_{s}^{\top}q_{s}|^{2}\,ds.
\end{aligned}
\]
Using \(\Gamma_{s}\le R\), \(|B_{s}|\le\bar B\), and
\(|B_{s}^{\top}q_{s}|\le |B_{s}||q_{s}|\), we get
\begin{equation}
\label{eq:drift-estimate}
\mathbb{E}\sup_{0\le t\le T}|A_{t}|^{2}
\le
R^{2}T\bar B^{2}
\mathbb{E}\int_{0}^{T}|q_{s}|^{2}\,ds
=
R^{2}T\bar B^{2}\|q\|_{\mathbb{H}^{2}}^{2}.
\end{equation}
The process \(M\) is a continuous square-integrable martingale because
\[
\mathbb{E}\int_{0}^{T}
\Gamma_{s}^{2}|q_{s}^{\top}\sigma_{s}|^{2}\,ds
\le
R^{2}\bar\sigma^{2}\|q\|_{\mathbb{H}^{2}}^{2}
<\infty.
\]
Doob's \(L^{2}\) maximal inequality and It\^o's isometry give
\[
\begin{aligned}
\mathbb{E}\sup_{0\le t\le T}|M_{t}|^{2}
&\le
4\mathbb{E}|M_{T}|^{2}\\
&=
4\mathbb{E}\int_{0}^{T}
\Gamma_{s}^{2}|q_{s}^{\top}\sigma_{s}|^{2}\,ds.
\end{aligned}
\]
Since
\[
|q_{s}^{\top}\sigma_{s}|
=
|\sigma_{s}^{\top}q_{s}|
\le
\|\sigma_{s}\||q_{s}|
\le
\bar\sigma|q_{s}|,
\]
it follows that
\begin{equation}
\label{eq:martingale-estimate}
\mathbb{E}\sup_{0\le t\le T}|M_{t}|^{2}
\le
4R^{2}\bar\sigma^{2}\|q\|_{\mathbb{H}^{2}}^{2}.
\end{equation}
Using \((a+b)^{2}\le 2a^{2}+2b^{2}\), equations
\eqref{eq:discounted-integral-form},
\eqref{eq:drift-estimate}, and
\eqref{eq:martingale-estimate} imply
\[
\begin{aligned}
\mathbb{E}\sup_{0\le t\le T}
|\Gamma_{t}\widetilde{X}_{t}|^{2}
&\le
2\mathbb{E}\sup_{0\le t\le T}|A_{t}|^{2}
+
2\mathbb{E}\sup_{0\le t\le T}|M_{t}|^{2}\\
&\le
2R^{2}
\left(
T\bar B^{2}+4\bar\sigma^{2}
\right)
\|q\|_{\mathbb{H}^{2}}^{2}.
\end{aligned}
\]
Finally, by \(\Gamma_{t}^{-1}\le R\), we have
\[
\begin{aligned}
\|\widetilde{X}\|_{\mathbb{S}^{2}}^{2}
&=
\mathbb{E}\sup_{0\le t\le T}|\widetilde{X}_{t}|^{2}\\
&\le
R^{2}
\mathbb{E}\sup_{0\le t\le T}
|\Gamma_{t}\widetilde{X}_{t}|^{2}\\
&\le
2R^{4}
\left(
T\bar B^{2}+4\bar\sigma^{2}
\right)
\|q\|_{\mathbb{H}^{2}}^{2}.
\end{aligned}
\]
Taking square roots and recalling \(R=e^{\bar rT}\), we obtain
\[
\|\widetilde{X}\|_{\mathbb{S}^{2}}
\le
\sqrt{2}\,e^{2\bar rT}
\sqrt{T\bar B^{2}+4\bar\sigma^{2}}\,
\|q\|_{\mathbb{H}^{2}}.
\]
This proves \eqref{eq:forward-S2-estimate} with
\(C_{X}\) given by \eqref{eq:explicit-CX}. The terminal estimate
\eqref{eq:terminal-forward-estimate} follows immediately from
\(|\widetilde{X}_{T}|\le \sup_{0\le t\le T}|\widetilde{X}_{t}|\).
\end{proof}

\begin{lemma}
\label{lem:linear-BSDE-stability}
Let \(r\) be a bounded predictable scalar process with $
|r_t|\le \bar r, \text{a.e. }t\in[0,T],
$
and let \(\xi\in L^2(\mathcal F_T;\mathbb R)\). Suppose that
\((Y,Z)\in\mathbb S^2(0,T;\mathbb R)\times\mathbb H^2(0,T;\mathbb R^d)\)
solves
\begin{equation}
\label{eq:generic-linear-BSDE}
dY_t=-r_tY_t\,dt+Z_t^\top dW_t,
\qquad
Y_T=\xi.
\end{equation}
Then
\begin{equation}
\label{eq:generic-BSDE-estimate}
\|Y\|_{\mathbb S^2}+\|Z\|_{\mathbb H^2}
\le
C_Y(T,\bar r)\|\xi\|_{L^2},
\end{equation}
where one may take
\begin{equation}
\label{eq:explicit-CY}
C_Y(T,\bar r)=3e^{2\bar rT}.
\end{equation}
In particular, if
\begin{equation}
\label{eq:tilde-BSDE}
d\widetilde Y_t
=
-r_t\widetilde Y_t\,dt+\widetilde Z_t^\top dW_t,
\qquad
\widetilde Y_T=2\widetilde X_T,
\end{equation}
then
\begin{equation}
\label{eq:BSDE-increment-estimate}
\|\widetilde Y\|_{\mathbb S^2}
+\|\widetilde Z\|_{\mathbb H^2}
\le
2C_Y\|\widetilde X_T\|_{L^2}.
\end{equation}
Combining this with the forward SDE estimate
\(\|\widetilde X_T\|_{L^2}\le C_X\|q\|_{\mathbb H^2}\) gives
\begin{equation}
\label{eq:combined-linear-estimate}
\|\widetilde Y\|_{\mathbb S^2}
+\|\widetilde Z\|_{\mathbb H^2}
\le
2C_XC_Y\|q\|_{\mathbb H^2}.
\end{equation}
\end{lemma}

\begin{proof}
Define the appreciation factor
\begin{equation}
\label{eq:appreciation-factor}
\beta_t
=
\exp\left(\int_0^t r_s\,ds\right).
\end{equation}
Then $
d\beta_t=r_t\beta_t\,dt.
$
Applying It\^o's product rule to \(\beta_tY_t\), we obtain
\[
\begin{aligned}
d(\beta_tY_t)
&=
\beta_t\,dY_t+Y_t\,d\beta_t\\
&=
\beta_t\left(-r_tY_t\,dt+Z_t^\top dW_t\right)
+r_t\beta_tY_t\,dt\\
&=
\beta_tZ_t^\top dW_t.
\end{aligned}
\]
Hence
\begin{equation}
\label{eq:discounted-BSDE-martingale}
N_t:=\beta_tY_t.
\end{equation}
Since $
e^{-\bar rT}\le \beta_t\le e^{\bar rT},
$
the assumptions \(\xi\in L^2(\mathcal F_T)\) and
\(Z\in\mathbb H^2(0,T;\mathbb R^d)\) imply
\[
\mathbb E|\beta_T\xi|^2
\le e^{2\bar rT}\|\xi\|_{L^2}^2<\infty,
\qquad
\mathbb E\int_0^T\beta_t^2|Z_t|^2\,dt
\le e^{2\bar rT}\|Z\|_{\mathbb H^2}^2<\infty.
\]
Hence \(N_t=\beta_tY_t\) is a continuous square-integrable martingale,
rather than merely a local martingale.

Integrating from \(t\)
to \(T\) and taking conditional expectation yields
\begin{equation}
\label{eq:conditional-representation}
\beta_tY_t
=
\mathbb E[\beta_T\xi\mid\mathcal F_t].
\end{equation}

\paragraph{Estimate for \(Y\).}
By Doob's \(L^2\) maximal inequality, $
\mathbb E\sup_{0\le t\le T}|N_t|^2
\le
4\mathbb E|N_T|^2.
$
Since \(N_T=\beta_T\xi\) and $
\beta_T\le e^{\bar rT},
$ we have
\[
\mathbb E\sup_{0\le t\le T}|N_t|^2
\le
4e^{2\bar rT}\|\xi\|_{L^2}^2.
\]
Moreover,
\[
\beta_t^{-1}
=
\exp\left(-\int_0^t r_s\,ds\right)
\le e^{\bar rT},
\]
and therefore
\[
\begin{aligned}
\|Y\|_{\mathbb S^2}^2
&=
\mathbb E\sup_{0\le t\le T}|Y_t|^2\\
&\le
e^{2\bar rT}
\mathbb E\sup_{0\le t\le T}|N_t|^2\\
&\le
4e^{4\bar rT}\|\xi\|_{L^2}^2.
\end{aligned}
\]
Thus
\begin{equation}
\label{eq:Y-S2-estimate}
\|Y\|_{\mathbb S^2}
\le
2e^{2\bar rT}\|\xi\|_{L^2}.
\end{equation}

\paragraph{Estimate for \(Z\).}
Since
\[
N_T-N_0
=
\int_0^T\beta_sZ_s^\top dW_s,
\]
It\^o's isometry gives
\[
\mathbb E\int_0^T\beta_s^2|Z_s|^2\,ds
=
\mathbb E|N_T-N_0|^2.
\]
Because \(N\) is a martingale and
\(N_0=\mathbb E[N_T\mid\mathcal F_0]\),
\[
\mathbb E|N_T-N_0|^2
=
\mathbb E|N_T|^2-\mathbb E|N_0|^2
\le
\mathbb E|N_T|^2.
\]
Using \(N_T=\beta_T\xi\) and \(\beta_T\le e^{\bar rT}\), we obtain
\[
\mathbb E\int_0^T\beta_s^2|Z_s|^2\,ds
\le
e^{2\bar rT}\|\xi\|_{L^2}^2.
\]
On the other hand,
\[
\beta_s
=
\exp\left(\int_0^s r_u\,du\right)
\ge e^{-\bar rT},
\]
so \(\beta_s^2\ge e^{-2\bar rT}\). Consequently,
\[
e^{-2\bar rT}\|Z\|_{\mathbb H^2}^2
\le
\mathbb E\int_0^T\beta_s^2|Z_s|^2\,ds
\le
e^{2\bar rT}\|\xi\|_{L^2}^2.
\]
Therefore
\begin{equation}
\label{eq:Z-H2-estimate}
\|Z\|_{\mathbb H^2}
\le
e^{2\bar rT}\|\xi\|_{L^2}.
\end{equation}
Combining \eqref{eq:Y-S2-estimate} and
\eqref{eq:Z-H2-estimate}, we get
\[
\|Y\|_{\mathbb S^2}+\|Z\|_{\mathbb H^2}
\le
3e^{2\bar rT}\|\xi\|_{L^2},
\]
which proves \eqref{eq:generic-BSDE-estimate} with
\(C_Y=3e^{2\bar rT}\).

Finally, for \(\xi=2\widetilde X_T\), equation
\eqref{eq:generic-BSDE-estimate} gives
\[
\|\widetilde Y\|_{\mathbb S^2}
+\|\widetilde Z\|_{\mathbb H^2}
\le
2C_Y\|\widetilde X_T\|_{L^2}.
\]
Using the forward SDE estimate
\[
\|\widetilde X_T\|_{L^2}
\le
\|\widetilde X\|_{\mathbb S^2}
\le
C_X\|q\|_{\mathbb H^2},
\]
we conclude that
\[
\|\widetilde Y\|_{\mathbb S^2}
+\|\widetilde Z\|_{\mathbb H^2}
\le
2C_XC_Y\|q\|_{\mathbb H^2}.
\]
This completes the proof.
\end{proof}

\begin{theorem}
\label{thm:homotopy-FBSDE}
Fix $\varepsilon>0$, $x_0\in\R$, and $\lambda>0$.
Assume that $r$ is a bounded adapted scalar process, that
$B\in\R^m$ and $\sigma\in\R^{m\times d}$ are bounded adapted
processes, and that $K_t(\omega)\subset\R^m$ is a adapted
nonempty closed convex box satisfying
\[
\sup_{t,\omega}\sup_{v\in K_t(\omega)}|v|
\le \kappa_K<\infty.
\]
Define
\[
\Phi_{\varepsilon,t}(y,z)
=
\Proj_{K_t}
\left[
-\varepsilon^{-1}\bigl(B_ty+\sigma_tz\bigr)
\right].
\]
Then the FBSDE
\begin{equation}
\label{eq:regularized-FBSDE}
\begin{cases}
dX_t
=
\left[r_tX_t+B_t^\top\Phi_{\varepsilon,t}(Y_t,Z_t)\right]dt
+
\Phi_{\varepsilon,t}(Y_t,Z_t)^\top\sigma_t\,dW_t,\\[1mm]
dY_t=-r_tY_t\,dt+Z_t^\top dW_t,\\[1mm]
X_0=x_0,\qquad Y_T=\lambda (X_{T}-\mathbb{E}X_{T})-1,
\end{cases}
\end{equation}
has a unique adapted solution
\[
(X,Y,Z)
\in
\mathbb S^2(0,T;\R)
\times
\mathbb S^2(0,T;\R)
\times
\mathbb H^2(0,T;\R^d).
\]
Moreover, the associated control
$u^{\varepsilon}_t=\Phi_{\varepsilon,t}(Y_t,Z_t)$
belongs to $\mathbb H^2(0,T;\R^m)$, takes values in $K_t$ almost
everywhere, and is unique.
\end{theorem}

\begin{proof}
We divide the proof into nine steps.

\paragraph{Step 1: Projection inequalities.}

For a nonempty closed convex set $K$, the Euclidean projection
$\Proj_K$ is characterized by
\[
\langle a-\Proj_K(a),v-\Proj_K(a)\rangle\le 0,
\qquad \forall v\in K.
\]
Consequently, if
$u=\Proj_K(-\varepsilon^{-1}A)$,
then
\begin{equation}
\label{eq:projection-VI}
\langle \varepsilon u+A,v-u\rangle\ge 0,
\qquad \forall v\in K.
\end{equation}
The projection is firmly nonexpansive:
\[
\langle \Proj_K(a)-\Proj_K(a'),a-a'\rangle
\ge
|\Proj_K(a)-\Proj_K(a')|^2.
\]
Thus, if
\[
u=\Proj_K(-\varepsilon^{-1}A),\qquad
u'=\Proj_K(-\varepsilon^{-1}A'),
\qquad q=u-u',
\]
then
\begin{equation}
\label{eq:firm-nonexpansiveness}
\langle q,A-A'\rangle\le -\varepsilon |q|^2.
\end{equation}
For a time-dependent box $K_t=[\ell_t,h_t]$, the projection is the
componentwise clipping map and is therefore predictable whenever
$\ell_t$ and $h_t$ are predictable.

\paragraph{Step 2: Linear SDE--BSDE stability estimate.}

Let $q\in\mathbb H^2(0,T;\R^m)$, and let $\widetilde X$ solve
\begin{equation}
\label{eq:linear-forward-increment}
d\widetilde X_t
=
r_t\widetilde X_t\,dt+B_t^\top q_t\,dt+q_t^\top\sigma_t\,dW_t,
\qquad
\widetilde X_0=0.
\end{equation}
Applying It\^o's formula to the discounted process
\[
\exp\!\left(-\int_0^t r_s\,ds\right)\widetilde X_t
\]
and using the It\^o isometry and the boundedness of the coefficients
gives a constant $C_X=C_X(T,\bar r,\bar B,\bar\sigma)$ such that
\begin{equation}
\label{eq:forward-stability}
\|\widetilde X\|_{\mathbb S^2}
\le C_X\|q\|_{\mathbb H^2},
\end{equation}
where
\[
\bar r=\esssup_{t,\omega}|r_t|,\qquad
\bar B=\esssup_{t,\omega}|B_t|,\qquad
\bar\sigma=\esssup_{t,\omega}\|\sigma_t\|.
\]
Next, let $(\widetilde Y,\widetilde Z)$ solve the linear BSDE
\begin{equation}
\label{eq:linear-backward-increment}
d\widetilde Y_t
=
-r_t\widetilde Y_t\,dt+\widetilde Z_t^\top dW_t,
\qquad
\widetilde Y_T=\xi,
\end{equation}
for a generic terminal condition $\xi\in L^2(\mathcal F_T)$. The
standard linear BSDE estimate (discounting together with the It\^o
isometry) gives a constant $C_Y=C_Y(T,\bar r)$ such that
\begin{equation}
\label{eq:bsde-stability}
\|\widetilde Y\|_{\mathbb S^2}
+
\|\widetilde Z\|_{\mathbb H^2}
\le
C_Y\|\xi\|_{L^2}.
\end{equation}
For the actual difference system the terminal data is the centered
value $\xi=\lambda\Pz\widetilde X_T$. Since $\Pz$ is a contraction,
\[
\|\xi\|_{L^2}
=\lambda\|\Pz\widetilde X_T\|_{L^2}
\le\lambda\|\widetilde X_T\|_{L^2}
\le\lambda C_X\|q\|_{\mathbb H^2}.
\]
Set
$\widetilde G_t=B_t\widetilde Y_t+\sigma_t\widetilde Z_t$.
Then
\begin{equation}
\label{eq:G-stability}
\begin{aligned}
\|\widetilde G\|_{\mathbb H^2}
&\le
\bar B\sqrt T\,\|\widetilde Y\|_{\mathbb S^2}
+\bar\sigma\|\widetilde Z\|_{\mathbb H^2}\\
&\le
\bigl(\bar B\sqrt T+\bar\sigma\bigr)
\bigl(\|\widetilde Y\|_{\mathbb S^2}
      +\|\widetilde Z\|_{\mathbb H^2}\bigr)\\
&\le
L_G\|q\|_{\mathbb H^2},
\end{aligned}
\end{equation}
where
\begin{equation}
\label{eq:LG-definition}
L_G
=
\lambda C_XC_Y\bigl(\sqrt T\,\bar B+\bar\sigma\bigr).
\end{equation}

\paragraph{Step 3: Homotopic family.}

For $\alpha\in[0,1]$ and an arbitrary exogenous source
$\eta\in\mathbb H^2(0,T;\R^m)$, consider the family
\begin{equation}
\label{eq:homotopy-FBSDE}
\begin{cases}
dX_t^{\alpha,\eta}
=
\left[r_tX_t^{\alpha,\eta}
+B_t^\top u_t^{\alpha,\eta}\right]dt
+
(u_t^{\alpha,\eta})^\top\sigma_t\,dW_t,\\[1mm]
dY_t^{\alpha,\eta}
=
-r_tY_t^{\alpha,\eta}dt
+(Z_t^{\alpha,\eta})^\top dW_t,\\[1mm]
u_t^{\alpha,\eta}
=
\Proj_{K_t}
\left[
-\varepsilon^{-1}
\left(
\alpha G_t^{\alpha,\eta}+\eta_t
\right)
\right],\\[1mm]
X_0^{\alpha,\eta}=x_0,\qquad
Y_T^{\alpha,\eta}
=\lambda (X_T^{\alpha,\eta}
          -\mathbb{E}X_T^{\alpha,\eta})-1,
\end{cases}
\end{equation}
where
\[
G_t^{\alpha,\eta}
=
B_tY_t^{\alpha,\eta}
+\sigma_tZ_t^{\alpha,\eta}.
\]
When $\alpha=1$ and $\eta=0$, equation
\eqref{eq:homotopy-FBSDE} is exactly
\eqref{eq:regularized-FBSDE}. When $\alpha=0$, the control is
independent of $Y$ and $Z$:
\begin{equation}
\label{eq:base-control}
u_t^{0,\eta}
=
\Proj_{K_t}[-\varepsilon^{-1}\eta_t].
\end{equation}
Since $K_t$ is uniformly bounded,
$u^{0,\eta}\in\mathbb H^2(0,T;\R^m)$. One first solves the linear
forward SDE with this known control and then solves the linear BSDE
with the centered terminal value
$\lambda\Pz X_T^{0,\eta}-1$. Hence
\eqref{eq:homotopy-FBSDE} has a unique solution for every $\eta$
when $\alpha=0$.

\paragraph{Step 4: Uniqueness for every $\alpha$.}

Fix $\alpha\in[0,1]$ and $\eta$, and suppose that
$(X,Y,Z,u)$ and $(X',Y',Z',u')$ are two solutions. Define
\[
\widehat X=X-X',\qquad
\widehat Y=Y-Y',\qquad
\widehat Z=Z-Z',
\]
\[
q=u-u',\qquad
\widehat G=B\widehat Y+\sigma\widehat Z.
\]
Applying \eqref{eq:firm-nonexpansiveness} with
$
A=\alpha G+\eta, A'=\alpha G'+\eta,
$
gives
\begin{equation}
\label{eq:uniqueness-projection}
\alpha\langle q,\widehat G\rangle
\le -\varepsilon |q|^2.
\end{equation}
The difference processes satisfy
\begin{equation*}
\left\{
\begin{array}{rcl}
d\widehat{X}_{t} & = & r_{t}\widehat{X}_{t}\,dt+B_{t}^{\top
}q_{t}\,dt+q_{t}^{\top }\sigma _{t}\,dW_{t}, \\
d\widehat{Y}_{t} & = & -r_{t}\widehat{Y}_{t}\,dt+\widehat{Z}_{t}^{\top
}dW_{t},\qquad \widehat{Y}_{T}
=\lambda \left( \widehat{X}_{T}-\mathbb{E}\widehat{X}_{T}\right).
\end{array}
\right.
\end{equation*}
Applying It\^o's formula to $\widehat Y_t\widehat X_t$, the two
interest-rate terms cancel and
\[
d(\widehat Y_t\widehat X_t)
=
\langle q_t,\widehat G_t\rangle\,dt+dM_t,
\]
where $M$ is a stochastic integral. Note that
\begin{eqnarray*}
\mathbb{E}\left[ \lambda \left( \widehat{X}_{T}
-\mathbb{E}\widehat{X}_{T}\right) \widehat{X}_{T}\right]
&=&\lambda \mathbb{E}\left[ \left(
\widehat{X}_{T}^{2}
-\widehat{X}_{T}\mathbb{E}\widehat{X}_{T}\right) \right]
\\
&=&\lambda \left( \mathbb{E}\widehat{X}_{T}^{2}
-\left( \mathbb{E}\widehat{X}_{T}\right) ^{2}\right) \\
&=&\lambda \mathbb{E}\left\vert \widehat{X}_{T}
-\mathbb{E}\widehat{X}_{T}\right\vert ^{2}.
\end{eqnarray*}
After the usual localization argument, the expectation of the
martingale term vanishes, so
\begin{equation}
\lambda \mathbb{E}\left\vert \widehat{X}_{T}
-\mathbb{E}\widehat{X}_{T}\right\vert ^{2}
=\mathbb{E}\int_{0}^{T}\langle q_{t},\widehat G_{t}\rangle\,dt.
\label{eq:uniqueness-Ito}
\end{equation}
Combining \eqref{eq:uniqueness-projection} and
\eqref{eq:uniqueness-Ito} yields
\begin{equation}
\alpha \lambda \mathbb{E}\left\vert \widehat{X}_{T}
-\mathbb{E}\widehat{X}_{T}\right\vert ^{2}
+\varepsilon \mathbb{E}\int_{0}^{T}|q_{t}|^{2}\,dt\leq 0.
\label{eq:uniqueness-energy}
\end{equation}
Therefore $q=0$ and, when $\alpha>0$,
$\widehat{X}_{T}=\mathbb{E}\widehat{X}_{T}$. When $\alpha=0$,
\eqref{eq:uniqueness-energy} directly gives $q=0$.

Once $q=0$, the forward difference equation becomes
\[
d\widehat X_t=r_t\widehat X_t\,dt,\qquad \widehat X_0=0,
\]
and hence $\widehat X=0$. The homogeneous linear BSDE
\[
d\widehat Y_t=-r_t\widehat Y_t\,dt+\widehat Z_t^\top dW_t,
\qquad
\widehat Y_T=0
\]
then has the unique solution $\widehat Y=0$, $\widehat Z=0$.
This proves uniqueness for every $\alpha$ and $\eta$.

\paragraph{Step 5: Contraction of the Picard iterates.}

Assume that $\alpha_0\in[0,1)$ is such that
\eqref{eq:homotopy-FBSDE} is solvable for every exogenous source
$\eta$. We prove that the same property holds for
$\alpha_0+\delta$, provided $\delta>0$ is sufficiently small.

Fix $\eta\in\mathbb H^2(0,T;\R^m)$. Let
$U^{(0)}=(X^{(0)},Y^{(0)},Z^{(0)})$ solve the $\alpha_0$-problem
with source $\eta$. Inductively, given
$G^{(k)}=BY^{(k)}+\sigma Z^{(k)}$,
let $U^{(k+1)}$ solve the $\alpha_0$-problem with exogenous
source $\eta+\delta G^{(k)}$.
Equivalently,
\begin{equation}
\label{eq:Picard-projection}
u^{(k+1)}
=
\Proj_{K_t}
\left[
-\varepsilon^{-1}
\left(
\alpha_0G^{(k+1)}
+\eta
+\delta G^{(k)}
\right)
\right].
\end{equation}
Define
\[
q^{(k+1)}=u^{(k+1)}-u^{(k)},
\]
and
\[
\widehat X^{(k+1)}=X^{(k+1)}-X^{(k)},\quad
\widehat Y^{(k+1)}=Y^{(k+1)}-Y^{(k)},\quad
\widehat Z^{(k+1)}=Z^{(k+1)}-Z^{(k)},
\]
\[
\widehat G^{(k+1)}=G^{(k+1)}-G^{(k)}.
\]
The projection variational inequalities for $u^{(k+1)}$ and
$u^{(k)}$ are, respectively,
\[
\left\langle
\varepsilon u^{(k+1)}
+\alpha_0G^{(k+1)}
+\eta
+\delta G^{(k)},
u^{(k)}-u^{(k+1)}
\right\rangle\ge 0,
\]
and
\[
\left\langle
\varepsilon u^{(k)}
+\alpha_0G^{(k)}
+\eta
+\delta G^{(k-1)},
u^{(k+1)}-u^{(k)}
\right\rangle\ge 0.
\]
Adding these inequalities gives
\begin{equation}
\label{eq:iterated-VI}
\alpha_0\langle q^{(k+1)},\widehat G^{(k+1)}\rangle
+
\delta\langle q^{(k+1)},\widehat G^{(k)}\rangle
\le
-\varepsilon |q^{(k+1)}|^2.
\end{equation}
The It\^o cross-energy identity for the increment (the centered
terminal carries no factor $\alpha$) is
\begin{equation}
\lambda \mathbb{E}\left\vert
\Pz\widehat X_T^{(k+1)}\right\vert^{2}
=
\mathbb{E}\int_{0}^{T}
\langle q_t^{(k+1)},\widehat G_t^{(k+1)}\rangle\,dt.
\label{eq:iterated-Ito}
\end{equation}
Combining \eqref{eq:iterated-Ito} directly with
\eqref{eq:iterated-VI} gives
\begin{eqnarray*}
&&\lambda\alpha_0\mathbb{E}\left\vert
\Pz\widehat X_T^{(k+1)}\right\vert^{2}
+\varepsilon\mathbb{E}\int_{0}^{T}|q_t^{(k+1)}|^{2}\,dt\\
&\leq &\delta\mathbb{E}\int_{0}^{T}
|q_t^{(k+1)}|\,|\widehat G_t^{(k)}|\,dt\\
&\leq &\delta\,\|q^{(k+1)}\|_{\mathbb H^2}
\|\widehat G^{(k)}\|_{\mathbb H^2}.
\end{eqnarray*}
Since the left-hand side is nonnegative,
\begin{equation}
\label{eq:q-before-stability}
\|q^{(k+1)}\|_{\mathbb H^2}
\le
\frac{\delta}{\varepsilon}
\|\widehat G^{(k)}\|_{\mathbb H^2}.
\end{equation}
For $k\ge1$, the increment
$(\widehat X^{(k)},\widehat Y^{(k)},\widehat Z^{(k)})$
satisfies the linear equations
\eqref{eq:linear-forward-increment}--
\eqref{eq:linear-backward-increment} with $q=q^{(k)}$. Hence
\eqref{eq:G-stability} yields
\[
\|\widehat G^{(k)}\|_{\mathbb H^2}
\le
L_G\|q^{(k)}\|_{\mathbb H^2}.
\]
Therefore
\begin{equation}
\label{eq:Picard-contraction}
\|q^{(k+1)}\|_{\mathbb H^2}
\le
\kappa\|q^{(k)}\|_{\mathbb H^2},
\qquad
\kappa=\frac{\delta L_G}{\varepsilon}.
\end{equation}
Choose
\begin{equation}
\label{eq:delta0}
\delta_0
=
\min\left\{1,\frac{\varepsilon}{2L_G}\right\},
\end{equation}
with $\delta_0=1$ if $L_G=0$. For $0<\delta\le\delta_0$, one
has $\kappa\le 1/2$.

Because every iterate $u^{(k)}$ takes values in the uniformly
bounded set $K_t$,
\[
\|q^{(1)}\|_{\mathbb H^2}
\le 2\kappa_K\sqrt T.
\]
It follows that
\begin{equation}
\label{eq:geometric-decay}
\|q^{(k)}\|_{\mathbb H^2}
\le
2\kappa_K\sqrt T\,2^{-(k-1)},
\qquad k\ge1.
\end{equation}
Since $\kappa\le 1/2<1$, the right-hand side tends to zero as
$k\to\infty$. Thus $\{u^{(k)}\}$ is Cauchy in
$\mathbb H^2$. By completeness, there exists
$u\in\mathbb H^2(0,T;\R^m)$ such that
\begin{equation}
\label{eq:u-convergence}
u^{(k)}\longrightarrow u
\qquad\text{in }\mathbb H^2.
\end{equation}

\paragraph{Step 6: Convergence of $X,Y,Z$ and $G$.}
The linear forward SDE estimate gives
\[
\|X^{(k+1)}-X^{(j+1)}\|_{\mathbb S^2}
\le
C_X
\|u^{(k+1)}-u^{(j+1)}\|_{\mathbb H^2}.
\]
Thus $\{X^{(k)}\}$ is Cauchy in $\mathbb S^2$, and there exists
$X\in\mathbb S^2$ such that
$X^{(k)}\longrightarrow X$ in $\mathbb S^2$.
The terminal condition for two backward increments is
\begin{eqnarray*}
Y_T^{(k+1)}-Y_T^{(j+1)}
&=&\lambda\bigl(X_T^{(k+1)}
-\mathbb{E}X_T^{(k+1)}
-\bigl(X_T^{(j+1)}
-\mathbb{E}X_T^{(j+1)}\bigr)\bigr)\\
&=&\lambda\bigl(X_T^{(k+1)}-X_T^{(j+1)}
-\mathbb{E}(X_T^{(k+1)}-X_T^{(j+1)})\bigr).
\end{eqnarray*}
Since
\begin{equation*}
\|X_T^{(k+1)}-X_T^{(j+1)}
-\mathbb{E}(X_T^{(k+1)}-X_T^{(j+1)})\|_{L^2}
\le
\|X_T^{(k+1)}-X_T^{(j+1)}\|_{L^2},
\end{equation*}
the linear BSDE estimate gives
\begin{eqnarray*}
&&\|Y^{(k+1)}-Y^{(j+1)}\|_{\mathbb S^2}
+\|Z^{(k+1)}-Z^{(j+1)}\|_{\mathbb H^2}\\
&\leq &\lambda C_Y\|X_T^{(k+1)}-X_T^{(j+1)}
-\mathbb{E}(X_T^{(k+1)}-X_T^{(j+1)})\|_{L^2}\\
&\leq &\lambda C_Y\|X_T^{(k+1)}-X_T^{(j+1)}\|_{L^2}\\
&\leq &\lambda C_XC_Y\|u^{(k+1)}-u^{(j+1)}\|_{\mathbb H^2}.
\end{eqnarray*}
Hence there exist $Y\in\mathbb S^2$ and
$Z\in\mathbb H^2$ such that
\[
Y^{(k)}\longrightarrow Y\quad\text{in }\mathbb S^2,
\qquad
Z^{(k)}\longrightarrow Z\quad\text{in }\mathbb H^2.
\]
It follows that
\[
\begin{aligned}
\|G^{(k+1)}-G^{(j+1)}\|_{\mathbb H^2}
&\le
\bar B\sqrt T\,\|Y^{(k+1)}-Y^{(j+1)}\|_{\mathbb S^2}\\
&\quad+\bar\sigma
\|Z^{(k+1)}-Z^{(j+1)}\|_{\mathbb H^2},
\end{aligned}
\]
so $\{G^{(k)}\}$ is Cauchy in $\mathbb H^2$. Define
$G:=BY+\sigma Z$.
Then
\begin{equation}
\label{eq:G-convergence}
G^{(k)}\longrightarrow G
\qquad\text{in }\mathbb H^2.
\end{equation}

\paragraph{Step 7: Passage to the limit in the projection.}
Define
\[
\bar u_t
=
\Proj_{K_t}
\left[
-\varepsilon^{-1}
\left(
(\alpha_0+\delta)G_t+\eta_t
\right)
\right].
\]
The projection is nonexpansive:
\[
|\Proj_{K_t}(a)-\Proj_{K_t}(b)|
\le |a-b|.
\]
Applying this pointwise in time and integrating gives
\[
\begin{aligned}
\|u^{(k+1)}-\bar u\|_{\mathbb H^2}
&\le
\frac{1}{\varepsilon}
\left\|
\alpha_0(G^{(k+1)}-G)
+\delta(G^{(k)}-G)
\right\|_{\mathbb H^2}\\
&\le
\frac{\alpha_0}{\varepsilon}
\|G^{(k+1)}-G\|_{\mathbb H^2}
+
\frac{\delta}{\varepsilon}
\|G^{(k)}-G\|_{\mathbb H^2}.
\end{aligned}
\]
By \eqref{eq:G-convergence}, the right-hand side tends to zero.
Hence $u^{(k)}\to\bar u$ in $\mathbb H^2$. Together with
\eqref{eq:u-convergence} and uniqueness of the limit, this gives
$u=\bar u$, namely
\begin{equation}
\label{eq:limit-projection}
u_t
=
\Proj_{K_t}
\left[
-\varepsilon^{-1}
\left(
(\alpha_0+\delta)G_t+\eta_t
\right)
\right].
\end{equation}

\paragraph{Step 8: Passage to the limit in the SDE and BSDE.}
The forward integral equation is
\[
X_t^{(k)}
=
x_0
+
\int_0^t
\left(r_sX_s^{(k)}+B_s^\top u_s^{(k)}\right)ds
+
\int_0^t
(u_s^{(k)})^\top\sigma_s\,dW_s.
\]
The drift terms converge in $L^2$ because
\[
\mathbb{E}\left|
\int_0^T B_t^\top(u_t^{(k)}-u_t)\,dt\right|^2
\le
T\bar B^2
\|u^{(k)}-u\|_{\mathbb H^2}^2
\longrightarrow0,
\]
and the stochastic integrals converge by the It\^o isometry:
\[
\mathbb{E}\int_0^T
\left|(u_t^{(k)}-u_t)^\top\sigma_t\right|^2dt
\le
\bar\sigma^2
\|u^{(k)}-u\|_{\mathbb H^2}^2
\longrightarrow0.
\]
The terms involving $rX^{(k)}$ converge because
$X^{(k)}\to X$ in $\mathbb S^2$. Passing to the limit gives
\[
X_t
=
x_0+\int_0^t(r_sX_s+B_s^\top u_s)\,ds
+\int_0^t u_s^\top\sigma_s\,dW_s,
\]
or equivalently
\[
dX_t=(r_tX_t+B_t^\top u_t)\,dt
+u_t^\top\sigma_t\,dW_t,
\qquad X_0=x_0.
\]
Similarly, the backward integral equation is
\[
Y_t^{(k)}
=
\lambda\Pz X_T^{(k)}-1
+
\int_t^T r_sY_s^{(k)}\,ds
-
\int_t^T (Z_s^{(k)})^\top dW_s.
\]
The convergences
\[
X_T^{(k)}\to X_T\quad\text{in }L^2,\qquad
Y^{(k)}\to Y\quad\text{in }\mathbb S^2,\qquad
Z^{(k)}\to Z\quad\text{in }\mathbb H^2
\]
allow passage to the limit in the terminal condition, the drift
integral, and the stochastic integral, respectively. Therefore
\[
dY_t=-r_tY_t\,dt+Z_t^\top dW_t,
\qquad
Y_T=\lambda (X_T-\mathbb{E}X_T)-1.
\]
Moreover, \eqref{eq:limit-projection} shows that $u_t\in K_t$ and
that the limit satisfies the projection relation with homotopy
parameter $\alpha_0+\delta$.

\paragraph{Step 9: Finite continuation to $\alpha=1$.}
The number $\delta_0$ in \eqref{eq:delta0} depends only on
$\varepsilon,T,\bar r,\bar B,\bar\sigma$, and not on
$\alpha_0$ or $\eta$. Starting from the explicitly solvable
case $\alpha=0$, repeatedly apply the local continuation lemma
along a finite partition
\[
0=\alpha_0<\alpha_1<\cdots<\alpha_N=1,
\qquad
\alpha_{k+1}-\alpha_k\le\delta_0.
\]
Such a finite partition exists because $\delta_0>0$ is fixed.

At $\alpha=1$, choose $\eta=0$. Then
\[
u_t
=
\Proj_{K_t}
\left[-\varepsilon^{-1}(B_tY_t+\sigma_tZ_t)\right]
=
\Phi_{\varepsilon,t}(Y_t,Z_t),
\]
and the resulting triple $(X,Y,Z)$ solves
\eqref{eq:regularized-FBSDE}. Uniqueness was proved in Step 4.
This completes the proof.
\end{proof}

Instead of random box, for the control domain $\Gamma \subset \mathbb{R}^{m}$ is a closed convex set. Typical examples of
such set is $\Gamma =\mathbb{R}_{+}^{m}$ which represents the positive control. We can get the similar result as Theorem \ref{thm:homotopy-FBSDE} (cf. Hu et al. \cite{hhl18}).
\begin{remark}
\label{rem:bsde-vs-hjb-random-coefficients}

The adjoint BSDE formulation remains valid when the market coefficients are
adapted random processes. More precisely, the coefficients $r_t(\omega), B_t(\omega), \sigma_t(\omega)$ may be progressively measurable and uniformly bounded, subject to the usual
integrability and non-degeneracy assumptions. The adjoint equation is then
\[
dY_t
=
-\Bigl(r_tY_t+\text{the corresponding adjoint drift}\Bigr)\,dt
+Z_t^\top dW_t,
\qquad
Y_T=\partial_\mu \Phi(\mathbb P_{X_T})(X_T).
\]
Both \(Y_t\) and \(Z_t\) are adapted processes, so the Hamiltonian gradient is
well defined path by path:
\[
G_t(\omega)=B_t(\omega)Y_t(\omega)
+\sigma_t(\omega)Z_t(\omega).
\]
Consequently, the regularized projection response still has the form
\[
u_t^\varepsilon
=
\operatorname{Proj}_{K_t}
\left[
-\varepsilon^{-1}
\bigl(B_tY_t+\sigma_tZ_t\bigr)
\right].
\]
Thus the BSDE and variational-inequality formulation does not require the
coefficients to be deterministic.

By contrast, the finite-dimensional HJB reduction relies on
deterministic coefficients. If \(r_t\), \(B_t\), or \(\sigma_t\) are general
adapted random processes, the classical deterministic value function
\(V(t,x)\) and the associated deterministic Riccati ordinary differential
equations generally no longer form a closed system: the coefficients of the
HJB equation depend on both \(\omega\) and the information history. A
deterministic PDE reduction is recovered only when the random coefficients
are functions of an observable Markov state, in which case that state must be
added to the HJB state variables. For general adapted coefficients, one must
instead use a stochastic or path-dependent HJB formulation (see Peng \cite{p92}), enlarge the state
space, or, in the linear--quadratic case, solve a backward stochastic Riccati
equation. Therefore, the adjoint BSDE approach is the more robust and direct framework
under random market coefficients, while the finite-dimensional deterministic
HJB reduction is a convenient special case that yields explicit feedback
controls and deterministic Riccati formulas.
\end{remark}

\subsection{Bounded control set}
In this subsection, we first consider the bounded control set. Under (A1), the constraint set is a predictable box
$K_t=[\ell_t,h_t]\subset\mathbb R^{m}$ which is nonempty, closed, convex and
uniformly bounded:
\begin{equation}\label{eq:box-bound}
\kappa_K:=\operatorname*{ess\,sup}_{t,\omega}\sup_{x\in K_t}|x|<\infty .
\end{equation}
For
$\varepsilon>0$ define the Tikhonov-regularized objective
\begin{equation}\label{eq:reg-obj}
J_{\lambda,\varepsilon}(u)
:=J_\lambda(u)+\frac{\varepsilon}{2}\|u\|_{\Htwo}^{2},
\end{equation}
and let
\begin{equation}\label{eq:sets}
S:=\arg\min_{u\in\mathcal U_{\mathrm{ad}}}J_\lambda(u),
\qquad
u^{\varepsilon}:=\arg\min_{u\in\mathcal U_{\mathrm{ad}}}
J_{\lambda,\varepsilon}(u).
\end{equation}

\begin{theorem}
\label{thm:vanishing-reg}
Under Assumption (A1), $S$ is nonempty, closed and convex,
$u^{\varepsilon}$ is uniquely defined for every $\varepsilon>0$, and
\begin{equation}\label{eq:strong-conv}
u^{\varepsilon}\xrightarrow[\varepsilon\downarrow0]{\ \text{strongly in }\Htwo\ }
u^{\dagger},
\qquad
u^{\dagger}:=\Proj_{S}(0),
\end{equation}
where $\Proj_S(0)$ is the metric projection of the origin onto $S$, i.e.\ the
unique admissible optimizer of minimal $\Htwo$-norm:
\begin{equation}\label{eq:minnorm}
\|u^{\dagger}\|_{\Htwo}
=\min_{u\in S}\|u\|_{\Htwo}.
\end{equation}
Moreover,
\begin{enumerate}
\item[(a)] for every $v\in S$ and every $\varepsilon>0$,
\begin{equation}\label{eq:norm-bound}
\|u^{\varepsilon}\|_{\Htwo}\le \|v\|_{\Htwo};
\end{equation}
\item[(b)] $u^{\dagger}$ is characterized by the variational inequality
\begin{equation}\label{eq:proj-VI}
(u^{\dagger},v-u^{\dagger})_{\Htwo}\ge 0
\qquad\text{for all }v\in S;
\end{equation}
\item[(c)] writing $G_t=B_tY_t+\sigma_tZ_t$ for the Hamiltonian gradient
produced by the adjoint BSDE, $u^{\varepsilon}$ solves the regularized
projection variational inequality
\begin{equation}\label{eq:reg-VI}
u_t^{\varepsilon}
=\Proj_{K_t}\!\left[-\varepsilon^{-1}G_t^{\varepsilon}\right]
\quad dt\otimes d\mathbb P\text{-a.e.},
\end{equation}
and $u^{\dagger}$ satisfies the unregularized FBSDVI
\begin{equation}\label{eq:limit-VI}
(G,\,v-u^{\dagger})_{\Htwo}\ge 0
\qquad\text{for all }v\in\mathcal U_{\mathrm{ad}}.
\end{equation}
\end{enumerate}
\end{theorem}

\begin{proof}
We proceed in six steps.
\paragraph{Step~1: Geometry of the admissible set.}
Because the box $K_t$ is almost surely nonempty, closed and convex, the
set $\mathcal U_{\mathrm{ad}}$ is nonempty, closed and convex in the Hilbert
space $\Htwo$. By \eqref{eq:box-bound},
\begin{equation}\label{eq:Uad-bounded}
\|u\|_{\Htwo}^{2}
=\E\!\int_0^T |u_t|^2\,dt
\le m\,\kappa_K^{2}T,
\qquad u\in\mathcal U_{\mathrm{ad}},
\end{equation}
so $\mathcal U_{\mathrm{ad}}$ is bounded and hence weakly sequentially
compact. Closed convex subsets of a Hilbert space are weakly closed.

\paragraph{Step~2: Nonemptiness and convexity of $S$.}
By Lemma \ref{lem:affine-n},
$u\mapsto X_T^{u}$ is bounded affine, and
$x\mapsto \frac{\lambda}{2}\Var(x)-\E[x]$ is convex and continuous on
$L^{2}(\mathcal F_T)$; therefore $J_\lambda$ is continuous and convex on
$\Htwo$, hence weakly sequentially lower semicontinuous. The direct method
of minimization applied to the weakly compact set
$\mathcal U_{\mathrm{ad}}$ yields $S\neq\varnothing$. Moreover, with
$J^{*}:=\min_{\mathcal U_{\mathrm{ad}}}J_\lambda$,
\begin{equation}\label{eq:S-sublevel}
S=\mathcal U_{\mathrm{ad}}\cap\{u:\ J_\lambda(u)\le J^{*}\},
\end{equation}
which is an intersection of closed convex sets; hence $S$ is nonempty,
closed and convex. The Hilbert projection theorem therefore gives a unique
$u^{\dagger}=\Proj_S(0)$, characterized by \eqref{eq:proj-VI}.

\paragraph{Step~3: Existence and uniqueness of $u^{\varepsilon}$.}
The map $u\mapsto \frac{\varepsilon}{2}\|u\|_{\Htwo}^{2}$ is strongly convex
with modulus $\varepsilon$, and $\mathcal U_{\mathrm{ad}}$ is nonempty,
closed and convex. Hence $J_{\lambda,\varepsilon}$ is strongly convex and
continuous, and it admits a unique minimizer $u^{\varepsilon}$. Since
pointwise Euclidean projection onto a closed convex set is characterized by
a variational inequality, the first-order optimality condition of
\eqref{eq:reg-obj} is
\begin{equation}\label{eq:reg-VI-integral}
\bigl(\varepsilon u^{\varepsilon}+G^{\varepsilon},\,
v-u^{\varepsilon}\bigr)_{\Htwo}\ge 0
\quad\text{for all }v\in\mathcal U_{\mathrm{ad}},
\end{equation}
with $G_t^{\varepsilon}=B_tY_t^{\varepsilon}+\sigma_tZ_t^{\varepsilon}$.
Decomposability of the predictable box reduces this to the pointwise
projection formula \eqref{eq:reg-VI}.

\smallskip
\paragraph{Step~4: Uniform comparison of norms.}
Let $v\in S$ be arbitrary. Optimality of $u^{\varepsilon}$ gives
\begin{equation}\label{eq:opt-ineq}
J_\lambda(u^{\varepsilon})
+\frac{\varepsilon}{2}\|u^{\varepsilon}\|_{\Htwo}^{2}
\le J_\lambda(v)+\frac{\varepsilon}{2}\|v\|_{\Htwo}^{2}
=J^{*}+\frac{\varepsilon}{2}\|v\|_{\Htwo}^{2}.
\end{equation}
Because $u^{\varepsilon}\in\mathcal U_{\mathrm{ad}}$, feasibility implies
$J_\lambda(u^{\varepsilon})\ge J^{*}$. Combining this with
\eqref{eq:opt-ineq} and cancelling $\varepsilon/2>0$ yields
\begin{equation}\label{eq:norm-compare}
\|u^{\varepsilon}\|_{\Htwo}^{2}\le \|v\|_{\Htwo}^{2},
\end{equation}
which proves \eqref{eq:norm-bound}. In particular,
$\{u^{\varepsilon}\}_{\varepsilon>0}$ is bounded in $\Htwo$.

\paragraph{Step~5: Weak cluster points lie in $S$, and the weak limit
is $u^{\dagger}$.}
Let $\varepsilon_n\downarrow0$ and choose a weakly convergent subsequence
(relabelled $u_n:=u^{\varepsilon_n}$) with $u_n\rightharpoonup\bar u$ in
$\Htwo$. Weak closedness of $\mathcal U_{\mathrm{ad}}$ gives
$\bar u\in\mathcal U_{\mathrm{ad}}$. Taking $v\in S$ in
\eqref{eq:opt-ineq},
\begin{equation}\label{eq:J-conv}
J^{*}\le J_\lambda(u_n)
\le J^{*}+\frac{\varepsilon_n}{2}
\bigl(\|v\|_{\Htwo}^{2}-\|u_n\|_{\Htwo}^{2}\bigr)
\xrightarrow[n\to\infty]{}J^{*}.
\end{equation}
Weak lower semicontinuity of $J_\lambda$ therefore gives
$J_\lambda(\bar u)\le\liminf_n J_\lambda(u_n)=J^{*}$, and hence
$\bar u\in S$. Next, weak lower semicontinuity of the norm and
\eqref{eq:norm-compare} imply, for every $v\in S$,
\begin{equation}\label{eq:bar-minnorm}
\|\bar u\|_{\Htwo}
\le \liminf_{n}\|u_n\|_{\Htwo}
\le \|v\|_{\Htwo}.
\end{equation}
Thus $\bar u$ is a minimum-norm element of $S$. Such an element is unique
(it is $\Proj_S(0)$), so $\bar u=u^{\dagger}$. Since every sequence
$\varepsilon_n\downarrow0$ has a subsequence with this property, the whole
family converges weakly:
\begin{equation}\label{eq:weak-conv}
u^{\varepsilon}\rightharpoonup u^{\dagger}
\quad\text{in }\Htwo .
\end{equation}

\paragraph{Step~6: Strong convergence.}
Apply \eqref{eq:norm-compare} with $v=u^{\dagger}$:
\begin{equation}\label{eq:norm-limsup}
\limsup_{\varepsilon\downarrow0}\|u^{\varepsilon}\|_{\Htwo}
\le \|u^{\dagger}\|_{\Htwo}.
\end{equation}
Weak convergence gives the reverse inequality
$\|u^{\dagger}\|_{\Htwo}\le\liminf_{\varepsilon\downarrow0}
\|u^{\varepsilon}\|_{\Htwo}$, hence
$\|u^{\varepsilon}\|_{\Htwo}\to\|u^{\dagger}\|_{\Htwo}$. Together with
\eqref{eq:weak-conv}, the Hilbert-space identity
\begin{equation}\label{eq:parallelogram}
\|u^{\varepsilon}-u^{\dagger}\|_{\Htwo}^{2}
=\|u^{\varepsilon}\|_{\Htwo}^{2}+\|u^{\dagger}\|_{\Htwo}^{2}
-2\,(u^{\varepsilon},u^{\dagger})_{\Htwo}
\end{equation}
shows $\|u^{\varepsilon}-u^{\dagger}\|_{\Htwo}\to0$, proving
\eqref{eq:strong-conv}.

Finally, \eqref{eq:affine}--\eqref{eq:Lbounded} and the standard linear
forward--backward estimates imply
$\mathcal L u^{\varepsilon}\to\mathcal L u^{\dagger}$ in $L^{2}$ and
$(Y^{\varepsilon},Z^{\varepsilon})\to(Y,Z)$ in the appropriate BSDE
spaces, hence $G^{\varepsilon}\to G$ in $\Htwo$. Passing to the limit in
\eqref{eq:reg-VI-integral}, using $\varepsilon u^{\varepsilon}\to0$ and
strong convergence of the controls, yields
\begin{equation}\label{eq:limit-VI-proof}
(G,\,v-u^{\dagger})_{\Htwo}\ge 0
\quad\text{for all }v\in\mathcal U_{\mathrm{ad}},
\end{equation}
which is \eqref{eq:limit-VI}. This completes the proof.
\end{proof}

The assumption of uniform boundedness is overly restrictive in practical applications. To accommodate short-selling, we consider an unbounded admissible control set in the following subsection.
\subsection{Unbounded control set}

\label{ass:convex}
Very different from the case of bounded control set, the condition \(S\neq\varnothing\) is imposed in this subsection, because, in infinite dimension
with an unbounded feasible set, the existence of an unregularized optimal
control is not automatic. The direct method of the calculus of variations
proceeds by taking a minimizing sequence and requires that it be bounded and
possess a convergent subsequence: boundedness is immediate when the feasible
set \(\mathcal U_{\mathrm{ad}}\) is bounded, and in finite dimension
boundedness together with closedness already yields a strongly convergent
subsequence by the Bolzano--Weierstrass theorem. In the infinite-dimensional,
unbounded case neither ingredient is available (more details see Table \ref{tab:direct-method}). The closed unit ball is not
compact in the strong topology (the standard orthonormal basis \(e_n\) is
bounded but has no strongly convergent subsequence), and a continuous convex
functional with a finite infimum may admit a minimizing sequence that is
unbounded and escapes to infinity, leaving no weak cluster point and hence no
minimizer.
\begin{table}[ht]
\centering
\caption{Availability of the boundedness and compactness ingredients for the
direct method of the calculus of variations.}
\label{tab:direct-method}
\begin{tabular}{@{}l l p{4.6cm} l@{}}
\toprule
Situation & Boundedness & Compactness & Conclusion \\
\midrule
Feasible set \(\mathcal U_{\mathrm{ad}}\) bounded
  & Bounded
  & ---
  & Applies \\
Finite dimension
  & Boundedness suffices
  & Bounded + closed \(\Rightarrow\) Bolzano--Weierstrass; strongly convergent subsequence
  & Applies \\
Infinite dimension, unbounded
  & Not guaranteed
  & Unit ball not compact
  & No \\
\bottomrule
\end{tabular}
\end{table}

This is precisely what can occur for \(J_\lambda\). The control enters the
terminal wealth only through the bounded affine map $
X_T=a+\mathcal L u,
$
whose range need not be closed and which may have a nontrivial null space (as
in the redundant-asset example, where distinct controls produce the same
terminal wealth). Since \(J_\lambda\) depends on \(u\) only through
\(\mathcal L u\), it is constant along the null-space directions and is
therefore not coercive on \(\mathcal U_{\mathrm{ad}}\); consequently,
finiteness of \(\inf J_\lambda\) does not imply that the infimum is attained.
By contrast, the Tikhonov term is coercive for every \(\varepsilon>0\), since
\[
J_{\lambda,\varepsilon}(u)\ge J^{*}+\frac{\varepsilon}{2}\|u\|^{2},
\]
and this single inequality guarantees existence (the minimizing sequence is
bounded and, in the reflexive control space, has a weakly convergent
subsequence) as well as uniqueness (the quadratic term is strictly convex),
with no boundedness assumption on \(\mathcal U_{\mathrm{ad}}\). The coercive
term disappears exactly at \(\varepsilon=0\), where the existence of \(S\)
must instead be postulated.

The condition is indispensable:
Example~\ref{ex:no-solution} below constructs a continuous convex quadratic
on $\ell^2$ with finite infimum $-\pi^{2}/12$ whose minimizing
sequence satisfies $\|u_N\|\to\infty$ and for which
$S=\varnothing$, so the vanishing-regularization conclusion has no
limit point to select.

\begin{example}
\label{ex:no-solution}
Let $H_{1}=H_{2}=\ell^{2}$, and let $(e_{n})_{n\ge 1}$ be the standard orthogonal basis. Define the diagonal operator
\[
T e_{n}=\frac1n\,e_{n},
\qquad
T(x_{1},x_{2},\ldots)
 =\Bigl(x_{1},\frac{x_{2}}{2},\frac{x_{3}}{3},\ldots\Bigr).
\]
It is bounded and self-adjoint: indeed
\[
\|Tx\|^{2}=\sum_{n=1}^{\infty}\frac{x_{n}^{2}}{n^{2}}
\le \sum_{n=1}^{\infty}x_{n}^{2}=\|x\|^{2},
\qquad T^{*}=T,
\]
and $\|T\|=1$ (attained on $e_{1}$). Its singular values are
$1,1/2,1/3,\ldots$, hence accumulate at $0$, and
\[
\operatorname{Ran}T
=\Bigl\{y\in\ell^{2}:\sum_{n=1}^{\infty}n^{2}y_{n}^{2}<\infty\Bigr\}
\]
is a proper dense subspace, which is not closed.

Set
\[
z:=\sum_{n=1}^{\infty}\frac1n\,e_{n}
=\Bigl(1,\frac12,\frac13,\ldots\Bigr).
\]
Since $
\|z\|^{2}=\sum_{n=1}^{\infty}\frac1{n^{2}}
=\frac{\pi^{2}}{6},
$
we have $z\in\ell^{2}$, and the truncations
\[
z_{N}:=\sum_{n=1}^{N}\frac1n\,e_{n},
\qquad
u_{N}:=\sum_{n=1}^{N}e_{n},
\]
satisfy $z_{N}=Tu_{N}$ and $
\|z-z_{N}\|^{2}
=\sum_{n>N}\frac1{n^{2}}\longrightarrow 0,
$
so $z\in\overline{\operatorname{Ran}T}$. However,
$z\notin\operatorname{Ran}T$: if $z=Tu$ for some $u\in\ell^{2}$,
then $u_{n}=n z_{n}=1$ for every $n$, giving
$u=(1,1,1,\ldots)$, whose squared norm is
$\sum_{n=1}^{\infty}1=\infty$, a contradiction. Equivalently, $
T^{-1}z=\sum_{n=1}^{\infty}e_{n}\notin\ell^{2}.
$
Consider the continuous convex quadratic $
J(u):=\frac12\|Tu\|^{2}-\langle T^{*}z,u\rangle .
$
Completing the square and using $T^{*}=T$ gives
\[
J(u)=\frac12\|Tu-z\|^{2}-\frac12\|z\|^{2}
\ge -\frac12\|z\|^{2}.
\]
The lower bound is approached along $(u_{N})$, since
\[
J(u_{N})
=\frac12\|Tu_{N}-z\|^{2}-\frac12\|z\|^{2}
=\frac12\|z-z_{N}\|^{2}-\frac12\|z\|^{2}
\longrightarrow -\frac12\|z\|^{2}.
\]
Hence
\[
\inf_{u\in H_{1}}J(u)
=-\frac12\|z\|^{2}
=-\frac{\pi^{2}}{12}.
\]
Equality holds if and only if $Tu=z$, i.e. $z\in\operatorname{Ran}T$;
this equation has no solution. Therefore the infimum, although finite,
is never attained, and $
S=\arg\min_{u\in H_{1}}J(u)=\varnothing .
$
Moreover, the minimizing sequence escapes to infinity: $
\|u_{N}\|^{2}=N\longrightarrow\infty .
$
\end{example}
Example \ref{ex:no-solution} is an abstract $\ell^2$ example. The following independent construction shows precisely the same
obstruction in the market model with square-integrable dollar positions.
\begin{example}
\label{ex:unattained}
Take $T=1$, $x_0=r=0$, one risky asset and one Brownian motion, and set
$\sigma_t=t$ and $B_t=\theta t$ with constant $\theta\ne0$. Both coefficients
are bounded and predictable on $[0,1]$. Writing the dollar control as
$u\in\Htwo(0,1;\R)$ and the exposure $z_t=t\,u_t$, the wealth obeys
\begin{equation}\label{eq:unattained-wealth}
\dif X_t=\theta z_t\,\dif t+z_t\,\dif W_t.
\end{equation}
Note that \(\sigma_0=0\), the "responsiveness" of assets to investment at the opening moment is zero. The control $u$ represents how much money you have invested, and z represents the actual risk borne by this sum of money, \(z=\sigma u\). The adducibility requirement is that the square of the dollar position is integrable \(u\in\Htwo(0,1;\R)\). The market price of risk is $B_t/\sigma_t=\theta$, so the state-price density
(again written $\pi$, not $\rho$) is
\begin{equation}\label{eq:unattained-spd}
\pi_t=\exp\!\left(-\theta W_t-\frac12\theta^2t\right),
\qquad \dif\pi_t=-\theta\pi_t\,\dif W_t,
\end{equation}
and the product rule gives the budget identity $\E[\pi_1X_1]=0$ for every
admissible control.

For any terminal claim $\xi\in L^2$ satisfying this budget identity, put
$y=\Pz\xi$. Since $\E\pi_1=1$,
$\E\xi=-\E[\pi_1y]=-\langle\pi_1-1,y\rangle$, and completing the square gives
\begin{equation}\label{eq:unattained-square}
\frac{\lambda}{2}\operatorname{Var}(\xi)-\E\xi
=\frac{\lambda}{2}\left\|y+\frac{\pi_1-1}{\lambda}\right\|^2
-\frac{\norm{\pi_1-1}_{L^2}^2}{2\lambda}.
\end{equation}
Hence the unique budget-feasible $L^2$ optimizer and the finite infimum are
\begin{equation}\label{eq:unattained-opt}
\xi_*=\frac{e^{\theta^2}-\pi_1}{\lambda},
\qquad
\inf J=-\frac{e^{\theta^2}-1}{2\lambda},
\end{equation}
using $\E(\pi_1-1)^2=e^{\theta^2}-1$.

The process that replicates $\xi_*$ is
\begin{equation}\label{eq:unattained-process}
X_t^*=\frac{e^{\theta^2}-e^{\theta^2(1-t)}\pi_t}{\lambda},
\qquad X_0^*=0,\quad X_1^*=\xi_*,
\end{equation}
and differentiation,
$\dif\!\bigl(e^{\theta^2(1-t)}\pi_t\bigr)
=-\theta^2e^{\theta^2(1-t)}\pi_t\,\dif t
-\theta e^{\theta^2(1-t)}\pi_t\,\dif W_t$,
shows that it is driven by
\begin{equation}\label{eq:unattained-exposure}
z_t^*=\frac{\theta}{\lambda}e^{\theta^2(1-t)}\pi_t.
\end{equation}
Because $\E\pi_t^2=e^{\theta^2t}$,
\begin{equation}\label{eq:unattained-second}
\E|z_t^*|^2
=\frac{\theta^2}{\lambda^2}e^{2\theta^2-\theta^2t}
\ge\frac{\theta^2}{\lambda^2}>0.
\end{equation}
The dollar control that would replicate $\xi_*$ is $u_t=z_t^*/t$, and
\begin{equation}\label{eq:unattained-divergence}
\E\int_0^1|u_t|^2\dif t
=\E\int_0^1\frac{|z_t^*|^2}{t^2}\dif t=\infty,
\end{equation}
so $u\notin\Htwo(0,1;\R)$ is inadmissible. Nevertheless, for each $n\ge1$ the
truncated control $u_t^{(n)}=\mathbf 1_{\{t\ge1/n\}}z_t^*/t$
(with $u_0^{(n)}=0$) is admissible, and the corresponding terminal wealths
converge to $\xi_*$ in $L^2$ by Cauchy--Schwarz for the drift and the It\^o
isometry for the martingale term. Continuity of the mean--variance functional
therefore yields the finite infimum in \eqref{eq:unattained-opt} without any
admissible control attaining it.
\end{example}
\begin{remark}
The functional $J$ in Example~\ref{ex:no-solution} is a continuous
convex quadratic with a finite infimum, yet it has no minimizer. The
failure is precisely that its singular values tend to zero, so
$T$ is not bounded below on $(\ker T)^\perp$ and its range is not
closed. Thus $S\neq\varnothing$ cannot be dropped in
Theorem~\ref{thm:vanishing-reg-unbounded}: when $S=\varnothing$ there
is no point to which the regularized minimizers could converge. For
every $\varepsilon>0$ the Tikhonov functional
$J_{\varepsilon}(u)=J(u)+\frac{\varepsilon}{2}\|u\|^{2}$ still has a
unique minimizer (it is coercive and strongly convex), but
$\|u^{\varepsilon}\|\to\infty$ as $\varepsilon\downarrow0$, and the
family has no weak cluster point in $\ell^{2}$.

In Example \ref{ex:unattained}, the obstruction is the vanishing volatility $\sigma_0=0$: the required dollar
position $z_t^*/t$ behaves like $1/t$ near the origin and is not square
integrable. Thus bounded, predictable coefficients alone do not guarantee an
optimizer on the unbounded space $\Htwo$; a general theorem must add an
attainment (or closed-range on the feasible fiber) hypothesis. Closed convex
boxes avoid the issue and always attain their minimum, which is the setting used
by the corrected well-posedness theorem in Section~3.
\end{remark}

The following lemma provides a set of sufficient conditions that guarantee \(S \neq \varnothing\).
\begin{lemma}
\label{lem:restricted-strong}
Let \(H_1,H_2\) be Hilbert spaces, \(T:H_1\to H_2\) bounded linear,
\(\ell\in H_1^*\), and
\[
J(u)=\frac12\|Tu\|_{H_2}^{2}+\ell(u).
\]
Assume:
\begin{enumerate}
\item[\rm (a)] there exists \(\alpha>0\) such that
\[
\|Tu\|_{H_2}\ge\alpha\|u\|_{H_1}
\qquad\text{for all }u\in(\ker T)^\perp,
\]
equivalently, \(T\) has closed range and is bounded below on
\((\ker T)^\perp\);
\item[\rm (b)] \(\ell\) vanishes on \(\ker T\), i.e.
\(\ell\in(\ker T)^\perp\) (this is automatic when
\(\ell\in\operatorname{range}(T^*)\)).
\end{enumerate}
Then \(J\) attains its minimum.
\end{lemma}
\begin{proof}
By the orthogonal decomposition \(H_1=(\ker T)^\perp\oplus\ker T\), write
\(u=u^\perp+u^0\). Since \(T(u^0)=0\) and, by assumption~(b),
\(\ell(u^0)=0\), both \(Tu\) and \(J(u)\) depend only on \(u^\perp\); hence
it suffices to minimize on the closed subspace \(H^\perp:=(\ker T)^\perp\).
There,
\[
J(u^\perp)
=\frac12\|Tu^\perp\|^2+\ell(u^\perp)
\ge\frac{\alpha^2}{2}\|u^\perp\|^2-\|\ell\|\,\|u^\perp\|,
\]
so \(J\) is coercive on \(H^\perp\), and its quadratic part is strongly
convex with modulus \(\alpha^2\). Let \(\{u_n\}\subset H^\perp\) be a
minimizing sequence. Coercivity makes it bounded; reflexivity of the
Hilbert space \(H^\perp\) gives a weakly convergent subsequence
\(u_n\rightharpoonup u^\dagger\in H^\perp\). The continuous convex
functional \(J\) is weakly sequentially lower semicontinuous, so
\[
J(u^\dagger)\le\liminf_n J(u_n)=\inf J.
\]
Thus \(u^\dagger\) is a minimizer; setting its \(\ker T\)-component to zero
gives a canonical optimizer.
\end{proof}

\begin{remark}
\label{rem:ell-kernel-necessary}
Condition~(b) of Lemma~\ref{lem:restricted-strong} cannot be dropped. Take
\(H_1=\mathbb R^{2}\), \(H_2=\mathbb R\), and
\[
T(u_1,u_2)=u_1,\qquad \ell(u_1,u_2)=u_2,
\]
so that \(\ker T=\{(0,u_2):u_2\in\mathbb R\}\) and
\((\ker T)^\perp=\{(u_1,0)\}\) is the \(u_1\)-axis. On this orthogonal
complement,
\[
|T(u_1,0)|=|u_1|=\|(u_1,0)\|\ge\alpha\,\|(u_1,0)\|
\]
with \(\alpha=1\), so condition~(a) holds. Nevertheless,
\[
J(u_1,u_2)=\frac12|T(u_1,u_2)|^{2}+\ell(u_1,u_2)
=\frac12u_1^{2}+u_2
\]
is unbounded below: along the kernel direction \((0,u_2)\),
\[
J(0,u_2)=u_2\longrightarrow-\infty\qquad (u_2\to-\infty).
\]
Hence \(J\) has no minimizer even though condition~(a) is satisfied. The
mechanism is that the quadratic term \(\frac12\|Tu\|^2\) vanishes on
\(\ker T\), leaving a linear functional \(\ell\) that is nonzero there
unrestricted along that direction; moving along the kernel then drives the
objective to \(-\infty\). Condition~(b), \(\ell|_{\ker T}=0\), is therefore
indispensable.
\end{remark}

\begin{corollary}
\label{cor:nonempty}
Assume $\sigma\ne0$ and set $\theta=b/\sigma$. Then
\begin{equation}
\ker(P_0\mathcal L)=\ker\mathcal L
=\{(c,-c):c\in \mathbb H^2(0,T;\R)\},
\label{eq:kernel}
\end{equation}
with orthogonal complement
\begin{equation}
\bigl(\ker(P_0\mathcal L)\bigr)^\perp
=\{(c,c):c\in \mathbb H^2(0,T;\R)\},
\label{eq:complement}
\end{equation}
and there is a constant $\alpha=|\sigma|e^{-3\theta^2T/2}>0$ such that
\begin{equation}
\|(P_0\mathcal L)u\|_{L^2}\ge\alpha\|u\|_{\mathbb H}
\qquad\text{for all }u\in\bigl(\ker(P_0\mathcal L)\bigr)^\perp.
\label{eq:coercivity}
\end{equation}
Equivalently, $P_0\mathcal L$ has closed range and is bounded below on the orthogonal
complement of its kernel. Lemma~\ref{lem:restricted-strong} therefore applies and gives
$S\ne\varnothing$.
\end{corollary}

\begin{proof}
We proceed in five steps.

\smallskip
\paragraph{Step 1: Aggregate dynamics and the state-price density.} Note that $r=0$, the variance of $X_{T}$ does not dependent on $x_{0}$, even if $r_{t}$ is deterministic. So we let $x_{0}=0$.
For the aggregate $v=u_1+u_2$, define the controlled wealth (starting from zero) and the
process $z$ by
\begin{equation}
F_t:=\int_0^t bv_s\,\dif s+\int_0^t\sigma v_s\,\dif W_s,
\qquad z_t:=\sigma v_t,
\end{equation}
so that $\dif F_t=bv_t\,\dif t+z_t\,\dif W_t$, $F_0=0$. It\^o's product rule gives
\begin{equation}
\dif(\rho_tF_t)=\rho_t\,\dif F_t+F_t\,\dif\rho_t+\dif[\rho,F]_t
=\rho_t\bigl(z_t-\theta F_t\bigr)\dif W_t,
\end{equation}
because the drift terms cancel:
$\rho_tbv_t-\theta\rho_tz_t=\rho_tbv_t-\theta\sigma\rho_tv_t=\rho_tbv_t-b\rho_tv_t=0$.
Hence $\rho F$ is a continuous local martingale. Since $v\in\Htwo$, the standard linear
SDE estimates give $F\in\Stwo$ and $z\in\Htwo$, while $\rho\in\Stwo$ because it is a
geometric Brownian motion with bounded moments. Cauchy--Schwarz then yields
\begin{equation}
\E\!\left[\sup_{0\le t\le T}|\rho_tF_t|\right]
\le\left(\E\sup_t|\rho_t|^2\right)^{1/2}
\left(\E\sup_t|F_t|^2\right)^{1/2}<\infty,
\end{equation}
so the local martingale is a true martingale, dominated by its integrable supremum. Using
$\rho_0F_0=0$,
\begin{equation}
\E[\rho_TF_T]=0.
\label{eq:rhoF-zero}
\end{equation}

\smallskip
\paragraph{Step 2: Centering identities.}
Write $P_0F_T=F_T-\E F_T$. Since $\E\rho_T=1$, equation \eqref{eq:rhoF-zero} gives
\begin{equation}
0=\E[\rho_TF_T]
=\E[(\rho_T-1)P_0F_T]+\E F_T,
\end{equation}
and hence
\begin{equation}
\E F_T=-\E[(\rho_T-1)P_0F_T].
\end{equation}
Cauchy--Schwarz, together with
\begin{equation}
\E(\rho_T-1)^2=\E\rho_T^2-2\E\rho_T+1=e^{\theta^2T}-1,
\end{equation}
implies $|\E F_T|^2\le(e^{\theta^2T}-1)\|P_0F_T\|_{L^2}^2$. Therefore
\begin{equation}
\E|F_T|^2=\|P_0F_T\|_{L^2}^2+|\E F_T|^2
\le e^{\theta^2T}\|P_0F_T\|_{L^2}^2,
\label{eq:center-bound}
\end{equation}
or equivalently $\|P_0F_T\|_{L^2}^2\ge e^{-\theta^2T}\E|F_T|^2$.

\smallskip
\paragraph{Step 3: It\^o energy estimate.}
Apply It\^o's formula to $e^{kt}F_t^2$ with $k=2\theta^2$. Since $bv_t=\theta z_t$,
\begin{equation}
\dif\!\left(e^{kt}F_t^2\right)
=e^{kt}\bigl(kF_t^2+2\theta F_tz_t+z_t^2\bigr)\dif t
+2e^{kt}F_tz_t\,\dif W_t.
\end{equation}
The drift integrand satisfies
\begin{equation}
kF_t^2+2\theta F_tz_t+z_t^2
=2\left(\theta F_t+\frac{z_t}{2}\right)^2+\frac{z_t^2}{2}
\ge\frac{z_t^2}{2}.
\end{equation}
Taking expectations (the stochastic integral has mean zero, first by localization and then
by the integrability established in Step~1; the deterministic factor $e^{kt}$ is bounded on
$[0,T]$) and using $F_0=0$ gives
\begin{equation}
\frac12\E\!\int_0^T|z_t|^2\dif t
\le e^{kT}\E|F_T|^2=e^{2\theta^2T}\E|F_T|^2.
\label{eq:energy-bound}
\end{equation}

\smallskip
\paragraph{Step 4: The coercivity estimate.}
Combining \eqref{eq:energy-bound} with \eqref{eq:center-bound},
\begin{equation}
\E\!\int_0^T|z_t|^2\dif t
\le 2e^{3\theta^2T}\|P_0F_T\|_{L^2}^2,
\end{equation}
hence
\begin{equation}
\|P_0F_T\|_{L^2}
\ge\frac{e^{-3\theta^2T/2}}{\sqrt2}
\left(\E\!\int_0^T|z_t|^2\dif t\right)^{1/2}.
\end{equation}
Since $z=\sigma v$ and $P_0F_T=(P_0\mathcal L)u$,
\begin{equation}
\|(P_0\mathcal L)u\|_{L^2}
\ge\frac{|\sigma|e^{-3\theta^2T/2}}{\sqrt2}\,\|v\|_{\Htwo},
\label{eq:bound-aggregate}
\end{equation}
where $\|v\|_{\Htwo}^2=\E\int_0^Tv_t^2\dif t$.

\smallskip
\paragraph{Step 5: Kernels and the bound on the full space.}
The state depends on $u$ only through $v=u_1+u_2$. If $v=0$, then $u_1=-u_2$, so
$u=(c,-c)$ and $\mathcal Lu=0$. Conversely, if $(P_0\mathcal L)u=0$, then $F_T$ is
deterministic; its centered stochastic part $\int_0^T\sigma v_t\dif W_t$ is a
mean-zero martingale with quadratic variation $\sigma^2\int_0^Tv_t^2\dif t$, which must
vanish. Because $\sigma\ne0$, this forces $v=0$. Thus
\begin{equation}
\ker(P_0\mathcal L)=\ker\mathcal L
=\{(c,-c):c\in \mathbb H^2(0,T;\R)\}.
\end{equation}
The orthogonal complement for the $\Htwo$ inner product is
$\{(c,c):c\in \mathbb H^2(0,T;\R)\}$, since $\int_0^T(cd-cd)\dif t=0$. For $u=(c,c)$ in this
complement, $v=2c$ and
\begin{equation}
\|v\|_{\Htwo}^2=4\E\!\int_0^Tc_t^2\dif t,\qquad
\|u\|_{\Htwo}^2=2\E\!\int_0^Tc_t^2\dif t,
\end{equation}
so $\|v\|_{\Htwo}=\sqrt2\,\|u\|_{\Htwo}$. Substituting into
\eqref{eq:bound-aggregate} yields
\begin{equation}
\|(P_0\mathcal L)u\|_{L^2}
\ge|\sigma|e^{-3\theta^2T/2}\|u\|_{\Htwo}
=\alpha\|u\|_{\Htwo},
\end{equation}
with $\alpha=|\sigma|e^{-3\theta^2T/2}>0$. This is precisely the hypothesis of
Lemma~\ref{lem:restricted-strong} (with $T$ there equal to $P_0\mathcal L$), so $J_\lambda$
attains its minimum and $S\ne\varnothing$. The same kernel is the one required in the
later minimum-energy corollary. This completes the proof.
\end{proof}


We impose the following condition:
\begin{enumerate}
\item[\textbf{(A2)}] \label{ass:market2}Assume that $r$, $B$, and $\sigma $
are bounded $\left\{ \mathcal{F}_{t}\right\} _{t\geq 0}$-adapted processes. Let $\mathcal U_{\mathrm{ad}}\subseteq\Htwo$ be nonempty, closed and convex
(so $\mathcal U_{\mathrm{ad}}=\Htwo$ is allowed in the unconstrained case). The
initial wealth $x_{0}$ is deterministic.
\end{enumerate}

\begin{theorem}
\label{thm:vanishing-reg-unbounded}
Under Assumption (A2), and assuming only the nonemptiness
$S\neq\varnothing$ stipulated there, the following hold:
\begin{enumerate}
\item[\rm (i)] the solution set $S$ is closed and convex;
\item[\rm (ii)] the regularized minimizer $u^{\varepsilon}$ is uniquely
defined for every $\varepsilon>0$;
\item[\rm (iii)]
\[
u^{\varepsilon}\xrightarrow[\varepsilon\downarrow0]{\ \text{strongly in }\Htwo\ }
u^{\dagger},
\qquad
u^{\dagger}:=\Proj_S(0),
\]
where
\[
\|u^{\dagger}\|_{\Htwo}
=\min_{u\in S}\|u\|_{\Htwo},
\qquad
(u^{\dagger},v-u^{\dagger})_{\Htwo}\ge 0
\quad\text{for all }v\in S.
\]
\end{enumerate}
\end{theorem}
\begin{remark}
\label{rem:affine-fiber}
In the present quadratic case, once the linear functional $\ell$ vanishes on the kernel
(cf.\ Lemma \ref{lem:restricted-strong}(b)), every element of the unconstrained optimal set has the form
$u^*+q$ with $q\in\ker(P_0\mathcal L)$. The constrained solution set is therefore the affine
fiber intersected with the feasible set:
\begin{equation}
S=\mathcal U_{ad}\cap\big(u^*+\ker(P_0\mathcal L)\big).
\end{equation}
The minimum-norm element $u^\dagger=\Pi_S(0)$ is precisely the point of this
intersection whose kernel component vanishes.
\end{remark}

\begin{proof}
We proceed in five steps.

\paragraph{Step~1: Geometry of $S$ and existence of the projection.}
Because $J_\lambda$ is continuous and convex, its sublevel set
$\{J_\lambda\le J^{*}\}$ is closed and convex. Hence
\[
S=\mathcal U_{\mathrm{ad}}\cap\{J_\lambda\le J^{*}\}
\]
is an intersection of nonempty, closed and convex sets.
Therefore $S$ is nonempty, closed and convex. The Hilbert projection
theorem yields a unique $u^{\dagger}=\Proj_S(0)$, characterized by the
variational inequality
\[
(u^{\dagger},v-u^{\dagger})_{\Htwo}\ge0
\qquad\text{for every }v\in S.
\]

\paragraph{Step~2: Existence and uniqueness of $u^{\varepsilon}$.}
We first establish coercivity, then apply the direct method, and finally
use strong convexity for uniqueness.

\smallskip
\emph{Coercivity.}
Since $S\neq\varnothing$ and $J^{*}$ is the attained minimum of
$J_\lambda$ over $\mathcal U_{\mathrm{ad}}$, $J^{*}$ is finite and
$J_\lambda(u)\ge J^{*}$ for every $u\in\mathcal U_{\mathrm{ad}}$. Hence
\[
J_{\lambda,\varepsilon}(u)
=J_\lambda(u)+\frac{\varepsilon}{2}\|u\|_{\Htwo}^{2}
\ge J^{*}+\frac{\varepsilon}{2}\|u\|_{\Htwo}^{2}.
\]
As $\|u\|_{\Htwo}\to\infty$ the right-hand side diverges to $+\infty$,
so $J_{\lambda,\varepsilon}$ is coercive on $\mathcal U_{\mathrm{ad}}$; it
is continuous and convex because $J_\lambda$ is continuous and convex and
$\frac{\varepsilon}{2}\|\cdot\|_{\Htwo}^{2}$ is continuous and convex.

\smallskip
\emph{The direct method.}
Let $\{u_n\}\subset\mathcal U_{\mathrm{ad}}$ be a minimizing sequence:
\[
J_{\lambda,\varepsilon}(u_n)\longrightarrow
\inf_{u\in\mathcal U_{\mathrm{ad}}}J_{\lambda,\varepsilon}(u)=:m_\varepsilon.
\]
Coercivity implies that $\{u_n\}$ is bounded. Indeed, let $\{u^{(n)}\}_{n\ge1}\subset\mathcal U_{\mathrm{ad}}$ be a minimizing
sequence, meaning
\[
J_{\lambda,\varepsilon}(u^{(n)})\longrightarrow
m:=\inf_{u\in\mathcal U_{\mathrm{ad}}}J_{\lambda,\varepsilon}(u),
\]
where $m$ is finite: the lower bound above gives $m\ge J^{*}>-\infty$, while
evaluating at any feasible point gives $m<\infty$.

We claim that $\{u^{(n)}\}$ is bounded. Suppose, for contradiction, that it
is unbounded. Then there is a subsequence $\{u^{(n_k)}\}$ with
$\|u^{(n_k)}\|_{\Htwo}\to\infty$: for every integer $k$ one may choose an
index $n_k$ (with $n_1<n_2<\cdots$) such that $\|u^{(n_k)}\|>k$. Coercivity
would then force
\[
J_{\lambda,\varepsilon}(u^{(n_k)})\longrightarrow+\infty,
\]
which contradicts $J_{\lambda,\varepsilon}(u^{(n_k)})\to m<\infty$, since a
convergent sequence cannot have a subsequence tending to $+\infty$. Hence
every minimizing sequence is bounded.

The Hilbert space $\Htwo$ is reflexive, so by the
Eberlein--\v Smulian theorem (cf. \cite{Brezis 2010}) a bounded sequence has a weakly convergent
subsequence; relabelling, write
\[
u_n\rightharpoonup\bar u\qquad\text{in }\Htwo.
\]
The set $\mathcal U_{\mathrm{ad}}$ is closed and convex in a Hilbert space,
hence it is weakly closed (a consequence of the Hahn--Banach theorem and
Mazur's lemma, see \cite{Brezis 2010}), so $\bar u\in\mathcal U_{\mathrm{ad}}$.

A continuous convex function on a Hilbert space is weakly sequentially
lower semicontinuous. Both $J_\lambda$ and
$\frac{\varepsilon}{2}\|\cdot\|_{\Htwo}^{2}$ have this property, and so
does their sum:
\[
J_{\lambda,\varepsilon}(\bar u)
\le\liminf_{n\to\infty}J_{\lambda,\varepsilon}(u_n)
=m_\varepsilon.
\]
Since $\bar u\in\mathcal U_{\mathrm{ad}}$, the reverse inequality
$J_{\lambda,\varepsilon}(\bar u)\ge m_\varepsilon$ holds; hence
equality and $\bar u$ is a minimizer.

\emph{Uniqueness via strong convexity.}
The squared norm satisfies, for $u,v\in\Htwo$ and $\theta\in[0,1]$,
\[
\|\theta u+(1-\theta)v\|_{\Htwo}^{2}
\le \theta\|u\|_{\Htwo}^{2}
+(1-\theta)\|v\|_{\Htwo}^{2}
-\theta(1-\theta)\|u-v\|_{\Htwo}^{2},
\]
which is the strong convexity inequality with modulus $1$ for
$\|\cdot\|_{\Htwo}^{2}$. Hence
$\frac{\varepsilon}{2}\|\cdot\|_{\Htwo}^{2}$ is strongly convex with
modulus $\varepsilon$. Adding the merely convex functional $J_\lambda$
does not decrease this modulus, so $J_{\lambda,\varepsilon}$ is strongly
convex with modulus $\varepsilon$: for $u,v\in\mathcal U_{\mathrm{ad}}$ and
$\theta\in[0,1]$,
\[
J_{\lambda,\varepsilon}\bigl(\theta u+(1-\theta)v\bigr)
\le \theta J_{\lambda,\varepsilon}(u)
+(1-\theta)J_{\lambda,\varepsilon}(v)
-\frac{\varepsilon}{2}\theta(1-\theta)\|u-v\|_{\Htwo}^{2}.
\]
Suppose $u'$ and $u''$ were two distinct minimizers. Taking
$\theta=\tfrac12$ in this inequality gives
\[
J_{\lambda,\varepsilon}\!\left(\frac{u'+u''}{2}\right)
\le J_{\lambda,\varepsilon}(u')
-\frac{\varepsilon}{8}\|u'-u''\|_{\Htwo}^{2}
<J_{\lambda,\varepsilon}(u'),
\]
because $\frac{u'+u''}{2}\in\mathcal U_{\mathrm{ad}}$ by convexity. This
contradicts optimality. Therefore $u'=u''$, and the minimizer
$u^{\varepsilon}$ is unique.

\paragraph{Step~3: Uniform comparison bound.}
Fix any $v\in S$, i.e. $v$ is any optimizer of the unregularized problem
(so $J_\lambda(v)=J^{*}$). Optimality of $u^{\varepsilon}$ gives
\[
J_\lambda(u^{\varepsilon})
+\frac{\varepsilon}{2}\|u^{\varepsilon}\|_{\Htwo}^{2}
\le J_\lambda(v)
+\frac{\varepsilon}{2}\|v\|_{\Htwo}^{2}
=J^{*}+\frac{\varepsilon}{2}\|v\|_{\Htwo}^{2}.
\]
Since $u^{\varepsilon}\in\mathcal U_{\mathrm{ad}}$,
$J_\lambda(u^{\varepsilon})\ge J^{*}$. Combining these inequalities and
cancelling $\varepsilon/2>0$,
\[
\|u^{\varepsilon}\|_{\Htwo}^{2}\le \|v\|_{\Htwo}^{2}
\quad\text{for every }v\in S.
\]
In particular, the family $\{u^{\varepsilon}\}_{\varepsilon>0}$ is bounded
in $\Htwo$, without any boundedness assumption on
$\mathcal U_{\mathrm{ad}}$.

\paragraph{Step~4: Weak cluster points lie in $S$ and the weak limit
is $u^{\dagger}$.}
Let $\varepsilon_n\downarrow0$. By Step~3 the family is bounded, so there
is a subsequence (relabelled $u_n:=u^{\varepsilon_n}$) converging weakly,
$u_n\rightharpoonup\bar u$ in $\Htwo$. The closed convex set
$\mathcal U_{\mathrm{ad}}$ is weakly closed, hence
$\bar u\in\mathcal U_{\mathrm{ad}}$. Taking any $v\in S$ in the comparison
inequality of Step~3,
\[
J^{*}\le J_\lambda(u_n)
\le J^{*}+\frac{\varepsilon_n}{2}
\bigl(\|v\|_{\Htwo}^{2}-\|u_n\|_{\Htwo}^{2}\bigr)
\xrightarrow[n\to\infty]{}J^{*},
\]
where the bracket is nonnegative and bounded. Hence
$J_\lambda(u_n)\to J^{*}$. Weak lower semicontinuity of $J_\lambda$ gives
\[
J_\lambda(\bar u)\le\liminf_{n}J_\lambda(u_n)=J^{*},
\]
and feasibility gives $J_\lambda(\bar u)\ge J^{*}$, so
$\bar u\in S$. Furthermore, weak lower semicontinuity of the norm and
Step~3 yield, for every $v\in S$,
\[
\|\bar u\|_{\Htwo}
\le\liminf_{n}\|u_n\|_{\Htwo}
\le \|v\|_{\Htwo}.
\]
Thus $\bar u$ is the unique minimum-norm element of $S$, namely
$\bar u=\Proj_S(0)=u^{\dagger}$. Every subsequence of
$\{u^{\varepsilon}\}$ contains a further subsequence converging weakly to
$u^{\dagger}$; since the family is bounded, this characterizes weak
convergence of the whole family:
\[
u^{\varepsilon}\rightharpoonup u^{\dagger}
\quad\text{in }\Htwo .
\]

\paragraph{Step~5: Strong convergence.}
Apply the comparison inequality of Step~3 with $v=u^{\dagger}$:
\[
\|u^{\varepsilon}\|_{\Htwo}\le \|u^{\dagger}\|_{\Htwo},
\qquad
\limsup_{\varepsilon\downarrow0}\|u^{\varepsilon}\|_{\Htwo}
\le \|u^{\dagger}\|_{\Htwo}.
\]
Weak convergence gives the reverse inequality
$\|u^{\dagger}\|_{\Htwo}\le
\liminf_{\varepsilon\downarrow0}\|u^{\varepsilon}\|_{\Htwo}$, so
$\|u^{\varepsilon}\|_{\Htwo}\to\|u^{\dagger}\|_{\Htwo}$. The
parallelogram identity,
\[
\|u^{\varepsilon}-u^{\dagger}\|_{\Htwo}^{2}
=\|u^{\varepsilon}\|_{\Htwo}^{2}
+\|u^{\dagger}\|_{\Htwo}^{2}
-2\,(u^{\varepsilon},u^{\dagger})_{\Htwo},
\]
together with $(u^{\varepsilon},u^{\dagger})_{\Htwo}
\to\|u^{\dagger}\|_{\Htwo}^{2}$, shows
\[
\|u^{\varepsilon}-u^{\dagger}\|_{\Htwo}^{2}
\longrightarrow 0.
\]
This proves the claimed strong convergence and completes the proof.
\end{proof}
\begin{remark}
The BSDE argument and the
Hilbert-space argument establish the same object but answer different
questions. After adding the Tikhonov term, the pointwise Hamiltonian is
strictly convex in the control, so its minimizer is unique and admits an
explicit feedback expression in terms of the adjoint pair
$(Y^\varepsilon,Z^\varepsilon)$; this \emph{characterizes}
$u^\varepsilon$---it exhibits the optimal control as a functional of the
adjoint process. By contrast, the Hilbert-space route proves uniqueness
without producing a formula: it uses only the strict convexity and coercivity
of the regularized cost on $\mathbb H^2(0,T;\mathbb R^m)$ together with the
projection theorem. Uniqueness there is a statement about a strictly convex
functional on a Hilbert space, not about a particular optimality system.
The two viewpoints are therefore complementary: the BSDE supplies the
explicit feedback (the ``who''), while the Hilbert-space argument supplies
the abstract existence and uniqueness guarantee (the ``why''). This
complementarity is decisive in the vanishing-regularization limit: the
Hilbert-space projection remains well defined on the closed convex solution
set and selects the minimum-norm portfolio, whereas the BSDE yields the
regularized feedback whose strong limit realizes precisely that projection.
\end{remark}

\begin{corollary}
\label{cor:redundant}
Suppose $m=2$, both assets have identical excess return and volatility,
\begin{equation}\label{eq:redundant-wealth}
dX_t=b(u_{1,t}+u_{2,t})\,dt+\sigma(u_{1,t}+u_{2,t})\,dW_t,
\end{equation}
and the unregularized optimizer determines the aggregate
$v_t^{*}=u_{1,t}^{*}+u_{2,t}^{*}$. Then
$\ker\mathcal L=\operatorname{span}\{(1,-1)\}$ and, in the unconstrained
case (and, by symmetric pointwise projection, also for a symmetric box),
\begin{equation}\label{eq:equal-split}
u_{1,t}^{\dagger}=u_{2,t}^{\dagger}=\frac12 v_t^{*},
\qquad
u^{\varepsilon}\xrightarrow[\varepsilon\downarrow0]{\Htwo}u^{\dagger}.
\end{equation}
\end{corollary}

\begin{proof}
The solution set is the affine optimal face
$S=\{u\in\mathcal U_{\mathrm{ad}}:\ u_{1,t}+u_{2,t}=v_t^{*}\ \text{a.e.}\}$.
Minimizing $\|u\|_{\Htwo}^{2}$ subject to $u_1+u_2=v^{*}$ is a quadratic
program; its Lagrange multiplier equations read
$2u_{1,t}^{\dagger}=\mu_t$ and $2u_{2,t}^{\dagger}=\mu_t$, hence
$u_{1,t}^{\dagger}=u_{2,t}^{\dagger}=v_t^{*}/2$. This is precisely
$\Proj_S(0)$, and \eqref{eq:equal-split} follows from
Theorem~\ref{thm:vanishing-reg}. For a symmetric box the same conclusion
holds because the two coordinates share the same gradient and the same
constraint interval, so the pointwise projection of
$-\varepsilon^{-1}G_t^{\varepsilon}$ (with identical components) is
coordinatewise identical.
\end{proof}

\begin{remark}
The theorem is the Hilbert-space vanishing-regularization principle in the
FBSDVI setting: the BSDE/Hamiltonian gradient determines the optimal face
(the aggregate exposure), while the quadratic penalty selects its unique
element of minimal energy. Strong (not merely weak) convergence is the key
conclusion and follows from convergence of the norms together with weak
convergence via the parallelogram identity.
\end{remark}

\begin{remark}
\label{rem:affine-min-norm}
Let $N:=\ker\mathcal L$, which is a closed linear subspace of the control
space. The set of optimal controls is the affine subspace
\begin{equation}
S=u^*+N
\end{equation}
for any particular optimal control $u^*$. In the two redundant-asset case
where the optimality condition reads
\begin{equation}
S=\{u\in\mathcal U_{\mathrm{ad}}:\ u_{1,t}+u_{2,t}=v_t^*\ \text{a.e.}\},
\end{equation}
the homogeneous part and a convenient particular solution are
\begin{equation}
N=\{c_t(1,-1):c\in\Htwo\},\qquad
u_t^*=\Bigl(\tfrac{v_t^*}{2},\tfrac{v_t^*}{2}\Bigr),
\end{equation}
so every element of $S$ admits the explicit parametrization
\begin{equation}
u_t=\Bigl(\tfrac{v_t^*}{2}+c_t,\ \tfrac{v_t^*}{2}-c_t\Bigr),
\qquad c\in\Htwo,
\end{equation}
and the two components always sum to $v_t^*$.

For any $u\in S$, write the orthogonal decomposition
$u=u^\perp+u^0$ with $u^0\in N$ and $u^\perp\in N^\perp$. The Pythagorean
identity gives
\begin{equation}
\|u\|^2=\|u^\perp\|^2+\|u^0\|^2,
\end{equation}
where the orthogonal component $u^\perp$ is the same for every $u\in S$,
while $\|u^0\|^2\ge0$. Consequently,
\begin{equation}
u^0=0
\quad\Longleftrightarrow\quad
\|u\|=\min_{v\in S}\|v\|.
\end{equation}
Indeed, if $u^0=0$ then $\|u\|^2=\|u^\perp\|^2\le\|v\|^2$ for every
$v\in S$; conversely, a minimum-norm point must satisfy
$\|u^0\|^2=0$, hence $u^0=0$. This unique point is precisely
\begin{equation}
u^\dagger=\Proj_S(0)\in N^\perp,
\qquad
(u^\dagger,v-u^\dagger)\ge0\quad\text{for all }v\in S.
\end{equation}
In the parametrization above it corresponds to $c=0$, namely the
equal-weight control $u^\dagger=(v^*/2,v^*/2)$. Thus vanishing of the
null-space component and selection of the minimum-norm optimizer are
equivalent, and the quadratic penalty realizes exactly this selection.
\end{remark}

\subsection{Multiple accounts}
\label{ma}
In this subsection, we examine the case of the multiple accounts, i.e. $n>1$. As demonstrated below, the proposed method is independent of the dimensionality of the wealth vector.

Let the $n$-dimensional wealth vector $X$ satisfy
\begin{equation}
\label{eq:wealth-n}
dX_t=\bigl(r_t I_n X_t+U_t\beta_t\bigr)\,dt+U_t\Sigma_t\,dW_t,
\qquad X_0=x_0,
\end{equation}
where
$
U_t\in\mathbb R^{n\times m},
\beta_t\in\mathbb R^{m},
\Sigma_t\in\mathbb R^{m\times d},
$
and $r_t$ is scalar. They are all adapted and bounded. The mean-variance objective for the $n$ accounts is the sum of the
$n$ scalar criteria,
\begin{equation}
\label{eq:cost-n}
J(u)
= \frac{\lambda}{2}\,\mathbb E\!\left|X_T-\mathbb E[X_T]\right|^2
  -\mathbf 1_n^\top \mathbb E[X_T]
= \sum_{i=1}^{n}\left[
  \frac{\lambda}{2}\operatorname{Var}(X_T^i)
  -\mathbb E[X_T^i]\right],
\end{equation}
where $\mathbf 1_n^\top=[1,...,1]$.
\paragraph{The vector mean-field adjoint BSDE.}

The adjoint processes are lifted to the vector and matrix dimensions
$
Y\in \mathbb S^2(0,T;\mathbb R^n),
Z\in \mathbb H^2(0,T;\mathbb R^{n\times d}),
$
and solve the linear mean-field adjoint BSDE
\begin{equation}
\label{eq:mf-bsde-n}
\left\{
\begin{aligned}
dY_t &= -r_t\,Y_t\,dt+Z_t\,dW_t,\\
Y_T &= \lambda\bigl(X_T-\mathbb E[X_T]\bigr)-\mathbf 1_n,
\end{aligned}
\right.
\end{equation}
where $Z_t\,dW_t\in\mathbb R^n$. The scalar ``$-1$'' becomes the vector
``$-\mathbf 1_n$'', while the mean-field dependence enters only through the terminal
condition; the driver itself contains no law-dependent term because the wealth SDE
\eqref{eq:wealth-n} does not depend on the law of $X$.

\paragraph{Hamiltonian and the switching matrix.}

The Hamiltonian in matrix form is
\begin{equation}
\label{eq:hamiltonian-n}
H_t(X,U,Y,Z)
= Y^\top\!\bigl(r_t I_n X+U\beta_t\bigr)
  +\operatorname{tr}\!\bigl(Z^\top U\Sigma_t\bigr)
= \big\langle Y,\,r_t I_n X+U\beta_t\big\rangle
  +\big\langle Z,\,U\Sigma_t\big\rangle_F ,
\end{equation}
where $\langle\cdot,\cdot\rangle_F$ denotes the Frobenius inner product. Differentiating
with respect to $U\in\mathbb R^{n\times m}$ gives the switching matrix
\begin{equation}
\label{eq:switching-n}
G_t:=\nabla_U H_t
= Y_t\beta_t^\top+Z_t\Sigma_t^\top\in\mathbb R^{n\times m}.
\end{equation}
For $n=1$, equation \eqref{eq:switching-n} reduces to the scalar switching vector
$G_t=B_tY_t+\sigma_tZ_t$.

\paragraph{The FBSDVI and the regularized response.}

Without regularization, optimality is the pointwise variational inequality
\begin{equation}
\label{eq:vi-n}
U_t\in K_t,\qquad
\big\langle G_t,V-U_t\big\rangle_F\ge 0,
\qquad \forall V\in K_t.
\end{equation}
With the Tikhonov $L^2$ penalty, the unique response is the Frobenius projection
\begin{equation}
\label{eq:response-n}
U_t^\varepsilon
= \Pi_{K_t}\!\left[-\varepsilon^{-1}G_t^\varepsilon\right],
\end{equation}
equivalently,
\begin{equation}
\label{eq:vi-eps-n}
\big\langle \varepsilon U_t^\varepsilon+G_t^\varepsilon,
V-U_t^\varepsilon\big\rangle_F\ge 0,
\qquad \forall V\in K_t.
\end{equation}

\paragraph{Vectorized form.}

Setting $u_t:=\operatorname{vec}(U_t)\in\mathbb R^{nm}$, the vectorized switching
quantity is
\begin{equation}
\label{eq:vec-n}
g_t:=\operatorname{vec}(G_t)
= (\beta_t\otimes I_n)\,Y_t
  +(\Sigma_t\otimes I_n)\,\operatorname{vec}(Z_t).
\end{equation}
The Frobenius projection \eqref{eq:response-n} then becomes the standard Euclidean
projection in $\mathbb R^{nm}$, which for a box reduces to componentwise clipping.

Because the state drift is $r_t I_n$ and the terminal cost \eqref{eq:cost-n} is separable
across the accounts, the components of $Y$ satisfy $n$ decoupled copies of the scalar
adjoint BSDE. Hence the Hilbert-space operator, projection, and strong-convergence
arguments apply unchanged for every finite $n$, which is precisely the dimension-free
mechanism that covers the $n$-account problem.

\section{Redundant-risky-asset model}
It suffices to consider the case $r=0$. Indeed, for $r>0$, introduce the
discounted wealth and position
$
\widetilde X_t:=e^{-rt}X_t, \widetilde v_t:=e^{-rt}v_t.
$
Since $e^{-rt}$ is of finite variation, It\^o's product rule yields
\[
d\widetilde X_t=b\,\widetilde v_t\,dt+\sigma\widetilde v_t\,dW_t,
\]
and the objective becomes
\[
J_\varepsilon
=\frac{\lambda e^{2rT}}{2}\Var(\widetilde X_T)
-e^{rT}\mathbb E[\widetilde X_T]
+\frac{\varepsilon}{4}\mathbb E\!\int_0^T e^{2rt}\widetilde v_t^2\,dt.
\]
Thus the problem is identical in form to the zero-rate one, up to the
deterministic time weight $e^{2rt}$ on the running cost and the rescaled
terminal data ($A_T=\lambda e^{2rT}$ and $B_T=-e^{rT}$). In particular, the
geometry of the feasible set, the pointwise projection, and the
sign-preserving stochastic-exponential structure of the tracking deviation
are unchanged; only deterministic weights and boundary values change.
Hence, throughout the sequel, we take $r=0$ without loss of generality.

Let $W$ be a one-dimensional Brownian motion. There are two risky assets with the same excess return $b>0$ and the same volatility $\sigma>0$. Writing
$
u_t=(u_{1,t},u_{2,t})^\top,
$
and taking the interest rate equal to zero, the wealth process satisfies
\begin{equation}
\label{eq:wealth-r}
\dif X_t
=
b(u_{1,t}+u_{2,t})\,\dif t
+
\sigma(u_{1,t}+u_{2,t})\,\dif W_t,
\qquad X_0=x_0.
\end{equation}
Here $b$ is the excess return over the risk-free asset. Since $r=0$, the regime $b>0$ is
the economically meaningful one: if $b\le0$, the risky asset carries no positive risk
premium and the optimal long-only position is $v_t^*=0$ (hold cash).

The mean-variance objective is
\begin{equation}
\label{eq:objective}
J(u)=\frac{\lambda}{2}\Var(X_T)-\E[X_T],
\qquad \lambda>0.
\end{equation}
Introduce the aggregate position and the long--short difference
\begin{equation}
\label{eq:vw}
v_t:=u_{1,t}+u_{2,t},
\qquad
w_t:=u_{1,t}-u_{2,t}.
\end{equation}
Then
\begin{equation}
\label{eq:inverse-vw}
u_{1,t}=\frac{v_t+w_t}{2},
\qquad
u_{2,t}=\frac{v_t-w_t}{2},
\end{equation}
and the wealth equation becomes
\begin{equation}
\label{eq:wealth-v}
\dif X_t=bv_t\,\dif t+\sigma v_t\,\dif W_t.
\end{equation}
Thus the state is completely independent of $w_t$. Moreover,
\begin{equation}
\label{eq:energy-vw}
u_{1,t}^2+u_{2,t}^2
=
\frac{v_t^2+w_t^2}{2}.
\end{equation}
Hence the control energy is
\begin{equation}
\label{eq:H2-energy-vw}
\|u\|_{\Htwo(0,T;\R^2)}^2
=
\E\int_0^T(u_{1,t}^2+u_{2,t}^2)\,\dif t
=
\frac12\E\int_0^T(v_t^2+w_t^2)\,\dif t.
\end{equation}
The explicit calculations below are written for the unconstrained problem. For a symmetric box constraint $K=[\ell,h]^2$, the same argument applies after replacing the first-order condition by the Euclidean projection; because both assets have the same adjoint gradient and the same coordinate constraint, the regularized controls remain equal whenever the constraints are treated symmetrically.

\subsection{The Lions derivative and the adjoint BSDE}

The terminal cost is the law-dependent functional
\[
f(\mu)=\frac{\lambda}{2}\operatorname{Var}_\mu(X)-\E_\mu[X].
\]
Its lift to $L^2(\mathcal F_T;\R)$ is
\[
\widetilde f(\xi)
=
\frac{\lambda}{2}
\left(\E[\xi^2]-(\E\xi)^2\right)-\E\xi.
\]
For a random variable $\xi$ with law $\mu$, the Lions derivative at the point $x$ is
\begin{equation}
\label{eq:Lions-terminal}
\partial_\mu f(\mu)(x)
=
\lambda\left(x-\int y\,\mu(\dif y)\right)-1.
\end{equation}
Therefore the terminal condition of the adjoint process is
\begin{equation}
\label{eq:Y-terminal}
Y_T
=
\lambda(X_T-\E X_T)-1.
\end{equation}
Since the interest rate is zero and the coefficients of \eqref{eq:wealth-v} do not depend on $X$, the adjoint BSDE is
\begin{equation}
\label{eq:adjoint-BSDE}
\dif Y_t=Z_t\,\dif W_t,
\qquad
Y_T=\lambda(X_T-\E X_T)-1.
\end{equation}
For a nonzero interest rate $r$, the backward equation would instead contain the usual discounting term $-rY_t\,\dif t$.

The Hamiltonian is
\begin{equation}
\label{eq:Hamiltonian}
H(u,y,z)
=
by(u_1+u_2)+z\sigma(u_1+u_2).
\end{equation}
Writing
$
G_t:=bY_t+\sigma Z_t,
$
we have
$
H=G_t(u_{1,t}+u_{2,t}).
$
The two control gradients are identical:
\begin{equation}
\label{eq:identical-gradients}
\frac{\partial H}{\partial u_1}
=
\frac{\partial H}{\partial u_2}
=
G_t.
\end{equation}
This identity has a direct interpretation. The adjoint BSDE determines the aggregate position $v_t$, but it contains no information about the difference $w_t$. In particular, a variation of the form
$
q_t=(a_t,-a_t)
$
produces the state variation $\eta_t$ satisfying
\[
\dif\eta_t
=
b(a_t-a_t)\,\dif t+\sigma(a_t-a_t)\,\dif W_t
=0.
\]
Thus $\eta_T=0$, and the first variation of the cost is
$
\delta J(q)=\E[Y_T\eta_T]=0.
$
Hence
\begin{equation}
\label{eq:kernel}
\ker \mathcal L=\operatorname{span}\{(1,-1)\}
\end{equation}
is the null space of the linearized control-to-state map. This null space is the source of nonuniqueness.

\subsection{Explicit solution of the unregularized aggregate problem}

Define the market price of risk and the state-price density by
\begin{equation}
\label{eq:theta-rho}
\theta:=\frac{b}{\sigma},
\qquad
\rho_t:=\exp\left(-\theta W_t-\frac12\theta^2t\right).
\end{equation}
Then
\begin{equation}
\label{eq:rho-SDE}
\dif\rho_t=-\theta\rho_t\,\dif W_t,
\qquad \rho_0=1.
\end{equation}
In the unconstrained case, the Hamiltonian stationarity condition is
\begin{equation}
\label{eq:stationarity}
G_t=bY_t+\sigma Z_t=0.
\end{equation}
Consequently,
\begin{equation}
\label{eq:Z-vs-Y}
Z_t=-\theta Y_t.
\end{equation}
Substituting this into the adjoint BSDE gives
\[
\dif Y_t=-\theta Y_t\,\dif W_t.
\]
Therefore
\begin{equation}
\label{eq:Y-Y0}
Y_t=Y_0\rho_t.
\end{equation}
Using the terminal condition \eqref{eq:Y-terminal},
\[
Y_0\rho_T=\lambda(X_T-\E X_T)-1,
\]
or equivalently
\begin{equation}
\label{eq:XT-rho-unsolved}
X_T
=
\E X_T+\frac1\lambda+\frac{Y_0}{\lambda}\rho_T.
\end{equation}
Taking expectations and using $\E\rho_T=1$ yields
\[
\E X_T
=
\E X_T+\frac1\lambda+\frac{Y_0}{\lambda},
\]
and hence
\begin{equation}
\label{eq:Y0}
Y_0=-1.
\end{equation}
Thus
\begin{equation}
\label{eq:Y-Z-explicit}
\boxed{Y_t=-\rho_t,\qquad Z_t=\theta\rho_t.}
\end{equation}
It follows immediately that
\[
bY_t+\sigma Z_t
=
-b\rho_t+\sigma\theta\rho_t
=0.
\]
Let
$
m:=\E X_T.
$
Equation \eqref{eq:XT-rho-unsolved} becomes
\begin{equation}
\label{eq:XT-rho}
X_T=m+\frac1\lambda-\frac1\lambda\rho_T.
\end{equation}
The state-price density removes the risk-premium drift of wealth:
\[
\begin{aligned}
\dif(\rho_tX_t)
&=\rho_t\dif X_t+X_t\dif\rho_t+\dif[\rho,X]_t\\
&=b\rho_tv_t\,\dif t+\sigma\rho_tv_t\,\dif W_t
-\theta\rho_tX_t\,\dif W_t
-\theta\sigma\rho_tv_t\,\dif t\\
&=\sigma\rho_tv_t\,\dif W_t-\theta\rho_tX_t\,\dif W_t,
\end{aligned}
\]
because $\theta\sigma=b$. Hence $\rho_tX_t$ is a martingale and
\begin{equation}
\label{eq:budget}
\E[\rho_TX_T]=x_0.
\end{equation}
Moreover,
\begin{equation}
\label{eq:rho-moments}
\E\rho_T=1,
\qquad
\E\rho_T^2=e^{\theta^2T}.
\end{equation}
Substituting \eqref{eq:XT-rho} into \eqref{eq:budget} gives
$
m+\frac1\lambda-\frac1\lambda e^{\theta^2T}=x_0,
$
so
\begin{equation}
\label{eq:m-star}
\boxed{
m^*=x_0+\frac{e^{\theta^2T}-1}{\lambda}.
}
\end{equation}
The optimal terminal wealth is therefore
\begin{equation}
\label{eq:XT-star}
\boxed{
X_T^*
=
x_0+\frac{e^{\theta^2T}}{\lambda}
-\frac1\lambda\rho_T.
}
\end{equation}
At time $t$, the replicating wealth is
\begin{equation}
\label{eq:Xt-star}
X_t^*
=
m^*+\frac1\lambda
-\frac1\lambda e^{\theta^2(T-t)}\rho_t.
\end{equation}
Indeed,
\[
\E\left[\frac{\rho_T}{\rho_t}\,\middle|\mathcal F_t\right]=1,
\qquad
\E\left[\frac{\rho_T^2}{\rho_t}\,\middle|\mathcal F_t\right]
=
\rho_t e^{\theta^2(T-t)}.
\]
Applying It\^o's formula to \eqref{eq:Xt-star} gives
\begin{equation}
\label{eq:dX-star}
\dif X_t^*
=
\frac{\theta^2}{\lambda}e^{\theta^2(T-t)}\rho_t\,\dif t
+
\frac{\theta}{\lambda}e^{\theta^2(T-t)}\rho_t\,\dif W_t.
\end{equation}
Comparing \eqref{eq:dX-star} with
\[
\dif X_t=bv_t^*\,\dif t+\sigma v_t^*\,\dif W_t
\]
gives the optimal aggregate position
\begin{equation}
\label{eq:v-star}
\boxed{
v_t^*
=
\frac{b}{\lambda\sigma^2}
e^{\theta^2(T-t)}\rho_t.
}
\end{equation}

\begin{remark}
In the zero-interest-rate setting $r=0$, the threshold in Theorem~6.3 in \cite{hz05} reduces
to $z\ge x_0$, since $x_0e^{\int_0^T r(s)\,ds}=x_0$. Taking expectations in
\eqref{eq:XT-star} and using $\mathbb E[\rho_T]=1$ together with
$\mathbb E[\rho_T^2]=e^{\theta^2T}$ gives
\begin{equation}
\mathbb E[X_T^*]
= x_0+\frac{e^{\theta^2T}-1}{\lambda}.
\end{equation}
As $e^{\theta^2T}-1>0$, the efficient-frontier condition
$\mathbb E[X_T^*]\ge x_0$ holds exactly when $\lambda>0$. Equivalently,
\begin{equation}
\lambda=\frac{e^{\theta^2T}-1}{z-x_0},
\end{equation}
so the efficient branch $z>x_0$ corresponds to $\lambda>0$, the all-bank
endpoint $z=x_0$ to $\lambda\to\infty$, and $\lambda<0$ to the dominated,
inefficient branch $z<x_0$. Hence the terminal wealth \eqref{eq:XT-star} is
consistent with Theorem~6.3 precisely on the branch $\lambda>0$.
\end{remark}

\subsection{The minimum-energy split}

Every original optimal control has the form
\begin{equation}
\label{eq:original-solution-set}
u_{1,t}=\frac{v_t^*+w_t}{2},
\qquad
u_{2,t}=\frac{v_t^*-w_t}{2},
\end{equation}
where $w_t$ is an arbitrary square-integrable adapted process. By \eqref{eq:energy-vw},
\begin{equation}
\label{eq:energy-on-solution-set}
\|u\|_{\Htwo}^2
=
\frac12\E\int_0^T(v_t^*)^2\,\dif t
+
\frac12\E\int_0^Tw_t^2\,\dif t.
\end{equation}
The first term is fixed by the state equation. Therefore the energy is minimized uniquely by
\begin{equation}
\label{eq:w-dagger}
w_t^\dagger=0.
\end{equation}
Consequently, the minimum-norm element of the original solution set is
\begin{equation}
\label{eq:u-dagger}
\boxed{
u_{1,t}^\dagger=u_{2,t}^\dagger
=
\frac{v_t^*}{2}
=
\frac{b}{2\lambda\sigma^2}
e^{\theta^2(T-t)}\rho_t.
}
\end{equation}
This also has a Hilbert-space interpretation. The control-to-state map depends only on the linear functional $u_1+u_2$. The minimum-norm right inverse of this linear functional divides the aggregate position equally between the two redundant coordinates.


\subsection{The \texorpdfstring{$L^2$}{L2}-regularized problem}

For $\eps>0$, consider
\begin{equation}
\label{eq:regularized-objective}
J_\eps(u)
=
J(u)
+
\frac{\eps}{2}
\E\int_0^T(u_{1,t}^2+u_{2,t}^2)\,\dif t.
\end{equation}
Using \eqref{eq:energy-vw},
\begin{equation}
\label{eq:regularized-vw}
J_\eps(u)
=
J(v)
+
\frac{\eps}{4}
\E\int_0^T(v_t^2+w_t^2)\,\dif t.
\end{equation}
The regularized Hamiltonian is
\begin{equation}
\label{eq:regularized-H}
H_\eps
=
G_t(u_1+u_2)
+
\frac{\eps}{2}(u_1^2+u_2^2).
\end{equation}
In the unconstrained case, the first-order conditions are
\begin{equation}
\label{eq:foc-individual}
G_t^\eps+\eps u_{1,t}^\eps=0,
\qquad
G_t^\eps+\eps u_{2,t}^\eps=0.
\end{equation}
Thus
\begin{equation}
\label{eq:equal-regularized}
\boxed{u_{1,t}^\eps=u_{2,t}^\eps=-\frac{G_t^\eps}{\eps}.}
\end{equation}
Equivalently,
\begin{equation}
\label{eq:w-zero}
\boxed{w_t^\eps=0\quad\text{for every }\eps>0.}
\end{equation}

With a symmetric closed convex constraint set $K_t$, the variational inequality gives the projection response
\begin{equation}
\label{eq:projection-equal}
u_{i,t}^\eps
=
\Proj_{K_{i,t}}
\left[-\eps^{-1}G_t^\eps\right],
\qquad i=1,2.
\end{equation}
Because the two projected arguments and the two coordinate constraint sets are identical, the regularized controls remain equal.

The regularization does not change the terminal Lions derivative; it changes the Hamiltonian and the control response. As $\eps\downarrow0$, the aggregate regularized control converges to the unregularized aggregate optimizer $v^*$, while the long--short component is fixed at zero. The Tikhonov regularization theorem therefore gives
\begin{equation}
\label{eq:Tikhonov-convergence}
u^\eps\longrightarrow u^\dagger
\quad\text{strongly in }\Htwo(0,T;\R^2).
\end{equation}

\subsection{Finite-\texorpdfstring{$\eps$}{epsilon} feedback solution}

It is useful to have an explicit finite-$\eps$ representation for the numerical experiment. For a fixed target mean $\rhobar$, consider the standard quadratic control representation
\begin{equation}
\label{eq:fixed-target}
J_{\eps,\rhobar}(v)
=
\E\left[
\frac{\lambda}{2}(X_T-\rhobar)^2-X_T
+
\frac{\eps}{4}\int_0^Tv_t^2\,\dif t
\right].
\end{equation}
At the self-consistent solution $\rhobar=\E X_T$, the terminal term equals $\frac{\lambda}{2}\Var(X_T)-\E X_T$. The following lemma indicates that the control set can be $\R^{2}_+$.
\begin{lemma}
\label{lem:measure-hjb-redundant-assets}
Let \(W\) be a one-dimensional Brownian motion and consider the controlled wealth process
\begin{equation}
\label{eq:wealth-lemma}
dX_t
=
b\bigl(u_{1,t}+u_{2,t}\bigr)\,dt
+
\sigma\bigl(u_{1,t}+u_{2,t}\bigr)\,dW_t,
\qquad X_0=x_0,
\end{equation}
where \(b,\sigma\in\mathbb R_{+}\), \(\sigma\neq 0\), and \(u_i\in\mathbb H^2(0,T;\mathbb R)\). For \(\varepsilon>0\) and \(\lambda>0\), consider the regularized mean--variance functional
\begin{equation}
\label{eq:objective-lemma}
J_\varepsilon(u)
=
\frac{\lambda}{2}\Var(X_T)
-\mathbb E[X_T]
+
\frac{\varepsilon}{2}
\mathbb E\int_0^T
\bigl(u_{1,t}^2+u_{2,t}^2\bigr)\,dt.
\end{equation}
Let \(\mu_t=\mathcal L(X_t)\), \(m(\mu)=\int_{\mathbb R}x\,\mu(dx)\), and define the terminal cost functional
\begin{equation}
\label{eq:terminal-functional}
\Phi(\mu)
=
\frac{\lambda}{2}
\int_{\mathbb R}\bigl(x-m(\mu)\bigr)^2\,\mu(dx)
-m(\mu).
\end{equation}
Then the measure-valued HJB equation
\begin{equation}
\label{eq:measure-hjb-lemma}
\begin{aligned}
&\partial_t\mathcal V^\varepsilon(t,\mu)
+
\int_{\mathbb R}
\inf_{u\in\mathbb R^2}
\left\{
\frac{\varepsilon}{2}|u|^2
+b(\mathbf 1^\top u)\,
\partial_\mu\mathcal V^\varepsilon(t,\mu)(x)\right.\\
&\hspace{3.8cm}\left.
+\frac{\sigma^2}{2}(\mathbf 1^\top u)^2
\partial_x\partial_\mu\mathcal V^\varepsilon(t,\mu)(x)
\right\}\mu(dx)
=0,\\
&\mathcal V^\varepsilon(T,\mu)=\Phi(\mu),
\qquad \mathbf 1=(1,1)^\top,
\end{aligned}
\end{equation}
has the following properties.

\begin{enumerate}
\item[(i)] The value functional has the finite-dimensional quadratic form
\begin{equation}
\label{eq:value-ansatz}
\mathcal V^\varepsilon(t,\mu)
=
\frac12 A_t^\varepsilon\,\Var_\mu(X)
+B_t^\varepsilon m(\mu)+C_t^\varepsilon,
\end{equation}
where
\begin{equation}
\label{eq:riccati-AB}
A_T^\varepsilon=\lambda,\qquad B_T^\varepsilon=-1,\qquad C_T^\varepsilon=0,
\end{equation}
and
\begin{equation}
\label{eq:riccati-A}
\dot A_t^\varepsilon
=
\frac{b^2(A_t^\varepsilon)^2}
{\sigma^2A_t^\varepsilon+\varepsilon/2},
\qquad
B_t^\varepsilon=-1.
\end{equation}

\item[(ii)] The unique regularized optimal control satisfies
\begin{equation}
\label{eq:equal-split-lemma}
u_{1,t}^\varepsilon=u_{2,t}^\varepsilon
=\frac12v_t^\varepsilon,
\end{equation}
where
\begin{equation}
\label{eq:aggregate-feedback}
v_t^\varepsilon
=
-\frac{b\bigl[A_t^\varepsilon(X_t-m_t)-1\bigr]}
{\sigma^2A_t^\varepsilon+\varepsilon/2}
=
\alpha_t^\varepsilon(c_t^\varepsilon-X_t)>0,
\end{equation}
with
\begin{equation}
\label{eq:coefficients-lemma}
m_t=\mathbb E[X_t],\qquad
\alpha_t^\varepsilon
=
\frac{bA_t^\varepsilon}
{\sigma^2A_t^\varepsilon+\varepsilon/2},
\qquad
c_t^\varepsilon=m_t+\frac1{A_t^\varepsilon}.
\end{equation}
Moreover, \(c_t^\varepsilon\) is constant in time.

\item[(iii)] Define
\begin{equation}
\label{eq:K-definition}
K_T^\varepsilon
=
\int_0^T
\frac{b^2A_t^\varepsilon}
{\sigma^2A_t^\varepsilon+\varepsilon/2}\,dt.
\end{equation}
Then the self-consistent terminal mean is
\begin{equation}
\label{eq:mean-epsilon}
\bar m_\varepsilon:=\mathbb E[X_T^\varepsilon]
=
x_0+\frac{e^{K_T^\varepsilon}-1}{\lambda}.
\end{equation}

\item[(iv)] Set
\begin{equation}
\label{eq:theta-rho}
\theta=\frac{b}{\sigma},
\qquad
\rho_t
=
\exp\left(-\theta W_t-\frac12\theta^2t\right).
\end{equation}
As \(\varepsilon\downarrow0\),
\begin{equation}
\label{eq:A-limit}
A_t^\varepsilon
\longrightarrow
A_t^0=\lambda e^{-\theta^2(T-t)},
\qquad
K_T^\varepsilon\longrightarrow\theta^2T,
\end{equation}
and
\begin{equation}
\label{eq:mean-limit}
\bar m_\varepsilon
\longrightarrow
m^*
=
x_0+\frac{e^{\theta^2T}-1}{\lambda}.
\end{equation}
The regularized controls converge strongly in \(\mathbb H^2\) to the minimum-energy control
\begin{equation}
\label{eq:minimum-energy-control}
u_{1,t}^\dagger=u_{2,t}^\dagger
=
\frac{b}{2\lambda\sigma^2}
e^{\theta^2(T-t)}\rho_t.
\end{equation}
\end{enumerate}
\end{lemma}

\begin{remark}
\label{rem:relaxation-verification}

The key observation is the elementary fact that \emph{a global minimizer
on a larger set, once it belongs to a smaller subset, is a global
minimizer on that smaller subset}. Indeed, suppose the verification
step shows that \(\widetilde u\in C\). Then, for every
\(u\in C\subseteq H\),
$
J_\varepsilon(\widetilde u)\le J_\varepsilon(u),
$
because \(\widetilde u\) minimizes \(J_\varepsilon\) over the whole
space \(H\). Thus \(\widetilde u\) is feasible and no feasible point
attains a smaller value; it is precisely the constrained minimizer
$
\widetilde u
=u_\varepsilon
=\arg\min_{u\in C}J_\varepsilon(u).
$
Strong convexity of \(J_\varepsilon\) on \(C\), inherited from the
\(\varepsilon>0\) regularization, also makes this constrained minimizer
unique.

This device avoids the constrained KKT system (equivalently, a
projection): one first solves the easier unconstrained stationarity
equation and afterwards checks the single feasibility condition
\(\widetilde u\in C\) (in the long-only setting, \(\widetilde u\ge0\)).
It is only when the relaxed minimizer falls outside \(C\) that the
constraint has to be imposed actively, through a projection or a KKT
multiplier.
\end{remark}

\subsection{Vanishing-regularization limit}

As $\eps\downarrow0$, $D_t^\eps\to\sigma^2A_t^0$ and
\[
A_t^\eps\to A_t^0=\lambda e^{-\theta^2(T-t)},\qquad
K_T^\eps\to\theta^2T,\qquad
\alpha_t^\eps\to\frac{b}{\sigma^2},
\]
so
\[
\bar m_\eps\to m^*=x_0+\frac{e^{\theta^2T}-1}{\lambda},
\qquad
u_{1,t}^\dagger=u_{2,t}^\dagger
=\frac{b}{2\lambda\sigma^2}e^{\theta^2(T-t)}
\exp\!\left(-\theta W_t-\frac12\theta^2t\right),
\]
which remains strictly positive.


Since $\chi_t$ is log-normal, its second moment is available in closed form:
\begin{equation}
\E[\chi_t^2]=\chi_0^2\,e^{-2K_t+\Gamma_t},
\qquad
\E[v_t^2]=\alpha_t^2\chi_0^2\,e^{-2K_t+\Gamma_t},
\end{equation}
where \(K_t:=\int_0^t k_s\,ds, \Gamma_t:=\int_0^t\gamma_s^2\,ds.\)
Hence the expected total control energy and the (aggregate) regularization
cost are
\begin{equation}
\mathcal E_\eps:=\E\!\int_0^T v_t^2\,dt
=\int_0^T\alpha_t^2\chi_0^2\,e^{-2K_t+\Gamma_t}\,dt,
\qquad
\mathcal R_\eps=\frac{\eps}{4}\,\mathcal E_\eps .
\end{equation}
The pathwise cumulative energy is $\int_0^t v_s^2\,ds$ (shown for each
realized path in the bottom-left panel).

\subsection{Realized paths, terminal wealth and energy}

\begin{figure}[h]
\centering
\includegraphics[width=\textwidth]{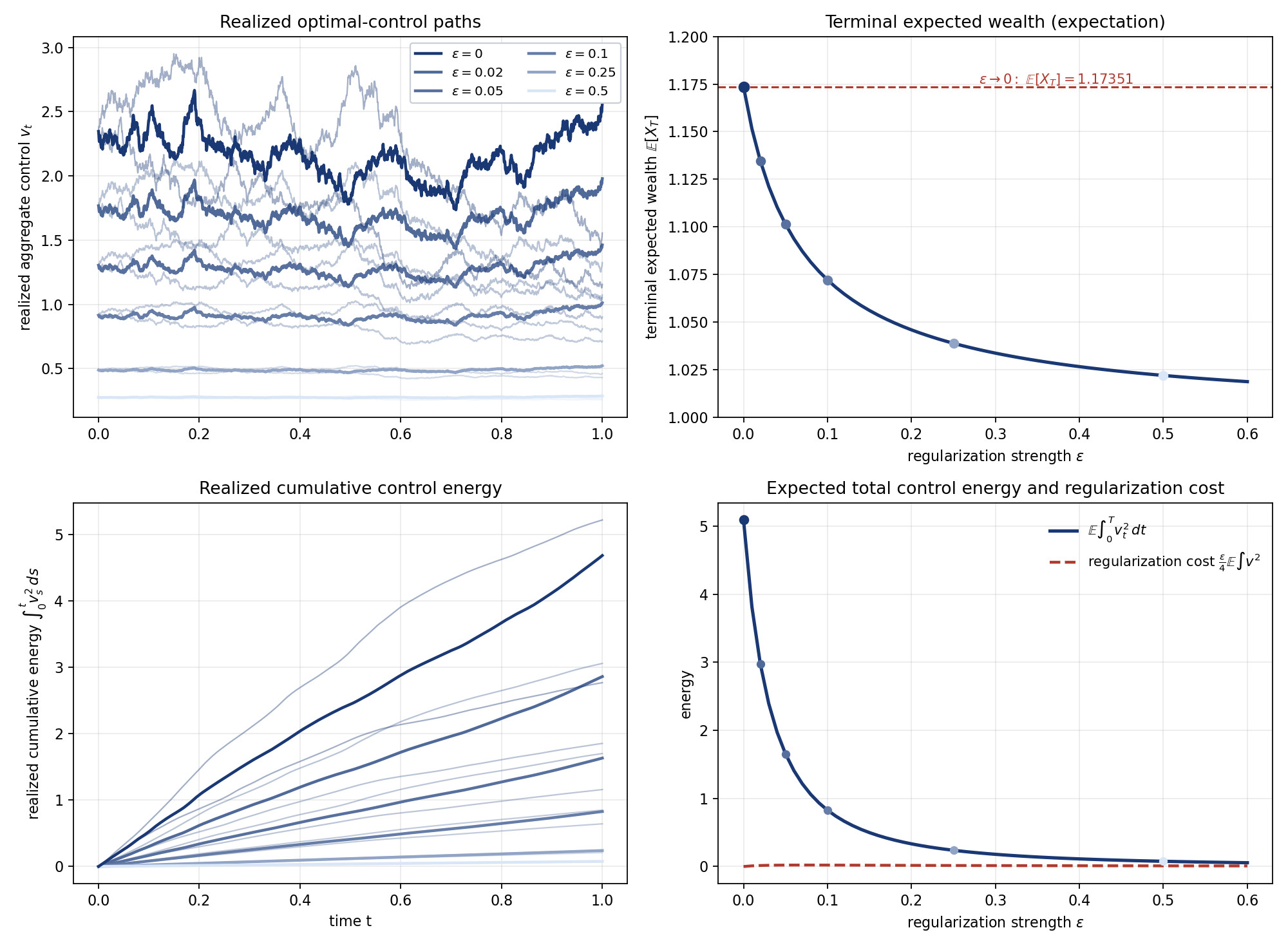}
\caption{Top-left: realized aggregate control $v_t$ (three paths per
$\eps$). Bottom-left: realized cumulative energy $\int_0^t v_s^2ds$.
Top-right: terminal expected wealth $\E[X_T]$ versus $\eps$.
Bottom-right: expected total energy $\E\int v_t^2dt$ and regularization
cost $\frac{\eps}{4}\E\int v^2$.}
\label{fig:paths}
\end{figure}

The top-left and bottom-left panels of Figure~\ref{fig:paths} show
\emph{realized} quantities along individual sample paths (the control
$v_t$ and the cumulative energy $\int_0^t v_s^2ds$); three Brownian paths are
shared across all $\eps$ so the comparison is pathwise. The right panels
report expectations: the terminal wealth $\bar m_\eps=\E[X_T]$ and the
expected total energy $\mathcal E_\eps$ with the regularization cost
$\mathcal R_\eps$.

\subsection{Numerical illustration}

Table~\ref{tab:coef} lists the numerical coefficients for $T=1$,
$b=0.08$, $\sigma=0.20$, $\lambda=1$, $x_0=1$ ($\theta=b/\sigma=0.4$).
All realized controls are strictly positive, so the long-only cone
constraint is slack throughout.

\begin{table}[h]
\centering
\small
\caption{Coefficients for each regularization strength. The last two
columns are the expected total control energy $\mathcal E_\eps$ and the
regularization cost $\mathcal R_\eps=\frac{\eps}{4}\mathcal E_\eps$.}
\label{tab:coef}
\begin{tabular}{@{}rrrrrrrrr@{}}
\toprule
$\eps$ & $c$ & $K_T$ & $A_0$ & $\alpha_0$ & $D_0$ & $\E[X_T]$ & $\mathcal E_\eps$ & $\mathcal R_\eps$\\
\midrule
0.00 & 2.17351 & 0.16000 & 0.85214 & 2.00000 & 0.03409 & 1.17351 & 5.0904 & 0.00000\\
0.02 & 2.13466 & 0.12633 & 0.88132 & 1.55804 & 0.04525 & 1.13466 & 2.9758 & 0.01488\\
0.05 & 2.10143 & 0.09661 & 0.90791 & 1.18456 & 0.06132 & 1.10143 & 1.6542 & 0.02068\\
0.10 & 2.07222 & 0.06973 & 0.93265 & 0.85460 & 0.08731 & 1.07222 & 0.8272 & 0.02068\\
0.25 & 2.03897 & 0.03823 & 0.96249 & 0.47095 & 0.16350 & 1.03897 & 0.2383 & 0.01490\\
0.50 & 2.02210 & 0.02186 & 0.97838 & 0.27070 & 0.28914 & 1.02210 & 0.0764 & 0.00956\\
\bottomrule
\end{tabular}
\end{table}

As $\eps\downarrow0$, $D_t\to\sigma^2A_t$, $A_t\to\lambda e^{-\theta^2(T-t)}$,
$K_T\to\theta^2T$, $\alpha_t\to b/\sigma^2$, and
$\bar m_\eps\to x_0+(e^{\theta^2T}-1)/\lambda=1.17351$. Control energy is
monotone decreasing in $\eps$; the regularization cost $\mathcal R_\eps$
starts at $0$, peaks for intermediate $\eps$ (around $\eps\approx0.05$--$0.10$),
and returns to $0$ as $\eps\to\infty$.

\begin{remark}
\label{rem:general-wealth}
Replace the wealth dynamics of the two redundant assets by
\[
dX_t=rX_t\,dt+\bigl(b_1u_{1,t}+b_2u_{2,t}\bigr)\,dt
+\bigl(\sigma_1u_{1,t}+\sigma_2u_{2,t}\bigr)\,dW_t,\qquad X_0=x_0,
\]
with possibly distinct coefficients $b_1\neq b_2$ and $\sigma_1\neq\sigma_2$.
A nonzero direction $(u_1,u_2)$ belongs to the null space of the
control-to-wealth map iff it annihilates both the drift and the diffusion:
\[
b_1u_1+b_2u_2=0,\qquad \sigma_1u_1+\sigma_2u_2=0.
\]
This homogeneous system admits a nonzero solution exactly when
\[
b_1\sigma_2-b_2\sigma_1=0
\quad\Longleftrightarrow\quad
\frac{b_1}{\sigma_1}=\frac{b_2}{\sigma_2},
\]
i.e. when the two assets have the same Sharpe ratio; the null-space direction is
then $n=(\sigma_2,-\sigma_1)$. The deterministic term $rX_t$ is removed by
discounting, $\widetilde X_t=e^{-rt}X_t$, and does not affect this conclusion.

\noindent\textbf{Case A: equal Sharpe ratios ($b_i=\theta\sigma_i$).}
Writing $b_i=\theta\sigma_i$ and defining the aggregate risk exposure
$v_t=\sigma_1u_{1,t}+\sigma_2u_{2,t}$, the wealth dynamics reduce to
\[
dX_t=\theta v_t\,dt+v_t\,dW_t,
\]
which is exactly the form of the redundant example after a change of aggregate
variable. Hence the aggregate position is unique while the allocation is not,
the regularized problem is uniquely solvable, and the strong $\mathbb H^2$
convergence together with the recovery of the optimal value as
$\varepsilon\downarrow0$ all persist. The minimum-norm allocation, i.e. the
orthogonal projection of the origin onto the line $\sigma_1u_1+\sigma_2u_2=v$,
is
\[
u_1=\frac{\sigma_1}{\sigma_1^2+\sigma_2^2}\,v,\qquad
u_2=\frac{\sigma_2}{\sigma_1^2+\sigma_2^2}\,v.
\]
It reduces to the equal split $u_1=u_2$ only when $\sigma_1=\sigma_2$;
correspondingly, the null space is
$\operatorname{span}\{(\sigma_2,-\sigma_1)\}$, which becomes
$\operatorname{span}\{(1,-1)\}$ only in that symmetric case.

\noindent\textbf{Case B: distinct Sharpe ratios.}
The direction $(\sigma_2,-\sigma_1)$ has zero diffusion,
$\sigma_1\sigma_2-\sigma_2\sigma_1=0$, but nonzero drift
$b_1\sigma_2-b_2\sigma_1\neq0$, and is therefore a riskless arbitrage. The
unregularized problem is unbounded below and admits no finite optimizer. For
$\varepsilon>0$ the quadratic penalty keeps the problem bounded, but as
$\varepsilon\downarrow0$ the controls escape to infinity and no finite limit
point exists. Under portfolio constraints that forbid the arbitrage (e.g.
long-only or box constraints), the optimum is instead a corner concentrated in
the asset with the larger Sharpe ratio, and the null-space selection mechanism
is no longer the relevant phenomenon.

\noindent The comparison with the original redundant-asset example is
summarized in Table~\ref{tab:general-vs-redundant}.
\end{remark}

\begin{table}[h]
\centering
\caption{Comparison of the redundant-asset example with the general wealth equation.}
\label{tab:general-vs-redundant}
\begin{tabular}{llll}
\toprule
Feature & Redundant example & Equal Sharpe ratio & Distinct Sharpe ratio\\
\midrule
Coefficients & $b_1=b_2,\ \sigma_1=\sigma_2$ & $b_i=\theta\sigma_i$ & general\\
Null space & $\operatorname{span}\{(1,-1)\}$ & $\operatorname{span}\{(\sigma_2,-\sigma_1)\}$ & $\{0\}$\\
No-arbitrage & yes & yes & no\\
Minimum-norm split & equal split & $\sigma_i$-weighted & not applicable\\
Finite optimizer (unconstrained) & yes, nonunique & yes, nonunique & no\\
\bottomrule
\end{tabular}
\end{table}

\begin{remark}
\label{rem:sharpe-scope}
The path-preserving long--short direction $w=u_1-u_2$ relies on the two assets sharing the
same Sharpe ratio $\theta=b/\sigma$. Along this direction both the drift and the diffusion
of wealth vanish, so it leaves the entire state path unchanged and belongs to the path
kernel $\mathcal K_{\mathrm{path}}$. If the two assets had distinct Sharpe ratios
$\theta_1\ne\theta_2$, the zero-diffusion direction $n=(\sigma_2,-\sigma_1)$ would carry a
nonzero drift $\sigma_1\sigma_2(\theta_1-\theta_2)$ and would no longer preserve the path;
the path kernel would then be trivial. This does not force the terminal kernel $\ker \mathcal L$ to
be trivial: controls may still produce the same terminal wealth without keeping the whole
path identical. The redundancy (path-kernel) statement is therefore specific to the
equal-Sharpe setup, whereas a nontrivial terminal kernel is a weaker property.
\end{remark}

\section{Conclusion}


This paper developed a FBSDVI for the constrained continuous-time mean-variance
portfolio problem with random coefficients. The formulation requires no
non-degeneracy of the volatility and accommodates unbounded control regions.
Working through the Hilbert-space operator $\mathcal L$---its null space,
orthogonal complement, and projection, rather than matrix Riccati equations,
the analysis is independent of the number of accounts. When $\mathcal L$ admits
a nontrivial null space, the optimal portfolio is not unique; vanishing
Tikhonov regularization selects the minimum-norm, minimum-energy solution.
Two redundant risky assets illustrate the mechanism: wealth depends only on
the aggregate position, fixed by the adjoint BSDE and state-price density,
while regularization drives the long--short direction to zero and yields the
equal-split portfolio. For finite $\varepsilon$, an infinite-dimensional HJB
equation reduces to finite-dimensional Riccati equations, and a numerical
example confirms the convergence of controls, terminal wealth, and energy as
$\varepsilon\downarrow0$.

\section{Appendix}
\label{app}

\begin{lemma}
\label{firor}
Let $K\subset\mathbb{R}^{m}$ be a nonempty convex set, and let
$J:U\to\mathbb{R}$ be differentiable on an open set $U$ containing $K$.
\begin{enumerate}
\item[(i)] If $J$ is convex on $K$, then a point $\bar u\in K$ is a global
minimizer of $J$ over $K$, i.e.
\begin{equation}
J(\bar u)\le J(v),\qquad \forall v\in K,
\end{equation}
if and only if
\begin{equation}\label{eq:VI}
\big\langle \nabla J(\bar u),\,v-\bar u\big\rangle\ge 0,
\qquad \forall v\in K.
\end{equation}
\item[(ii)] If, in addition, $J$ is strictly convex, then the global minimizer
$\bar u$ is unique whenever it exists.
\end{enumerate}
\end{lemma}

\begin{proof}
\emph{(Necessity.)} Suppose that $\bar u\in K$ is a global minimizer of $J$
over $K$. Fix any $v\in K$. Since $K$ is convex, for every $\lambda\in(0,1]$,
\begin{equation*}
w_\lambda:=\bar u+\lambda(v-\bar u)=(1-\lambda)\bar u+\lambda v\in K.
\end{equation*}
By the minimality of $\bar u$, $J(\bar u)\le J(w_\lambda)$, hence
\begin{equation*}
\frac{J\big(\bar u+\lambda(v-\bar u)\big)-J(\bar u)}{\lambda}\ge 0.
\end{equation*}
Letting $\lambda\downarrow 0$ and using the differentiability of $J$ at
$\bar u$, the left-hand side converges to the directional derivative
$\big\langle\nabla J(\bar u),v-\bar u\big\rangle$. Passing to the limit
preserves the non-strict inequality, so
\begin{equation*}
\big\langle\nabla J(\bar u),v-\bar u\big\rangle\ge 0.
\end{equation*}
Since $v\in K$ was arbitrary, \eqref{eq:VI} follows. Notice that this
implication only uses differentiability of $J$ and convexity of $K$;
convexity of $J$ is not required.

\emph{(Sufficiency.)} Assume \eqref{eq:VI}. We first recall the first-order
convexity inequality: if $J$ is differentiable and convex, then for all
$x,y\in K$,
\begin{equation}\label{eq:grad-lower}
J(y)\ge J(x)+\big\langle\nabla J(x),y-x\big\rangle.
\end{equation}
For completeness, \eqref{eq:grad-lower} follows from the definition of
convexity: for every $\lambda\in(0,1]$,
\begin{equation*}
J\big(x+\lambda(y-x)\big)
=J\big((1-\lambda)x+\lambda y\big)
\le (1-\lambda)J(x)+\lambda J(y),
\end{equation*}
which, after rearranging and letting $\lambda\downarrow0$, gives
\eqref{eq:grad-lower}.

Now take any $v\in K$ and apply \eqref{eq:grad-lower} with $x=\bar u$,
$y=v$:
\begin{equation*}
J(v)\ge J(\bar u)+\big\langle\nabla J(\bar u),v-\bar u\big\rangle.
\end{equation*}
By \eqref{eq:VI}, the inner-product term is nonnegative, and therefore
$J(v)\ge J(\bar u)$ for every $v\in K$. Thus $\bar u$ is a global minimizer
of $J$ over $K$.

\emph{(Uniqueness under strict convexity.)} Suppose $J$ is strictly convex
and $\bar u_1,\bar u_2\in K$, $\bar u_1\neq\bar u_2$, are both global
minimizers with common minimal value $m$. Then for $\lambda\in(0,1)$,
convexity of $K$ gives $(1-\lambda)\bar u_1+\lambda\bar u_2\in K$, and
strict convexity yields
\begin{equation*}
J\big((1-\lambda)\bar u_1+\lambda\bar u_2\big)
<(1-\lambda)J(\bar u_1)+\lambda J(\bar u_2)=m,
\end{equation*}
contradicting the minimality of $m$. Hence the minimizer is unique.
This completes the proof.
\end{proof}

\subsection{Proof of Lemma \ref{lem:measure-hjb-redundant-assets}}
\begin{proof}
We divide the proof into five steps.

\paragraph{Step 1: Lions derivatives of the terminal cost.}
Write $
m=m(\mu), M_2(\mu)=\int_{\mathbb R}x^2\,\mu(dx).
$
Then
\[
\Phi(\mu)=\frac{\lambda}{2}\bigl(M_2(\mu)-m^2\bigr)-m.
\]
The Lions derivative of the second moment is \(2x\), while the Lions derivative of \(m^2\) is \(2m\). Hence
\begin{equation}
\label{eq:terminal-lions}
\partial_\mu\Phi(\mu)(x)
=
\lambda(x-m)-1,
\qquad
\partial_x\partial_\mu\Phi(\mu)(x)=\lambda.
\end{equation}
This is consistent with the adjoint terminal condition
$
Y_T=\lambda(X_T-\mathbb E[X_T])-1.
$

\paragraph{Step 2: Finite-dimensional reduction.}
Because the dynamics and costs are linear--quadratic and the terminal cost depends only on the mean and variance, seek the value functional in the form
\[
\mathcal V^\varepsilon(t,\mu)
=
\frac12A_t\Var_\mu(X)+B_tm+C_t.
\]
The Lions derivatives are
\begin{equation}
\label{eq:lions-derivatives-ansatz}
p(t,\mu,x)
:=
\partial_\mu\mathcal V^\varepsilon(t,\mu)(x)
=
A_t(x-m)+B_t,
\end{equation}
and
\begin{equation}
\label{eq:second-lions-derivative}
q(t,\mu,x)
:=
\partial_x\partial_\mu\mathcal V^\varepsilon(t,\mu)(x)
=
A_t.
\end{equation}
The terminal condition gives $
A_T=\lambda, B_T=-1, C_T=0.
$

\paragraph{Step 3: Pointwise minimization and equal splitting.}
For a representative state \(x\sim\mu\), the inner minimization in \eqref{eq:measure-hjb-lemma} is
\begin{equation}
\label{eq:inner-problem}
\inf_{u_1,u_2\in\mathbb R}
\left\{
\frac{\varepsilon}{2}(u_1^2+u_2^2)
+bp(u_1+u_2)
+\frac{\sigma^2}{2}A_t(u_1+u_2)^2
\right\}.
\end{equation}
The first-order conditions are
\begin{equation}
\label{eq:foc-two-assets}
\varepsilon u_1+bp+\sigma^2A_t(u_1+u_2)=0,
\end{equation}
and
\begin{equation}
\label{eq:foc-two-assets-2}
\varepsilon u_2+bp+\sigma^2A_t(u_1+u_2)=0.
\end{equation}
Subtracting these equations yields $
\varepsilon(u_1-u_2)=0. $
Since \(\varepsilon>0\), $
u_1=u_2.
$
Thus the long--short component is identically zero:
\[
w_t^\varepsilon=u_{1,t}^\varepsilon-u_{2,t}^\varepsilon=0.
\]
Set \(v=u_1+u_2\). Since \(u_1=u_2=v/2\), the running penalty becomes $
\frac{\varepsilon}{2}(u_1^2+u_2^2)
=
\frac{\varepsilon}{4}v^2.
$
Writing $
D_t:=\sigma^2A_t+\frac{\varepsilon}{2},
$
the inner problem reduces to $
\label{eq:aggregate-inner}
\inf_{v\in\mathbb R}
\left\{
bpv+\frac12D_tv^2
\right\}.
$
The unique minimizer is
\begin{equation}
\label{eq:v-minimizer}
v=-\frac{bp}{D_t}.
\end{equation}
Using \eqref{eq:lions-derivatives-ansatz} and \(B_t=-1\), we obtain
\[
v
=
-\frac{b[A_t(x-m)+B_t]}{D_t}
=
-\frac{b[A_t(x-m)-1]}{D_t}.
\]
Along the controlled trajectory \(x=X_t\), this is
\begin{equation}
\label{eq:v-feedback-proof}
v_t^\varepsilon
=
-\frac{b[A_t^\varepsilon(X_t-m_t)-1]}{D_t}.
\end{equation}
Define $
\alpha_t=\frac{bA_t}{D_t}.
$
Since
\[
A_t(X_t-m_t)-1
=
A_t\left(X_t-m_t-\frac1{A_t}\right),
\]
we have
\begin{equation}
\label{eq:v-chi}
v_t^\varepsilon
=
\alpha_t(c_t-X_t),
\qquad
c_t:=m_t+\frac1{A_t}.
\end{equation}
The minimized value of \eqref{eq:aggregate-inner} is
$
\label{eq:minimized-value}
-\frac{b^2p^2}{2D_t}.
$
Moreover, $
\mathbb E_\mu[p]
=
B_t,
$
and
$
\mathbb E_\mu[p^2]
=
A_t^2\Var_\mu(X)+B_t^2.
$
Substitution into the measure-valued HJB and matching the coefficients of \(m\), \(\Var_\mu(X)\), and the constant term gives
\begin{equation}
\label{eq:riccati-system}
\dot A_t=\frac{b^2A_t^2}{D_t},
\qquad
\dot B_t=0,
\qquad
\dot C_t=\frac{b^2B_t^2}{2D_t}.
\end{equation}
Together with \(B_T=-1\), this gives \(B_t=-1\), proving \eqref{eq:riccati-A} and \eqref{eq:aggregate-feedback}.

\paragraph{Step 4: The mean is generated endogenously by the measure flow.}
Taking expectations in \eqref{eq:wealth-lemma} gives $
\dot m_t=b\mathbb E[v_t^\varepsilon].
$
From \eqref{eq:v-minimizer} and \(\mathbb E[p]=B_t=-1\),
$
\label{eq:mean-ode-proof}
\mathbb E[v_t^\varepsilon]
=
-\frac{bB_t}{D_t}
=
\frac{b}{D_t},
$
and hence
$
\label{eq:mean-dot}
\dot m_t=\frac{b^2}{D_t}.
$
On the other hand,
\begin{equation}
\label{eq:inverse-A-derivative}
\frac{d}{dt}\left(\frac1{A_t}\right)
=
-\frac{\dot A_t}{A_t^2}
=
-\frac{b^2}{D_t}.
\end{equation}
Therefore
$
\frac{d}{dt}\left(m_t+\frac1{A_t}\right)=0,
$
so \(c_t=m_t+1/A_t\) is constant.

Define
$
k_t:=\frac{b^2A_t}{D_t}.
$
The Riccati equation can be written as
$
\dot A_t=k_tA_t.
$
Consequently,
\begin{equation}
\label{eq:A-zero}
A_0
=
A_T\exp\left(-\int_0^Tk_t\,dt\right)
=
\lambda e^{-K_T},
\end{equation}
where
$
K_T=\int_0^Tk_t\,dt.
$
At the initial time,
\[
c=m_0+\frac1{A_0}=x_0+\frac1{A_0},
\]
because \(m_0=x_0\). At the terminal time,
\[
c=m_T+\frac1{A_T}
=
\bar m_\varepsilon+\frac1\lambda.
\]
Thus
\[
\bar m_\varepsilon+\frac1\lambda
=
x_0+\frac1{A_0}
=
x_0+\frac{e^{K_T}}{\lambda},
\]
which proves
$
\bar m_\varepsilon
=
x_0+\frac{e^{K_T}-1}{\lambda}.
$
No externally imposed mean-target equation is needed: the mean flow and the measure-valued HJB jointly generate the terminal mean. This is the measure-space version of the usual self-consistency condition \(\bar m=\mathbb E[X_T]\).

\paragraph{Step 5: The vanishing-regularization limit.}
When \(\varepsilon=0\),
$
D_t=\sigma^2A_t,
$
and the Riccati equation becomes
$
\dot A_t=\frac{b^2}{\sigma^2}A_t=\theta^2A_t.
$
Using \(A_T=\lambda\),
$
A_t^0=\lambda e^{-\theta^2(T-t)}.
$
Moreover,
$
k_t=\theta^2,
\alpha_t=\frac{b}{\sigma^2},
\gamma_t:=\sigma\alpha_t=\theta.
$
Define the tracking deviation
$
\chi_t^\varepsilon:=c_t^\varepsilon-X_t.
$
Then
$
v_t^\varepsilon=\alpha_t^\varepsilon\chi_t^\varepsilon,
$
and
\begin{equation}
\label{eq:chi-sde-epsilon}
d\chi_t^\varepsilon
=
-k_t^\varepsilon\chi_t^\varepsilon\,dt
-\gamma_t^\varepsilon\chi_t^\varepsilon\,dW_t.
\end{equation}
In the limit,
\begin{equation}
\label{eq:chi-sde-limit}
d\chi_t
=
-\theta^2\chi_t\,dt-\theta\chi_t\,dW_t,
\qquad
\chi_0=\frac1\lambda e^{\theta^2T}.
\end{equation}
The explicit solution is
\begin{equation}
\label{eq:chi-explicit}
\chi_t
=
\frac1\lambda e^{\theta^2(T-t)}
\exp\left(-\theta W_t-\frac12\theta^2t\right)
=
\frac1\lambda e^{\theta^2(T-t)}\rho_t.
\end{equation}
Therefore
\begin{equation}
\label{eq:v-star-proof}
v_t^*
=
\frac{b}{\sigma^2}\chi_t
=
\frac{b}{\lambda\sigma^2}
e^{\theta^2(T-t)}\rho_t.
\end{equation}
The coefficients \(A_t^\varepsilon,\alpha_t^\varepsilon,k_t^\varepsilon,\gamma_t^\varepsilon\) converge uniformly to their \(\varepsilon=0\) limits. Standard linear SDE estimates, using Gronwall's inequality and the It\^o isometry, then give
$
\label{eq:strong-convergence}
\mathbb E\int_0^T
|v_t^\varepsilon-v_t^*|^2\,dt
\longrightarrow0.
$
Since \(u_{1,t}^\varepsilon=u_{2,t}^\varepsilon=v_t^\varepsilon/2\) for every \(\varepsilon>0\), the same convergence holds for the individual controls.

Finally, writing
$
w_t=u_{1,t}-u_{2,t},
$
one has the algebraic identity
$
\label{eq:energy-identity-proof}
u_{1,t}^2+u_{2,t}^2
=
\frac12\bigl(v_t^2+w_t^2\bigr).
$
The state equation depends only on \(v_t\), so the original problem determines \(v_t^*\) but leaves \(w_t\) undetermined. Among all controls generating the same optimal wealth process, minimizing the \(\mathbb H^2\)-norm requires \(w_t=0\). Hence the strong limit is the unique minimum-energy control
\[
u_{1,t}^\dagger=u_{2,t}^\dagger
=
\frac{v_t^*}{2}
=
\frac{b}{2\lambda\sigma^2}
e^{\theta^2(T-t)}\rho_t.
\]

\paragraph{Step 6: The verification theorem.}
When $\lambda>0$, the terminal cost $
\Phi(\mu)=\frac{\lambda}{2}\int_{\mathbb R}(x-m(\mu))^2\,\mu(dx)-m(\mu)
$
does not satisfy the global bounded-derivative, since
$
\partial_\mu\Phi(\mu)(x)=\lambda(x-m(\mu))-1,
\partial_x\partial_\mu\Phi(\mu)(x)=\lambda.
$
The first Lions derivative grows linearly in $x$ and is therefore not globally
bounded. But the cited bounded-derivative verification theorem (see \cite{gpw23}) cannot be invoked as stated. Accordingly, we have removed the direct appeal to that theorem and replaced it
with a self-contained verification under explicit growth hypotheses. The
candidate value functional
\[
\mathcal V^\varepsilon(t,\mu)
=\tfrac12 A_t^\varepsilon\operatorname{Var}_\mu(X)
+B_t^\varepsilon m(\mu)+C_t^\varepsilon
\]
satisfies, with a constant $C$ independent of $t$ and $\mu$,
\[
|\mathcal V^\varepsilon(t,\mu)|
\le C\Bigl(1+\int_{\mathbb R}|x|^2\mu(dx)\Bigr),\qquad
|\partial_\mu\mathcal V^\varepsilon(t,\mu)(x)|\le C(1+|x|),\qquad
|\partial_x\partial_\mu\mathcal V^\varepsilon(t,\mu)(x)|\le C.
\]
For every admissible control $u\in\mathbb H^2(0,T;\mathbb R)$, the standard
energy estimate (It\^o isometry together with the BDG inequality and Gronwall)
gives $\mathbb E[\sup_{t\le T}|X_t|^2]<\infty$. In the lifted Lions It\^o
formula the only stochastic integral is
\[
M_t=\int_0^t\sigma\bigl[A_s^\varepsilon(\widetilde X_s-m_s)
+B_s^\varepsilon\bigr]\,dW_s,
\]
and
\[
\mathbb E\langle M\rangle_T
\le C\,\mathbb E\!\int_0^T(1+|\widetilde X_s|^2)\,ds<\infty.
\]
Thus $M$ is a true square-integrable martingale; equivalently, one may stop at
$\tau_n=\inf\{t:\int_0^t|\partial_x\mathcal V^\varepsilon\sigma|^2ds\ge n\}$
and pass to the limit by Fatou. Taking expectations in the Lions It\^o identity
and using the HJB inequality yields
\[
J_\varepsilon(u)\ge \mathcal V^\varepsilon(0,\mu_0)
\quad\text{for every admissible }u,
\]
with equality for the feedback
\[
v_t^\varepsilon
=-\frac{b[A_t^\varepsilon(X_t-m_t)-1]}
{\sigma^2A_t^\varepsilon+\varepsilon/2}
=\alpha_t^\varepsilon(c_t^\varepsilon-X_t).
\]
Hence the feedback is the unique optimal control, now established without any
bounded-derivative assumption. This completes the proof.
\end{proof}

\begin{remark}
Suppose the unconstrained control set \(\mathbb R^2\) is replaced by the symmetric box
$
K=[\ell,h]^2.
$
The unconstrained first-order condition is replaced by the Euclidean projection associated with the strictly convex pointwise minimization problem. The objective and the feasible set are both invariant under the exchange \(u_1\leftrightarrow u_2\), so the unique projected minimizer still satisfies
$
u_{1,t}^\varepsilon=u_{2,t}^\varepsilon.
$
When the bounds are inactive along the optimal trajectory, formulas \eqref{eq:riccati-A}--\eqref{eq:mean-epsilon} remain unchanged. If a bound is active, both coordinates are clipped identically, and the explicit interior Riccati feedback is replaced by the corresponding projected feedback.
\end{remark}


\subsection{Properties of projection}

We recall the following properties of projection $\mathbf{P}_{U}$ onto a
closed convex set $U$, see \cite{Brezis 2010}, Chapter 5.

\begin{theorem}
\label{projection theorem} For a nonempty closed convex set $U\subset
\mathbb{R}^{m}$, for every $x\in \mathbb{R}^{m}$, there exists a unique $%
x^{\ast }\in U$, such that
\begin{equation*}
|x-x^{\ast }|=\min_{y\in \Gamma }|x-y|=:dist(x,U).
\end{equation*}%
Moreover, $x^{\ast }$ is characterized by the property
\begin{equation}
x^{\ast }\in U,\quad \big<x^{\ast }-x,x^{\ast }-y\big>\leq 0\qquad \forall
y\in U.  \label{projection characterization}
\end{equation}%
The above element $x^{\ast }$ is called the projection of $x$ onto $U$ and
is denoted by $\mathbf{P}_{U}[x]$.
\end{theorem}

One can immediately obtain the following

\begin{proposition}
\label{pro1}Let $U\subset \mathbb{R}^{m}$ be a nonempty closed convex set,
then we have
\begin{equation}
\big|\mathbf{P}_{U}[x]-\mathbf{P}_{U}[y]\big|^{2}\leq \big<\mathbf{P}_{U}[x]-%
\mathbf{P}_{U}[y],x-y\big>.  \label{projection inequality}
\end{equation}
\end{proposition}

\begin{proposition}
\label{pro2}Let $U\subset \mathbb{R}^{m}$ be a nonempty closed convex set,
then the projection $\mathbf{P}_{U}$ does not increase the distance, i.e.
\begin{equation}
\big|\mathbf{P}_{U}[x]-\mathbf{P}_{U}[y]\big|\leq \big|x-y\big|.
\label{projection
lipschitz}
\end{equation}
\end{proposition}

\noindent \textbf{Competing Interests Statement.}
The authors declare that there are no known competing financial interests or
personal relationships that could have appeared to influence the work reported in this paper.

\end{document}